
\def\texversion{1}      
 
\ifnum\texversion=0
        \font\normal=Palatino at 10pt
        \font\ita=cmsl10 
        \font\bol=PalatinoB at 10pt
        \font\small=cmr6
\else
        \font\normal=cmr10 
        \font\ita=cmsl10
        \font\bol=cmbx10
        \font\small=cmr6
\fi
\normal

\overfullrule=0pt
\baselineskip 14pt		

\hsize 6truein
\vsize=8truein


\def\today{\ifcase\month\or
January\or February\or March\or April\or May\or June\or
July\or August\or September\or October\or November\or
December\fi
\space\number\day, \number\year}

\headline={\small (Final Version)\hfill\small Lattice diagrams and extended Pieri rules $\ess\ess\ess\ess\ess$ 
\hfill$\ess\ess\ess$\today$\ess\ess\ess\ess\ess\ess$ \folio } \footline={\hfil}


\def\multi#1{\vbox{\baselineskip=0pt\halign{\hfil$\scriptstyle\vphantom{(_)}##$\hfil\cr#1\crcr}}}

\def\heading#1{\noindent{#1}\par\nobreak}


\def \BA {{\bf A}}
\def \BB {{\bf B}}
\def \BH {{\bf H}}
\def \BK {{\bf K}}

\def \BM {{\bf M}}
\def \BQ {{\bf Q}}
\def \BR {{\bf R}}

\def \BX {{\bf X}}
\def \BY {{\bf Y}}
\def \BTM {{ \BM_\TDD}}
\def\BTMperp{{\BM_{\TDD}^{\;\;\perp}}}
\def\BMperp{{\BM_{\DD}^{\;\perp}}}
\def\BTK{{\bf\widetilde K}}
\def\BTA{{\bf\widetilde A}}


\def \CB {{\cal B}}
\def \CC {{\cal C}}
\def \CH {{\cal H}}
\def \CL {{\cal L}}
\def \CP {{\cal P}}
\def \CR {{\cal R}}
\def \CX {{\cal X}}
\def \CY {{\cal Y}}


\def \TK {{\tilde K}}
\def \TH {{\tilde H}}

\def\TDD{{\widetilde\Delta}}


\def \eee {\epsilon}

\def\aaa {\alpha}
\def\bbb {\beta}

\def\aa {\alpha}
\def\bb {\beta}

\def\DD {\Delta}

\def\om {\omega}

\def\la {\lambda}

\def\sig{\sigma}


\def \RA {{ \rightarrow }}

\def \DA {\downarrow}
\def \da {\downarrow}


\def \ses {\enskip = \enskip}
\def \sps {\enskip + \enskip}
\def \sms {\enskip -\enskip}
\def \scs {\ssp , \ssp}
\def \ess {\enskip}
\def \ssp {\hskip .25em}
\def \sssp {\hskip .15em}
\def \bigsp {\hskip .5truein}


\def \sas {\vskip .06truein}
\def\sa{{\vskip .125truein}}
\def\sap{{\vskip .25truein}}


\def\flip{\mathop{\bf f{}lip}\nolimits}
\def\sign{\mathop{\rm sign}\nolimits}
\def\weight{\mathop{\rm weight}\nolimits}
\def\CF{{\cal F}}
\def\ch{\mathop{\rm ch}\nolimits}
\def\Fch{\mathop{\CF\,\rm ch}\nolimits}
\def\gr{\mathop{\rm gr}\nolimits}
\def\nshadow{\mathop{\#\rm shadow}\nolimits}
\def\Pred{\mathop{{\CP}red}\nolimits}


\def\del {\partial}
\def \con {\subseteq} 
\def \part {\vdash}
\def \pmu {\rho_\mu}
\def \xon {x_1,\ldots ,x_n}
\def \yon {y_1, \ldots ,y_n}

\def\LL{\big\langle}
\def\RR {\big\rangle}
\def \LLL {\langle\langle}
\def \RRR {\rangle\rangle}


\newdimen\mylength
\newdimen\mythick
\mylength=.3pt		
\mythick=.4pt		
\def\boldthick{3\mythick}
\def\cornd{20}		

\newcount\drawcntA
\newdimen\drawdimX
\newdimen\drawdimY

\def\myput(#1,#2)#3{\raise#2\mylength\hbox to 0pt{\hskip #1\mylength #3\hss}\ignorespaces}

\def\multiputline(#1,#2)(#3,#4)#5#6{%
	\drawcntA=#5%
	\drawdimX=#1\mylength \drawdimY=#2\mylength%
	\loop\ifnum\drawcntA>0%
		\raise\drawdimY\hbox to 0pt{\hskip\drawdimX#6\hss}%
		\advance\drawdimX by #3\mylength%
		\advance\drawdimY by #4\mylength%
		\advance\drawcntA by -1%
	\repeat%
	\ignorespaces
}

\def\multiputgrid(#1,#2)(#3,#4)#5(#6,#7)#8#9{
	\multiputline(#1,#2)(#3,#4){#5}{{\multiputline(0,0)(#6,#7){#8}{#9}}}%
	\ignorespaces
}

\def\hseg(#1,#2;#3){\myput(#1,#2){\vrule width#3\mylength depth\mythick}\ignorespaces}

\def\vseg(#1,#2;#3){\myput(#1,#2){\vrule height#3\mylength width\mythick depth\mythick}\ignorespaces}

%
%

\def\adjustY{\kern-\mythick}
\def\rectangle(#1,#2)(#3,#4){\myput(#1,#2){%
	\drawdimX=\mythick 	\advance\drawdimX by #3\mylength
	\drawdimY=\mythick	\advance\drawdimY by #4\mylength
	\let\adjustH\adjustY \let\adjustV\relax
	\rectangleframe}\ignorespaces}
\def\rectangleframe{%
	\vbox to0pt{\hbox to0pt{\hss%
	\vbox to\drawdimY{%
		\adjustV%
		\hrule width \drawdimX height\mythick%
		\vss%
		\hbox to\drawdimX{%
			\vrule height\drawdimY width\mythick%
			\hss \vrule height\drawdimY width\mythick%
		}%
		\vss%
		\hrule width\drawdimX height\mythick%
	}%
	\adjustH}\vss}}

\def\corn(#1,#2){\rectangle(#1,#2)(\cornd,\cornd)}
\def\bcorn(#1,#2){\myput(#1,#2){%
	\drawdimX=\mythick%
	\advance\drawdimX by \cornd\mylength%
	\drawdimY=\drawdimX%
	\let\adjustV\relax%
	\edef\adjustH{\expandafter\kern-\the\mythick}%
	\mythick=\boldthick%
	\rectangleframe}\ignorespaces}
\def\brectangle(#1,#2)(#3,#4){\myput(#1,#2){%
	\drawdimX=\boldthick 	\advance\drawdimX by #3\mylength
	\drawdimY=\mythick	\advance\drawdimY by #4\mylength
	\let\adjustH\adjustY
	\let\adjustV\relax
	\mythick=\boldthick
	\rectangleframe}\ignorespaces}

\def\labUR(#1,#2;#3){\myput(#1,#2){\hskip 2pt$#3$}\ignorespaces}
\def\labLL(#1,#2;#3){\myput(#1,#2){\myput(-\cornd,-\cornd){%
	\vbox to0pt{\hbox to0pt{\hss $#3$\hskip 1pt}\vss}}}\ignorespaces}


\def\drawL(#1,#2,#3,#4;#5,#6){%
	\drawcntA=-#3 \advance\drawcntA#4%
	\vseg(#2,#3;\drawcntA)\corn(#2,#4)\labUR(#2,#4;#5)%
%
	\drawcntA=-#1 \advance\drawcntA#2%
	\hseg(#1,#4;\drawcntA)\corn(#1,#4)\labLL(#1,#4;#6)%
	\ignorespaces
}

\def\lhair(#1,#2)#3{\multiputline(#1,#2)(0,20){#3}{\hseg(-5,0;5)}}
\def\rhair(#1,#2)#3{\multiputline(#1,#2)(0,20){#3}{\hseg(0,0;5)}}
\def\dhair(#1,#2)#3{\multiputline(#1,#2)(20,0){#3}{\vseg(0,-5;5)}}
\def\uhair(#1,#2)#3{\multiputline(#1,#2)(20,0){#3}{\vseg(0,0;5)}}

\def\Brectangle(#1,#2)(#3,#4){\myput(#1,#2){%
	\drawdimX=\mythick 	\advance\drawdimX by #3\mylength
	\drawdimY=\mythick	\advance\drawdimY by #4\mylength
	\edef\adjustH{\expandafter\kern-\the\mythick}%
	\let\adjustV\relax%
	\mythick=\boldthick%
	\rectangleframe}\ignorespaces}


\def\filllinedensity{3}
\def\filllineskip{6.666666666}
\def\hrectangle(#1,#2)(#3,#4){%
	\drawcntA=#4%
	\multiply\drawcntA by\filllinedensity%
	\divide\drawcntA by\cornd%
	\advance\drawcntA by 1%
	\multiputline(#1,#2)(0,-\filllineskip){\the\drawcntA}{\hseg(-#3,0;#3)}%
	\ignorespaces
}
\def\vrectangle(#1,#2)(#3,#4){%
	\drawcntA=#3%
	\multiply\drawcntA by\filllinedensity%
	\divide\drawcntA by\cornd%
	\advance\drawcntA by 1%
	\multiputline(#1,#2)(-\filllineskip,0){\the\drawcntA}{\vseg(0,-#4;#4)}%
	\ignorespaces
}

\edef\savecatcodeat{\the\catcode`@}
\catcode`\@=11

\def\@IFNEXTCHAR#1#2#3{\let\@tempe #1\def\@tempa{#2}\def\@tempb{#3}\futurelet
    \@tempc\@IFNCH}
\def\@IFNCH{\ifx \@tempc \@sptoken \let\@tempd\@xifnch
      \else \ifx \@tempc \@tempe\let\@tempd\@tempa\else\let\@tempd\@tempb\fi
      \fi \@tempd}

\def\tb@ifSpecChars#1#2{#1}
\def\tb@ifNoSpecChars#1#2{#2}

\newcount\tb@rpN


\def\tableau{%
  \bgroup
%
  \@IFNEXTCHAR[{\tb@tableauC}{\tb@tableauC[]}}     

\def\tb@tableauC[#1]{\hbox\bgroup%
    \let\\=\cr
    \def\bl{\global\let\tbcellF\tb@cellNF}%
    \def\tf{\global\let\tbcellF\tb@cellH}
%
    \dimen2=\ht\strutbox \advance\dimen2 by\dp\strutbox%
    \ifx\baselinestretch\undefined\relax%
    \else%
       \dimen0=100sp \dimen0=\baselinestretch\dimen0%
       \dimen2=100\dimen2 \divide\dimen2 by\dimen0%
    \fi%
    \let\tpos\tb@vcenter
    \tb@initYoung
    \tb@options#1\eoo
    \let\arrow\tb@arrow%
    \dimen0=\Tscale\dimen2%
    \dimen1=\dimen0 \advance\dimen1 by \tb@fframe%
    \lineskip=0pt\baselineskip=0pt
%
    \def\tb@nothing{}%
    \def\endcellno{$\rss\egroup\bss\egroup}
    \def\endcell{\endcellno\kern-\dimen0}
    \def\begincell{\vbox to\dimen0\bgroup\vss\hbox to\dimen0\bgroup\hss$}%
    \let\overlay\tb@overlay%
    \let\fl\tb@fl%
    \let\fr\tb@fr%
    \let\lss\hss\let\rss\hss\let\tss\vss\let\bss\vss
    \def\mkcell##1{
        \let\tbcellF\tb@cellD
        \def\tb@cellarg{##1}
        \ifx\tb@cellarg\tb@nothing\let\tb@cellarg\tb@cellE\fi%
%
	        \begincell\tb@cellarg\endcellno
	        \tbcellF
    }%
    \let\savecellF\tbcellF
    \tb@tableauD%
}%

\let\tb@savetableauD\tableauD
{
\gdef\tableauD#1{%
  \tpos{\tabskip=0pt\halign{&\mkcell{##}\cr#1\crcr}}%
  \global\let\tbcellF\savecellF
  \egroup
  \egroup}
}
\let\tb@tableauD\tableauD
\let\tableauD\tb@savetableauD
\let\tb@savetableauD\undefined


\def\tb@options#1{\ifx#1\eoo\relax\else\tb@option#1\expandafter\tb@options\fi}

\def\tb@option#1{%
  \if#1t\let\tpos\tb@vtop\fi
  \if#1c\let\tpos\tb@vcenter\fi
  \if#1b\let\tpos\vbox\fi
  \if#1F\tb@initFerrers\fi
  \if#1Y\tb@initYoung\fi
  \if#1E\tb@initEmpty\fi
  \if#1s\tb@initSmall\fi
  \if#1m\tb@initMedium\fi
  \if#1l\tb@initLarge\fi
  \if#1p\tb@initPartition\fi
  \if#1a\tb@initArrow\fi
}

\def\tb@vcenter#1{\ifmmode\vcenter{#1}\else$\vcenter{#1}$\fi}

\def\tb@vtop#1{\hbox{\raise\ht\strutbox\hbox{\lower\dimen0\vtop{#1}}}}

\def\tb@initPartition{\def\Tscale{.3}}
\def\tb@initSmall{\def\Tscale{1}}
\def\tb@initMedium{\def\Tscale{2}}
\def\tb@initLarge{\def\Tscale{3}}

\def\tb@initArrow{\dimen2=1.25em}

\def\tb@initYoung{%
  \def\tb@cellE{}
  \let\tb@cellD\tb@cellN
  \def\sk{\global\let\tbcellF\tb@cellNF}}
\def\tb@initFerrers{%
  \def\tb@cellE{\bullet}
  \let\tb@cellD\tb@cellNF
  \def\sk{\bullet}}
\def\tb@initEmpty{%
  \def\tb@cellE{}
  \let\tb@cellD\tb@cellNF
  \def\sk{\global\let\tbcellF\tb@cellNF}}

\tb@initMedium

\def\tb@sframe#1{%
  \vbox to0pt{
    \vss
    \hbox to0pt{%
      \hss
      \vbox to\dimen1{
        \hrule depth #1 height0pt
        \vss
        \hbox to\dimen1{
          \vrule width #1 height\dimen1
          \hss
          \vrule width #1
          }%
        \vss
        \hrule height #1 depth 0in
        }%
      \kern-\tb@hframe
      }%
    \kern-\tb@hframe}}

\def\tb@hframe{.2pt}\def\tb@fframe{.4pt}\def\tb@bframe{2pt}
\def\tb@cellH{\tb@sframe{\tb@bframe}}       
\def\tb@cellNF{}                            
\def\tb@cellN{\tb@sframe{\tb@fframe}}       
\let\tbcellF\tb@cellN                       

\def\tb@Fsframe{%
  \vbox to0pt{
    \vss
    \hbox to0pt{%
      \hss
      \vbox to\dimen1{
        \fr@iftop{\hrule depth \fr@width height0pt}{\vskip \fr@width}
        \vss
        \hbox to\dimen1{
	  \fr@ifleft{\vrule width \fr@width height\dimen1}{\hskip \fr@width}
          \hss
          \fr@ifright{\vrule width \fr@width height\dimen1}{\hskip \fr@width}
          }%
        \vss
        \fr@ifbottom{\hrule height \fr@width depth 0in}{\vskip\fr@width}
        }%
      \kern-\tb@hframe
      }%
    \kern-\tb@hframe}}

\def\tb@fr{\@IFNEXTCHAR[{\tb@fra}{\global\let\tbcellF\tb@cellN}}
\def\tb@fra[#1]{%
	\global\let\fr@iftop\tb@IFNO
	\global\let\fr@ifbottom\tb@IFNO%
	\global\let\fr@ifleft\tb@IFNO%
	\global\let\fr@ifright\tb@IFNO%
	\global\let\fr@width\tb@fframe%
	\global\let\tbcellF\tb@Fsframe%
	\froptions#1\eoo
}
\def\froptions#1{\ifx#1\eoo\relax\else\froption#1\expandafter\froptions\fi}
\def\froption#1{
	\if#1t\global\let\fr@iftop\tb@IFYES\fi
	\if#1b\global\let\fr@ifbottom\tb@IFYES\fi
	\if#1l\global\let\fr@ifleft\tb@IFYES\fi
	\if#1r\global\let\fr@ifright\tb@IFYES\fi
	\if#1w\global\let\fr@width\tb@bframe\fi
}
\def\tb@IFYES#1#2{#1}
\def\tb@IFNO#1#2{#2}

\def\tb@rpad{1pt}
\def\tb@lpad{1pt}
\def\tb@tpad{1.8pt}
\def\tb@bpad{1.8pt}

\def\tb@overlay{\endcell\@IFNEXTCHAR[{\tb@overlaya}{\begincell}}
\def\tb@overlaya[#1]{\vbox to\dimen0\bgroup%
  \tb@overlayoptions#1\eoo%
  \tss\hbox to\dimen0\bgroup\lss$}
\def\tb@overlayoptions#1{\ifx#1\eoo\relax\else\tb@overlayoption#1\expandafter\tb@overlayoptions\fi}

\def\tb@overlayoption#1{
  \if#1t\def\tss{\vskip\tb@tpad}\let\bss\vss\fi
  \if#1c\let\tss\vss\let\bss\vss\fi
  \if#1b\def\bss{\vskip\tb@bpad}\let\tss\vss\fi
  \if#1l\def\lss{\hskip\tb@lpad}\let\rss\hss\fi
  \if#1m\let\lss\hss\let\rss\hss\fi
  \if#1r\def\rss{\hskip\tb@rpad}\let\lss\hss\fi
}

\def\tb@fl{\endcell\begincell\vrule depth 0pt width \dimen0 height \dimen0 \endcell\begincell}



\def\tbgobble#1{}
\def\Pscale{1}

\def\skewptn{%
  \@IFNEXTCHAR[{\tb@ptnC}{\tb@ptnC[]}}     

\def\tb@ptnC[#1](#2){%
	{%
    \let\Tscale\Pscale
    \let\\=\cr
   \def\tb@initYoung{%
	\def\tb@cell{\hskip\dimen0\tb@cellN}%
	\def\tb@kernA{\kern.5\dimen0}%
	\def\tb@kernB{\kern-.5\dimen0}%
   }%
   \def\tb@initFerrers{%
	\def\tb@cell{\hbox to\dimen0{\hss$\bullet$\hss}}%
	\def\tb@kernA{}%
	\def\tb@kernB{}%
   }%
%
    \dimen2=\ht\strutbox \advance\dimen2 by\dp\strutbox%
    \ifx\baselinestretch\undefined\relax%
    \else%
       \dimen0=100sp \dimen0=\baselinestretch\dimen0%
       \dimen2=100\dimen2 \divide\dimen2 by\dimen0%
    \fi%
    \let\tpos\tb@vcenter
    \tb@initYoung
    \tb@options#1\eoo
    \dimen0=\Tscale\dimen2%
    \dimen1=\dimen0 \advance\dimen1 by \tb@fframe%
    \lineskip=0pt\baselineskip=0pt
    \tpos{\skewptnDnewline#2|)}%
	}%
}%

\def\skewptnDnewline#1|{\vbox to\dimen0\bgroup\vss\tb@kernA\hbox\bgroup\skewptnEon#1,|}
\def\skewptnDendline|{\egroup\tb@kernB\vss\egroup\@IFNEXTCHAR{)}{\tbgobble}{\skewptnDnewline}}
\def\skewptnEon#1,{%
	\tb@rpN=#1%
	\ifnum#1>0
	        \loop%
		\tb@cell%
	        \ifnum\tb@rpN>1\advance\tb@rpN by-1%
        	\repeat%
	\fi%
	\@IFNEXTCHAR{|}{\skewptnDendline}{\skewptnEoff}}
\def\skewptnEoff#1,{\hskip #1\dimen0%
	\@IFNEXTCHAR{|}{\skewptnDendline}{\skewptnEon}}


\catcode`\@=\savecatcodeat
\let\savecatcodeat\undefined

\sap
\centerline {\bol  Lattice Diagram Polynomials}
\centerline {\bol and}
\centerline {\bol Extended Pieri Rules }
\sa
\centerline {\ita
	F.\ Bergeron\footnote{${}^*$}{Work carried out with support from NSERC and FCAR grant.}, N.\ Bergeron${}^*$,
	A.\ M.\ Garsia\footnote{${}^{\dag}$}{Work carried out under NSF grant support.},
	M.\ Haiman${}^{\dag}$,
	and G.\ Tesler${}^{\dag}$}
\sa

{\narrower

\noindent{\bol Abstract. }
The lattice cell in the ${i+1}^{st}$
row and ${j+1}^{st}$ column of the positive quadrant of the plane
is denoted $(i,j)$.  If $\mu$ is a partition of $n+1$, we
denote by $\mu/ij$ the diagram obtained by removing the cell $(i,j)$
from the (French) Ferrers diagram of $\mu$.  We set
$\Delta_{\mu/ij}=\det \|\, x_i^{p_j}y_i^{q_j}\, \|_{i,j=1}^n$,
where
$(p_1,q_1),\ldots ,(p_n,q_n)$ are the cells of $\mu/ij$,
and let ${\bf M}_{\mu/ij}$ be the linear span of the partial derivatives
of $\Delta_{\mu/ij}$. The bihomogeneity of $\Delta_{\mu/ij}$ and its
alternating nature under the diagonal action of $S_n$ gives
${\bf M}_{\mu/ij}$ the structure of a bigraded $S_n$-module. We conjecture
that ${\bf M}_{\mu/ij}$ is always a direct sum of $k$ left regular
representations of $S_n$, where $k$ is the number of cells that are
weakly north and east of $(i,j)$ in $\mu$.  We also make a number of
conjectures describing the precise nature of the bivariate Frobenius
characteristic of ${\bf M}_{\mu/ij}$ in terms of the theory
of Macdonald polynomials.  On the validity of these conjectures, we
derive a number of surprising identities.  In particular, we obtain a
representation theoretical interpretation of the coefficients
appearing in some Macdonald Pieri Rules.

}

\sa

\heading{\bol Introduction}

The lattice cells of the positive plane quadrant will be assigned coordinates $ i,j\geq 0$ as 
indicated in the figure below.
$$
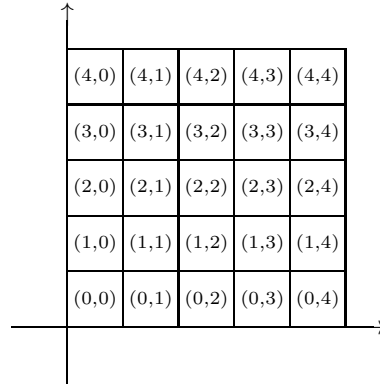
\newcount\drawcntB
\mylength=.037cm
\hskip 20\mylength%
\hbox{%
	\multiputgrid(100,100)(-20,0){5}(0,-20){5}{\rectangle(0,0)(20,20)}
	\hseg(-20,0;130) \myput(115,0){\llap{\lower 2.65pt\hbox to0pt{\hss$\rightarrow$}}}
	\vseg(0,-20;130) \myput(0,110){\vtop{\hbox to0pt{\kern.6pt\hss$\uparrow$\hss}}}
	\drawcntA=0\drawdimX=10.5\mylength\loop\ifnum\drawcntA<5{%
		\drawcntB=0\drawdimY=10\mylength\loop\ifnum\drawcntB<5%
		\raise\drawdimY\hbox to0pt{\hskip\drawdimX%
			\hbox to0pt{\hss${}_{(\the\drawcntB,\the\drawcntA)}$\hss}%
			\hss}%
		\advance\drawcntB 1\advance\drawdimY 20\mylength\repeat%
	}\advance\drawcntA 1\advance\drawdimX 20\mylength\repeat%
}\hskip 140\mylength
$$
A collection of distinct lattice cells will be briefly referred to as a ``{\ita lattice diagram}.''
Given a partition $\mu=(\mu_1\geq \mu_2\geq \cdots \geq \mu_k>0)$, the lattice diagram
with cells
$$
\{\ssp (i,j)\ssp :\ssp 0\leq i\leq k-1\ssp ;\ssp  0\leq j\leq \mu_{i+1}-1\ssp \}\ess ,
$$
as customary, will be called a ``{\ita Ferrers diagram}.'' It will be convenient to use the symbol
$\mu$ for the partition as well as its Ferrers diagram. 
\sa

Given any sequence of lattice cells 
$$
L\ses \{(p_1,q_1)\scs (p_2,q_2)\scs\ldots \scs (p_n,q_n)\}\ess ,
\eqno {\rm I}.1
$$
we define the ``{\ita lattice determinant}''
$$
\DD_L(x;y)\ses {1\over p!q!}\ess  \det \big\|\ssp  x_i^{p_j}y_i^{q_j}\big\|_{i,j=1}^n \ess,
\eqno {\rm I}.2
$$ 
where $p!=p_1!\,p_2!\cdots p_n!$ and $q!=q_1!\,q_2!\cdots q_n!$.
We can easily see that $\DD_L(x;y)$ is a polynomial different
from zero if and only if $L$ consists of $n$ distinct lattice cells. Note also that 
$\DD_L(x;y)$ is  bihomogeneous of degree $|p|=p_1+\cdots +p_n$ in $x$ 
and  degree $|q|=q_1+\cdots +q_n$  in  $y$. It will be good that the definition in I.2
associates a unique polynomial to $L$, as a geometric object. To this end we shall require
that the list of lattice cells in I.1 be given in increasing lexicographic order.
This amounts to listing the cells of $L$ in the order they are encountered as we proceed 
from left to right and from the lowest to the highest.
\sas

Given a polynomial $P(x;y)$,
the vector space spanned by all the partial derivatives of $P$ of all orders
will be denoted 
$\CL_\del[P]$. We recall that the ``{\ita diagonal action}'' of $S_n$ on a polynomial
$$
P(x;y)\ses P(\xon;\yon)
$$ 
is defined by setting for a permutation $\sig=(\sig_1,\sig_2,\ldots ,\sig_n)$
$$
\sig \ssp P(x;y)\ses P(x_{\sig_1},x_{\sig_2},\ldots ,x_{\sig_n};y_{\sig_1},y_{\sig_2},\ldots ,y_{\sig_n})
\ess .
$$
It is easily seen from the definition $\ssp {\rm I}.1\ssp $ that  $\DD_L$ is
an alternant under the diagonal action. This given, it follows that for any
lattice diagram $L$ with $n$ cells, the vector space
$$
\BM_L\ses \CL_\del[\DD_L]
$$
is an $S_n$-module. Since $\DD_L$ is  bihomogeneous, this module affords a natural  bigrading.
Denoting by $\CH_{r,s}[\BM_L]$ the subspace consisting of the bihomogeneous elements
of degree  $r$ in $x$ and degree $s$ in $y$, we have the direct sum decomposition
$$
\BM_L\ses \bigoplus_{r=0}^{|p|}\ssp \bigoplus_{s=0}^{|q|}\ssp \CH_{r,s}[\BM_L]\ess ,
$$
and the polynomial
$$
F_L(q,t)\ses \sum_{r=0}^{|p|}\ssp \sum_{s=0}^{|q|}\ssp t^r\ssp q^s \ssp \dim \CH_{r,s}[\BM_L]
$$
gives the ``{\ita bigraded  Hilbert series}'' of $\BM_L$. In this vein,
since each of the subspaces $\CH_{r,s}[\BM_L]$ is necessarily 
an $S_n$-submodule,  
we can also set
$$
C_L(x;q,t)\ses \sum_{r=0}^{|p|}\ssp \sum_{s=0}^{|q|}\ssp t^r\ssp q^s \ssp \Fch \CH_{r,s}[\BM_L]
\eqno {\rm I}.3
$$
where $\ch \CH_{r,s}[\BM_L]$ denotes the character of $\CH_{r,s}[\BM_L]$ and
$\ssp \Fch \CH_{r,s}[\BM_L]$ denotes the image of $\ch \CH_{r,s}[\BM_L]$ under the Frobenius map
$\CF$ which sends the irreducible character $\chi^\la$  into the Schur function $S_\la$.
The ``$x$'' in $C_L(x;q,t)$ is only to remind us that it is a symmetric function in the 
infinite alphabet $x_1,x_2,x_3,\ldots$ (as customary in [20]), and we should not confuse it with the 
``$x$'' appearing in $\DD_L(x;y)$. This may be unfortunate, but it is too much of an ingrained notation
to be altered at this point. This notation should create no problems since all computations with symmetric polynomials 
are seldom performed in terms of the variables, but rather in terms of the classical 
symmetric function bases. For instance, if $f$  is a symmetric polynomial, by writing
$$
\del_{p_1}\ssp f
$$
we mean the symmetric polynomial obtained by expanding $f$ in terms of the power basis 
and differentiating the result  with respect to $p_1$ as if $f$ were
a polynomial in the indeterminates $p_1,p_2,p_3,\ldots\,$.
Now it is known and easy to prove that for any Schur function $S_\la$ we have
$$
\del_{p_1}\ssp S_\la \ses \sum_{\nu\RA\la} \ssp S_\nu
$$
where ``$\nu \RA \la$'' is to mean that the sum is carried out over partitions $\nu$ 
that are obtained from $\la$ by removing one of its corners. Since, when $\la$ is a  partition of $n$,
we have the well-known ``branching rule'':

$$
\chi_\la\ess  \big\downarrow_{S{n-1}}^{S_n} \ses \sum_{\nu\RA\la} \ssp \chi^\nu\ess ,
$$
we  see that we must have
$$
\del_{p_1}\ssp C_L(x;q,t) \ses  \sum_{r=0}^{|p|}\ssp 
\sum_{s=0}^{|q|}\ssp t^r\ssp q^s \ssp \CF\ssp 
\Bigl( \ch \CH_{r,s}[\BM_L]\ssp \big\downarrow_{S{n-1}}^{S_n}\ssp  \Bigr)\ess .
$$
In other words, $\del_{p_1}\ssp C_L(x;q,t)$ gives the bigraded Frobenius characteristic of
$\BM_L$ restricted to $S_{n-1}$.
In particular we see that we must necessarily have (for any lattice diagram $L$ with $n$ cells)
$$
F_L(q,t)\ses \del_{p_1}^n\ssp C_L(x;q,t)\ess .
\eqno {\rm I}.4
$$
\sa

Computer experimentation with a limited number of cases suggests that the following may hold true:
\sa

\heading{\bol Conjecture I.1}

{\ita For any Lattice diagram $L$ with $n$ cells, the module $\BM_L$ decomposes
into a direct sum of left regular representations of $S_n$.
}
\sa

 Unfortunately, the complexity of computing 
$C_L(x;q,t)$ for large lattice diagrams prevents us from gathering sufficiently strong evidence 
in support of this  conjecture. 
However, the situation is quite different for lattice diagrams obtained by removing a single cell from
a partition diagram. It develops that in this case we have tools at our disposal which
allow us to convert our experimental evidence into a collection of conjectures asserting
that the Frobenius characteristics $C_L(x;q,t)$ satisfy some truly remarkable  recurrences.
Since the latter may be expressed as very precise and explicit symmetric function identities, 
we have been in a position to obtain overwhelming computational and theoretical evidence
in their support. To see how this comes about we need to state some auxiliary results
whose proofs will be found in the next section.
To begin with we have the following useful fact:
\sas

\heading{\bol Proposition I.1}

{\ita Let $L=\{  (p_1,q_1)\scs  (p_2,q_2)\scs \ldots \scs  (p_n,q_n)\ssp \}$ be a lattice diagram. 
Then for any integers $h,k\geq 0$ 
(with $h+k\geq 1$) we have
$$
\sum_{i=1}^n\ssp\del_{x_i}^h\ssp  \del_{y_i}^k\ssp \DD_L(x;y)\ses
\sum_{i=1}^n\ssp \eee(L\DA_{hk}^i)\ssp \DD_{L\DA_{hk}^i}(x;y)
$$
where
$$
L\DA_{hk}^i=\{  (p_1,q_1)\scs \ldots  (p_i-h,q_i-k)\scs \ldots \scs  (p_n,q_n)\ssp \}
\eqno {\rm I}.5
$$
and the coefficient $\eee(L\DA_{hk}^i)$ is different from zero only if $(p_i-h,q_i-k)$ is
in the positive quadrant and $L\DA_{hk}^i$ consists of $n$ distinct cells, in which case
it is given by the sign of the permutation that rearranges the pairs in I.5 in
increasing lexicographic order. 
}
\sas

If $\mu$ is a partition of $n+1$, we shall denote by $\mu/ij$ the lattice diagram obtained by removing
the cell $(i,j)$ from the diagram of $\mu$. We shall refer to the cell $(i,j)$ as the ``{\ita hole}'' of 
 $\mu/ij$. We can easily see that the Proposition  I.1 has the following immediate corollary:
\sas

\heading{\bol Proposition I.2}

{\ita For any partition $\mu$ and  $(i,j)\in \mu$ we have
$$
\sum_{i=1}^n\ssp\del_{x_i}^h\ssp  \del_{y_i}^k\ssp \DD_{\mu/ij}(x;y)\ses
\cases
{\pm\ssp \DD_{\mu/i+h,j+k}(x;y)& if $(i+h,j+k)\in \mu$\cr
\cr
0 & otherwise\cr
}
$$
where the sign is ``$+$'' if there is an odd number of cells (in the lex order) between
$(i,j)$ and $(i+h,j+k)$ and is ``$-$'' otherwise.}
\sa 
It will be convenient to write $(i,j)\leq (i',j')$ meaning
$\{i\leq i'\ssp \&\ssp j\leq j'\}$. This given, the collection of cells
$$
\{(i',j')\in \mu\ssp :\ssp (i,j)\leq (i',j')\ssp \}
$$
will be called the ``{\ita shadow}'' of $(i,j)$ in $\mu$.
It is a translation of the Ferrers diagram of a partition.
 Let us also set
$$
D_x\ses \sum_{i=1}^n\del_{x_i}
\scs\ess\ess
 D_y\ses \sum_{i=1}^n\del_{y_i}
\ess\ess\ess\ess\hbox { and}\ess\ess\ess 
D_{hk}\ses \sum_{i=1}^n\del_{x_i}^h\del_{y_i}^k\ess .
$$ 
Now we have the following important consequences of Proposition I.2:
\sa

\heading{\bol Proposition I.3}

{\ita Let $\mu$ be a partition of $n+1$. Then for any pair of cells  $(i,j)\scs(i+h,j+k)\in \mu$ 
 we have
$$
D_x^hD_y^k\ssp \BM_{\mu/ij}\ses D_{hk}\ssp  \BM_{\mu/ij}\ses \BM_{\mu/i+h,j+k}
\eqno {\rm I}.6
$$
meaning that both  $D_x^hD_y^k$ and $D_{hk}$ are surjective linear maps. 
In particular we have the inclusion
$$
\BM_{\mu/i'j'}\ess \con\ess \BM_{\mu/ij}
\eqno {\rm I}.7
$$ 
for all cells $(i',j')$ in the shadow of $(i,j)\ess .$
}
\sap

\heading{\bol Proposition I.4}

{\ita  The collection of polynomials
$$
\{\ssp \DD_{\mu/i'j'}(x;y)\ssp :\ssp (i',j')\in \mu \hbox{ and } (i',j')\geq (i,j)
\ess \}
$$
form a basis for the submodule of alternants of $\BM_{\mu/ij}\ess $.
}
\sas

Note that  Conjecture I.1, combined with this result, leads us to a 
more precise statement concerning our modules  $\BM_{\mu/ij}\ess $:
\sas

\heading{\bol Conjecture I.2}

{\ita For any $\mu\part n+1$ and any $(i,j)\in \mu\,$, the $S_n$-module $\BM_{\mu/ij}$
decomposes into the direct sum of $m$ left regular representations
of $S_n\,$, where 
$m$ gives the number of cells in the shadow of $(i,j)$.
}
\sa

This may be viewed as an extension of the conjecture made in [7] that 
for any $\mu\part n$ the module $\BM_\mu$ gives a bigraded version of the left regular representation
of $S_n$. It was also conjectured in [7] that the bivariate Frobenius characteristic of
$\BM_\mu$ is given by the the symmetric polynomial
$$
\TH_\mu(x;q,t)\ses \sum_{\la \part n}\ssp S_\la(x)\ssp \TK_{\la\mu}(q,t)
	\ess ,
\eqno {\rm I}.8
$$
where the coefficients $\TK_{\la\mu}(q,t)$ are related to the Macdonald [19] $q,t$-Kostka
coefficients $K_{\la\mu}(q,t)$ by the formula
$$
\TK_{\la\mu}(q,t)\ses t^{n(\mu)}\ssp K_{\la\mu}(q,1/t)\ess .
$$
Here as in [20], for any partition $\mu$  we set 
$$
n(\mu)\ses \sum_i\ssp (i-1)\ssp \mu_i\ess .
\eqno {\rm I}.9
$$
In the present notation, the latter conjecture may be expressed by writing
$$
C_\mu(x;q,t)\ses \TH_\mu(x;q,t)\ess .
\eqno {\rm I}.10
$$
For this reason, we shall refer to this equality as the $C=\TH$ conjecture. 
Macdonald conjectured in [19] that $K_{\la\mu}(q,t)$ is always a
polynomial in $q,t$ with positive integer coefficients. Though recently in
[12], [13], [15], [16] and [18]
it was shown that they are polynomials with integer coefficients, the positivity
still remains to be proved. Of course, the equality in I.10 would completely
settle the positivity conjecture. It follows from Macdonald's work that
$$
\TK_{\la\mu}(1,1)\ses f_\la\ses \#\{ \ess \hbox{standard tableaux of shape $\la$}\ess \}\ssp . 
$$
Thus I.10 is consistent with the statement that $\BM_\mu$ is a bigraded version
of the left regular representation of $S_n$. Now it develops that there is also 
a way of extending the $C=\TH$ conjecture to the lattice diagrams $\mu/ij$. The point of departure
here is the following remarkable fact.
\sa

\heading{\bol Proposition I.5}

{\ita  For any $\mu\part n+1$ we have 
$$
C_{\mu/00}(x;q,t)\ses  \sum_{r=0}^{|p|}\ssp 
\sum_{s=0}^{|q|}\ssp t^r\ssp q^s \ssp \CF\ssp 
\Bigl( \ch \CH_{r,s}[\BM_\mu]\ssp \big\downarrow_{S_{n}}^{S_{n+1}}\ssp  \Bigr)
\ses 
\del_{p_1}\ssp C_\mu(x;q,t)\ess . 
\eqno {\rm I}.11
$$
}
\sa

\noindent
Thus on the $C=\TH$ conjecture we should have
$$
C_{\mu/00}(x;q,t)\ses \del_{p_1}\ssp \TH_\mu(x;q,t)\ess .
\eqno {\rm I}.12
$$ 
Since the operator $\del_{p_1}$ is the adjoint of multiplication by the elementary symmetric function $e_1$ 
with respect to the Hall scalar product, it may be derived from one of the Macdonald Pieri rules
(see [6]) that we have
$$
\del_{p_1}\ssp \TH_\mu(x;q,t)\ses \sum_{\nu\RA \mu}\ssp c_{\mu\nu}(q,t)\ssp \TH_\nu(x;q,t)
\eqno {\rm I}.13
$$
with 
$$
{ c}_{\mu\nu}(q,t)\ses 
\prod_{s\in \CR_{\mu/\nu}}
{ t^{l_\mu (s)}-q^{a_\mu (s)+1}\over { t^{l_\mu (s)}-q^{a_\mu (s)}} } 
\ssp
\prod_{s\in \CC_{\mu/\nu}}
{
q^{a_\mu (s)}- t^{l_\mu (s)+1}
\over 
{q^{a_\mu (s)}-t^{l_\mu (s)}} 
}
\ess .
\eqno {\rm I}.14
$$ 
Here $\CR_{\mu/\nu} $ (resp. $\CC_{\mu/\nu} $) denotes the set of 
lattice squares of $\nu$ that are in the same row (resp. same column) as the
cell we must remove from $\mu$ to obtain $\nu$ and for any cell
$s\in \mu$, the parameter $l_\mu(s)$ gives the number of cells of $\mu$ that 
are strictly north of $s$ and $a_\mu(s)$ gives the number of cells that
are strictly east. In view of I.13, we may rewrite I.11 in the form
$$
C_{\mu/00}(x;q,t)\ses \sum_{\nu\RA \mu}\ssp c_{\mu\nu}(q,t)\ssp \TH_\nu(x;q,t)\ess .
\eqno {\rm I}.15
$$
Now extensive computations with the modules $\BM_{\mu/ij}$ have revealed that   
a truly remarkable analogue of this formula may hold true for all the Frobenius 
characteristics $C_{\mu/ij}(x;q,t)$; we can state it as follows: 
\sa

\heading{\bol Conjecture I.3}

{\ita For any $(i,j)\in \mu$ we have
$$
C_{\mu/ij}(x;q,t)\ses \sum_{\rho\RA \tau}\ssp c_{\tau\rho} (q,t)\ssp \TH_{\mu-\tau+\rho\ssp}(x;q,t)
	\ess ,
\eqno {\rm I}.16
$$
where $\tau$ denotes the Ferrers diagram contained in the shadow of $(i,j)$
and the symbol ``$\mu-\tau+\rho$'' is to represent replacing $\tau$ by $\rho$ in the shadow of $(i,j)$.
}
\sa

The following result not only reveals the true nature of I.16,
but sheds some surprising 
light on the Macdonald Pieri rule corresponding to the identity in I.13.
\sas

\heading{\bol Theorem I.1}

{\ita The validity of I.16 for all $(i,j)\in \mu$ is equivalent to 

\sas

\item{(a)} the four term recursion
$$
C_{\mu/ij}\ses {t^l-q^{a+1}\over t^l-q^a}\sssp  C_{\mu/i,j+1}\ssp +\ssp {t^{l+1}-q^{a}\over t^l-q^a}\sssp  C_{\mu/i+1,j}\ssp -\ssp
{t^{l+1}-q^{a+1}\over t^l-q^a} \sssp C_{\mu/i+1,j+1}\ess , 
\eqno {\rm I}.17
$$
where $l$ and $a$  give the number of cells that are respectively north and east of $(i,j)$ in  $\mu$,
\sa 

\item{(b)} together with the boundary conditions that the terms 
\hbox {$ C_{\mu/i,j+1},$ $C_{\mu/i,j+1}$  or $C_{\mu/i,j+1}$}
\item {} are equal to zero when the corresponding
cells $(i,j+1),$ $(i+1,j)$ or  $(i+1,j+1)$
\item {} 
fall outside of $\mu\ess ,$ and are equal to 
$ \TH_{\mu/i,j+1},$ $ \TH_{\mu/i,j+1}$ or  $ \TH_{\mu/i,j+1}$  when any of the
\item {} corresponding
cells \hbox {is a corner of $\mu$}.
}
\sa 

\noindent
Now a crucial development here is that I.17 has a representation  theoretical interpretation
that strongly suggests an inductive argument for proving both Conjectures I.2 and I.3. 
To present it we must introduce some notation. For a given $\ess (i,j)\in \mu\ess $, let
$\BK_{ij}^x$ denote the kernel of the operator $D_x$ as a map of $\BM_{ij}$ onto $\BM_{i+1,j}$.
Similarly, let $\BK_{ij}^y$ be the kernel of $D_y$ as a map of $\BM_{ij}$ onto $\BM_{i,j+1}$.
It will also be convenient to denote by $K_{ij}^x$ and  $K_{ij}^y$ the corresponding
bivariate Frobenius characteristics. Note that since $\ssp \BM_{i,j+1}\con\BM_{i,j}\ssp $ and
$\ssp \BM_{i+1,j}\con\BM_{i,j}\ssp $ we see that we must have
$$
\BK_{i,j+1}^x\ssp \con \ssp \BK_{ij}^x
\ess\ess\ess 
\hbox{ as well as 
}\ess\ess\ess 
\BK_{i+1,j}^y\ssp \con \ssp \BK_{ij}^y
	\ess .
$$
Note further that if $\mu\part n+1$ all of these vector spaces are $S_n$-invariant
and the quotients 
$$
\BA_{ij}^x\ses \BK_{ij}^x\ssp /   \BK_{i,j+1}^x\ssp 
\ess\ess\ess 
\hbox{and
}\ess\ess\ess 
\BA_{ij}^y\ses \BK_{ij}^y\ssp /  \BK_{i+1,j}^y\ssp  
\eqno {\rm I}.18
$$
are well-defined bigraded $S_n$-modules. Let $\ssp A_{ij}^x\ssp $ and $\ssp A_{ij}^y\ssp $ denote 
their respective Frobenius characteristics.  This given, a simple linear algebra argument gives that we have
the following relations:
\sa

\heading{ \bol  Proposition I.6}

$$
\eqalign{
a)\ess\ess\ess K_{ij}^x\ses C_{\mu/ij}- t\ssp C_{\mu/i+1,j}
\ess  
  &,
 \ess\ess\ess
K_{ij}^y\ses C_{\mu/ij}- q\ssp C_{\mu/i,j+1}\cr \cr
b)\ess\ess\ess 
\ess\ess\ess\ssp\ssp
A_{ij}^x\ses K_{ ij}^x-  \ssp K_{ i,j+1}^x
\ess 
 & , 
\ess\ess\ess
A_{ij}^y\ses K_{ ij}^y- \ssp K_{ i+1,j}^y
\cr
}
\eqno {\rm I}.19
$$
{\ita In particular, the recurrence in I.17 may be rewritten in the simple form
$$
t^l\ssp A_{ij}^x\ses q^a \ssp A_{ij}^y\ess .
\eqno {\rm I}.20
$$
 }
\sa

It develops that I.20 encapsulates a great deal of combinatorial and
representation theoretical information. Indeed, a proof of this identity 
may turn out to be the single most important result
in the present theory and in the theory of Macdonald polynomials.
For this reason we shall here and after refer to I.20 as the ``{\ita crucial identity}.''

To be precise, we shall show in Section 1 that I.20 is more than sufficient
to imply the validity of Conjectures I.2 and I.3 and the $q,t$-Kostka positivity conjecture. 
The argument also shows that for $\mu\part n+1$ the modules
$\BA_{ij}^x$ and  $\BA_{ij}^y$ are all left regular representations of $S_n$. It will then
result that in some sense the modules $\BA_{i'j'}^x$ and  $\BA_{i'j'}^y$ with $(i',j')\geq (i,j)$,
yield what may be viewed as an ``atomic'' decomposition of $\BM_{\mu/ij}$ into a direct sum
of left regular representations of $S_n$. 

This given, our basic goal here is to understand  
the representation theoretical significance of I.20 in the hope that it may lead to the construction
of a proof. Now it develops that the methods introduced in [1] can be extended to the present
situation to yield some very precise information concerning the behavior of the
Frobenius characteristics $\BA_{ij}^x$ and  $\BA_{ij}^y$ as $(i,j)$ varies in $\mu$.
One of the main results in [1], translated into the present language, is
an algorithm for decomposing $\BM_{\mu/00}$ as
a direct sum of appropriate intersections of the modules $\BM_\aaa\ssp $ with $\aaa\RA \mu$.
This algorithm is based
on a package of assumptions which have come to be referred to as the ``{\ita SF-heuristics}.''
We shall show here that the SF-heuristics can be extended to yield a similar decomposition 
for all the modules $\BM_{\mu/ij}$. We should mention that, as was the case in [1], 
all these decompositions, combined with the $C=\TH$ conjecture, yield (via the Frobenius map)
a variety of symmetric function identities for which we have
overwhelming experimental and theoretical confirmation through the theory of Macdonald polynomials.

To state our results we need to review and extend some of the  notation introduced in [1].
The reader is referred to [1] for the motivation underlying these definitions.
\sas

Here and after, if $P(x;y)=P(\xon;\yon)$ is a polynomial, we will let
$P(\del_x;\del_y)$, or even simply 
$P(\del)$, denote the differential operator obtained by replacing, for each $i$ and $j$,
$x_i$ by $\del_{x_i}$ and $y_j$ by $\del_{y_j}$. This given, we shall set for 
any two polynomials $P(x;y)$ and $Q(x;y)$
$$
\LL P\scs Q\RR\ses P(\del_x;\del_y)\ssp Q(x;y)\ssp \big|_{x=y=0}\ess . 
\eqno {\rm I}.21
$$
It easily seen  that this defines a scalar product which is invariant under the
diagonal action of $S_n$. That is, for each $\sig\in S_n$ we have
$$
\LL\sig  P\scs Q\RR\ses 
\LL P\scs\sig^{-1} Q\RR\ess . 
\eqno {\rm I}.22
$$
Moreover, since the monomials $\{x^py^q\}_{p,q}$ form an orthogonal set under this scalar product, 
pairs of polynomials of different bidegree will necessarily be orthogonal to each other. 
\sas
 
If $\DD(x;y)$ is any diagonally alternating polynomial, the space $\BM_\DD=\CL[\del_x^p \del_y^q\DD(x;y)]$
spanned by all partial derivatives of $\DD(x;y)$ will necessarily be $S_n$-invariant. If $\DD$ is
bihomogeneous of bidegree $(r_0,s_0)$, then  $\BM_\DD$ has a sign-twisting, bidegree-complementing
isomorphism we shall denote by $\flip_\DD$, which may be defined by setting for each $P\in \BM_\DD$
$$
\flip_\DD\ssp P(x;y)\ses P(\del_x;\del_y)\ssp \DD(x;y)\ess .
\eqno {\rm I}.23
$$
In particular, this implies that the bivariate Frobenius characteristic $\Phi_\DD(x;q,t)$ of $\BM_\DD$
will necessarily satisfy the identity
$$
\Phi_\DD(x;q,t)\ses t^{r_0}q^{s_0}\ssp \om \Phi_\DD(x;1/q,1/t)
$$
where $\om$, as customary, denotes the involution that sends the Schur function $S_\la$ into
$S_{\la '}$. It will be convenient to set, for any symmetric polynomial $\Phi(x;q,t)$
 with coefficients rational functions of $q$ and $t$:
$$ 
\da \Phi(x;q,t)\ses \om \Phi(x;1/q,1/t)\ess .
\eqno {\rm I}.24
$$
It can also be seen that if $\BM_1\con \BM_\DD$ is any bigraded $S_n$-invariant submodule of $\BM_\DD$ 
with bivariate Frobenius characteristic $\Phi_1(x;q,t)$ then the subspace
$$
\flip_\DD\ssp \BM_1\ses \{ \ssp \flip_\DD\ssp P\ssp :\ssp P\in \BM_1\ssp \}
$$ 
is also $S_n$-invariant, bigraded, and its bivariate Frobenius characteristic is given by the formula
$$
\Fch \flip_\DD\ssp \BM_1\ses t^{r_0}q^{s_0}\ssp \da \Phi_1(x;q,t)
\ess .
\eqno {\rm I}.25
$$ 
Both the flip map and our scalar product have a number of easily verified properties that will be used
in our development. To begin with, we should note that the orthogonal complement 
$\BMperp$ of $\BM_\DD$ with respect to $\LL \scs \RR$,
that is the space
$$
\BMperp=\{\ssp Q(x;y)\ssp :\ssp \LL P,Q\RR =0\ess\ess \forall\ess\ess P\in \BM_\DD\ssp \}\ssp , 
$$ 
 consists of all the polynomial
differential operators that kill $\DD(x;y)$. More precisely, 
$$
\BMperp\ses \{\ssp Q(x;y)\ssp :\ssp Q(\del_x;\del_y)\ssp \DD(x;y)\ses 0\ssp\}\ess . 
\eqno {\rm I}.26
$$ 
Note that since 
$$
\LL  P\scs \flip_\DD Q\RR\ses P(\del_x;\del_y)Q(\del_x;\del_y)\DD(x;y)\ssp \big|_{x=y=0}\ess , 
$$
we see that $\flip_\DD$ is self-adjoint. That is, for all $P$ and $Q$,  we have 
$$
\LL \flip_\DD P\scs  Q\RR\ses \LL  P\scs \flip_\DD Q\RR\ess .
\eqno {\rm I}.27
$$
For an element $P\in\BM_\DD$, the (necessarily) unique $P_1\in \BM_\DD$ such that
$$
P(x;y)\ses P_1(\del_x;\del_y)\ssp \DD(x;y) 
$$
will be denoted by $\flip_\DD^{-1}P$.
The following result will play a basic role in our development:
\sas

\heading{\bol Proposition I.7}

{\ita
Let $D(x;y)$ be a polynomial, $\Delta(x;y)$ be an alternant, and set
$\TDD(x;y)=D(\del_x;\del_y)\Delta(x;y)$.
Let $\BM_\DD$ (resp., $\BTM$) be the module spanned by all
partial derivatives of $\Delta$ (resp., $\TDD$).  
Then $\BTM$ is a submodule of $\BM_\DD$ and
$D(\del_x;\del_y)$ is a surjective map from $\BM_\DD$ to $\BTM$.
Letting $\BK$ denote the kernel of this map, we have that
$$
\BM_\DD\,\cap\,\BTMperp \ses \flip_\DD^{-1}\ssp \BK\ess .
\eqno {\rm I}.28
$$
This gives the direct sum decompositions 
$$
\eqalign{
&a)\ess \ess \BM_\DD \ses  \ess\BTM \ess \oplus_\perp\ess   \flip_\DD^{-1}\ssp \BK  ,\cr
&b)\ess \ess \BM_\DD \ses  \flip_\Delta \BTM \hskip .07truein\oplus\ess\hskip .04truein \BK  \ess
\ess .\cr
}
\eqno {\rm I}.29
$$
where the symbol ``$\oplus$'' denotes the direct sum of disjoint spaces,
and ``$\oplus_\perp$'' further denotes that these spaces are orthogonal to
each other.}
\sa

Now let $\mu$ be a fixed partition of $n+1$ and let 
$$
\Pred(\mu)\ses \bigl\{\ssp \nu^{(1)}\scs \nu^{(2)}\scs \ldots \scs \nu^{(d)}\ssp \bigr\}
\eqno {\rm I}.30
$$
be the collection of partitions obtained by removing one of the corners of $\mu$. 
For a pair  $\nu\RA \mu$,  it will be convenient to denote by $\mu/\nu$ the corner cell 
we must remove from $\mu$ to get $\nu$. To be specific, we shall assume that
the partitions in I.30 are ordered so that the corner $\mu/\nu^{(k)}$
is northwest of the corner $\mu/\nu^{(k+1)}$. Similarly, for a given cell $(i,j)\con \mu$ 
let
$$
\Pred_{ij}(\mu)\ses \bigl\{\ssp \aaa^{(1)}\scs \aaa^{(2)}\scs \ldots \scs \aaa^{(m)}\ssp \bigr\}
\eqno {\rm I}.31
$$
be the subset of $\Pred(\mu)$ consisting of the $\nu^{(k)}$ such that $\mu/\nu^{(k)}$ is
in the shadow of $(i,j)$. We again assume that the $\aaa^{(i)}$ are labelled so that, for $i=1,\ldots,m-1$, the corner  
$\mu/\aaa^{(i)}$ is northwest of the corner $\mu/\aaa^{(i+1)}$.
\sas

Following Macdonald [20] we call the ``{\ita coleg}'' and ``{\ita coarm}''
of a lattice cell $s\in \mu$ the numbers
$l'_\mu(s)$, and  $a'_\mu(s)$ of cells that are respectively strictly  south and strictly
west of $s$ in $\mu\,$. In our notation, if $s=(i,j)$ then $l'_\mu(s)=i$ and $a'_\mu(s)=j$. We
shall call the monomial $w(s)=t^{l'_\mu(s)}q^{a'_\mu(s)}$ the ``{\ita weight}'' of $s$. 
For any lattice diagram $L$ we set
$$
T_L= \prod_{s\in L}\ssp w(s)\ess .
$$
We shall also denote by $\nabla$ the linear operator defined by setting for every partition $\mu$
$$
\nabla\ssp \TH_\mu(x;q,t)\ses T_\mu\ssp \TH_\mu(x;q,t)\ess .
\eqno {\rm I}.32
$$
Since the polynomials $\TH_\mu(x;q,t)$ form a symmetric function basis, I.32 defines $\nabla$
as an operator acting on all symmetric polynomials. For two subsets $T\con S\con\Pred(\mu)\ess  $ set
$$
\BM_S^T\ses 
\left(\;\bigcap_{\aaa\in T} \BM_\aaa \right)
\cap \left( \Biggl(\;\sum_{\bbb\in S-T} \BM_\bbb \Biggr)
\cap \left(\; \bigcap_{\aaa\in T} \BM_\aaa \right) \right)^\perp
\eqno {\rm I}.33
$$
where the symbols
``$\bigcap$'' and ``$\sum$'' denote intersection and sum (not usually direct)
of vector spaces, and
``$\perp$'' denotes the operation of taking orthogonal complements
with respect to the scalar product defined in I.21.
Since this scalar product is invariant under the diagonal action of $S_n\ssp ,$ we  see that
$\BM_S^T$ is a well-defined $S_n$-module, and its bivariate Frobenius characteristic
will be denoted by $\phi_S^T$. One of the assertions of the SF-heuristics 
is that in the linear span 
$$
\CL[\ssp \TH_\aaa\ssp :\ssp \aa\in S\ssp]
$$
we  have $\ssp m= |S| \ssp $ Schur positive symmetric polynomials
$$
\phi^{(1)}_S\scs \phi^{(2)}_S\scs\ldots \scs  \phi^{(m)}_S 
$$  
such that for any $T\con S$ of cardinality $k$ we have
$$
\phi_S^T\ses {\phi^{(k)}_S\over \prod_{\aaa\in S-T}T_\aaa}\ess .
\eqno {\rm I}.34
$$
It is  also a consequence of the SF-heuristics  that for $k=1,\ldots ,m-1$ we can set
$$
\phi^{(k)}_S\ses (-\nabla )^{m-k}\ssp  \phi^{(m)}_S\ess ,
\eqno {\rm I}.35
$$
while $\phi^{(m)}_S$ itself can be computed from the formula
$$
\phi^{(m)}_S\ses \sum_{\aaa\in S}\ssp \Bigl(\prod_{\bbb\in S/\{\aaa\}}{1\over 1-T_\aaa/T_\bbb}\ssp \Bigr)\ssp \TH_\aaa
\ses \sum_{\aaa\in S}\ssp \Bigl(\prod_{\bbb\in S/\{\aaa\}}{1\over 1-\nabla/T_\bbb}\ssp \Bigr)\ssp \TH_\aaa\ess .
\eqno {\rm I}.36
$$
To be consistent with the notation we adopted in [1] we shall use the symbols $\phi_\mu$ or $\phi_\mu^{(k)}$ to denote
$\phi_S^{(m)}$ or $\phi_S^{(k)} $ when $S$ consists of all the predecessors of $\mu$. In this vein, it will also be convenient
to set, for any subset $S\con\Pred(\mu)$,
$$
\ ^cS\ses \Pred(\mu)- S\ess .
$$
By comparing the expansion of $\phi^{(m)}_S$ with that of
$\phi_\mu=\phi^{(m')}_{S'}$ (where $S'=\Pred(\mu)$ has cardinality $m'$)
in I.36, it follows that
$$
\phi^{(m)}_S\ses \biggl(\ssp \prod_{\bbb\in \ ^cS}\Bigl(\ssp 1-{\nabla\over  T_\bbb }\Bigr) \biggr)\ssp  \phi_\mu\ess .
\eqno {\rm I}.37
$$
In particular, when $S$ consists of a single partition
$\nu^{(i)}\in \Pred(\mu)$, this reduces to
$$
\TH_{\nu^{(i)}}(x;q,t)\ses 
\biggl(\ssp\prod_{j=1\,;\,j\neq i}^d\Bigl(\ssp 1-{\nabla\over T_{\nu^{(j)}} } \Bigr) \biggr)\ssp  \phi_\mu\ess ,
\eqno {\rm I}.38
$$
which may also be rewritten in the form (see 3.19 of [1])
$$
\TH_{\nu^{(i)}}(x;q,t)\ses\sum_{k=1}^d 
\phi_\mu^{(k)}\ssp e_{d-k}\Bigl[{1\over T_{\nu^{(1)}}}+{1\over T_{\nu^{(2)}}}
+\cdots +{1\over T_{\nu^{(d)}}}\sms {1\over T_{\nu^{(i)}}}   \Bigr]\ess . 
\eqno {\rm I}.39
$$
Finally note that if $\nu^{(i)}=\aaa\in S$ then by combining I.37 and I.38  we can also write
$$
\TH_\aaa(x;q,t)\ses  \prod_{\bbb\in S\,;\,\bbb\neq \aaa}\Bigl(\ssp 1-{\nabla\over T_\bbb }\Bigr) 
\prod_{\bbb\in \ ^cS}\Bigl(\ssp 1-{\nabla\over T_\bbb }\Bigr)\ssp \phi_\mu
\ses
  \prod_{\bbb\in S\,;\,\bbb\neq \aaa}\Bigl(\ssp 1-{\nabla\over T_\bbb }\Bigr)\ssp
\ssp   \phi^{(m)}_S
\eqno {\rm I}.40
$$
or equivalently, for $S=\bigl\{\aaa^{(1)},\aaa^{(2)}, \ldots , \aaa^{(m)}\bigr\}$ and $\aaa=\aaa^{(i)}$
$$
\TH_{\aaa^{(i)}}(x;q,t)\ses \sum_{k=1}^m \phi_S^{(k)} 
e_{m-k}\Bigl[{1\over T_{\aaa^{(1)}}}+{1\over T_{\aaa^{(2)}}}+\cdots +
{1\over T_{\aaa^{(m)}}}\sms {1\over T_{\aaa^{(i)}}}   \Bigr] \ess .
\eqno {\rm I}.41
$$
To complete our notation we need to recall that in [13] the weights  of the corners
$$
\mu/\nu^{(1)}\scs \mu/\nu^{(2)}\scs \ldots \scs \mu/\nu^{(d)}
$$
were respectively called
$$
x_1\scs x_2\scs \ldots \scs x_d \ess .
$$
Moreover, if  $\ess x_i=t^{l'_i}  q^{a'_i}\ess $ then we also let  
$$
u_i\ses t^{l'_{i+1}}q^{a'_{i}}\bigsp \hbox{( for $\ess  i=1,2,\ldots ,m-1\ssp $)}
\eqno {\rm I}.42 
$$
be the weights of what we might refer to as the ``{\ita inner corners}'' of $\mu$. 
The picture is completed by setting
$$
u_0\ses t^{l'_1}/q\ess\scs\ess 
u_m\ses q^{a'_m}/t\ess\ess\hbox{ and }\ess\ess 
x_0\ses 1/tq\ess .
\eqno {\rm I}.43 
$$
To appreciate the geometric significance of these weights, in the figure below we illustrate
a $4$-corner case with 
corner cells labelled $A_1$, $A_2$, $A_3$, $A_4$   
and inner corner cells
labelled  $B_0$, $B_1$, $B_2$, $B_3$, $B_4$.
$$
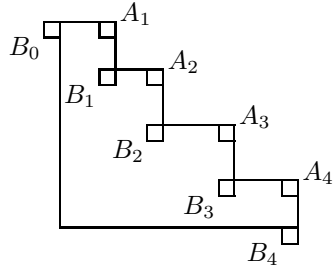
\hbox{%
	\hseg(0,0;300)\corn(300,0)\labLL(300,0;B_4)	
	\drawL(220,300,0,60;A_4,B_3)
	\drawL(130,220,60,130;A_3,B_2)
	\drawL(70,130,130,200;A_2,B_1)
	\drawL(0,70,200,260;A_1,B_0)
	\vseg(0,0;260)					
}\hskip 300\mylength
$$

\sa
 
It was shown in [13] that the products in I.14 giving the coefficients $c_{\mu\nu}(q,t)$ undergo
massive cancellations which reduce them to relatively simpler expressions in terms of the corner weights.
This results in the formula      
$$
c_{\mu\nu^{(i)}}\ses {1\over M}\ssp{1\over x_i}\ssp  
{\prod_{j=0}^d \ssp (x_i-u_j)\over \prod_{j=1\,;\,j\neq i}^d(x_i-x_j) }
\eqno {\rm I}.44 
$$
where for convenience we have set
$$
M=(1-1/t)(1-1/q)\ess .
\eqno {\rm I}.45 
$$
Taking account of the fact that $\ess x_i\ssp T_{\nu^{(i)}}=T_\mu\ssp ,\ess$ formula I.38 can also be written in the form
$$
\TH_{\nu^{(i)}}(x;q,t)\ses  \prod_{j=1\,;\,j\neq i}^d\Bigl(\ssp 1-\nabla{x_i\over T_\mu }\Bigr)\ssp \phi_\mu \ess .
\eqno {\rm I}.46 
$$
It was shown in [1] (Theorem 3.3) that 
using I.44 and I.46 in I.13 yields the following beautiful identities:  
$$
\eqalign{
a)\ess\ess \del_{p_1}\TH_\mu&\ses {1\over M}{T_\mu\over \nabla}\ess 
\biggl(\ssp \prod_{s=0}^d\ssp \Bigl(1-\nabla {u_s\over T_\mu}\ssp \Bigr)
 \biggr)\ssp \phi_\mu\cr
b)\ess\ess 
\del_{p_1}\TH_\mu&\ses
\sum_{k=1}^d\ssp {\phi^{(k)}_S\over T_\mu^{m-k}}\ssp 
{ e_{m+1-k} [x_0+\cdots +x_d]\sms e_{m+1-k} [u_0+\cdots +u_d]\over M}
\ess .\cr
}
\eqno {\rm I}.47 
$$
It develops that using the same argument we can obtain analogous identities for I.16.
To state them we need some notation. Let 
$$
S\ses \Pred_{ij}(\mu)\ses \bigl\{\ssp \aaa^{(1)}\scs \aaa^{(2)}\scs \ldots \scs \aaa^{(m)}\ssp \bigr\}\ess ,
\eqno {\rm I}.48 
$$
and let $\tau$ denote the partition that corresponds to the shadow of $(i,j)$ in $\mu$. That is,
$$
\tau\ses (\mu_{i+1}-j,\mu_{i+2}-j,\ldots ,\mu_{i+1+l}-j)\ssp ,
$$
where $l$ gives the number of cells above $(i,j)$ in $\mu$. 
Finally, let  $x_s^{ij}$ and $u_s^{ij}$  (for  $\ssp 0\leq s\leq m$) be the corner weights of $\tau$. This given,
we can rewrite I.16 in either of the two forms
\sas

\heading{\bol Proposition I.8}
$$
\eqalign{
a)\ess\ess 
C_{\mu/ij}&\ses {1\over M}{T_{\mu/ij}\over \nabla}\ess 
\biggl(\ssp \prod_{s=0}^m\ssp \Bigl(1-\nabla {u_s^{ij}\over T_{\mu/ij}}\ssp \Bigr)
 \biggr)\ssp \phi_S^{(m)}\cr
b)\ess\ess 
C_{\mu/ij}&\ses
\sum_{k=1}^m\ssp {\phi^{(k)}_S\over T_{\mu/ij}^{m-k}}\ssp 
{ e_{m+1-k} [x_0^{ij}+\cdots +x_m^{ij}]\sms e_{m+1-k} [u_0^{ij}+\cdots +u_m^{ij}]\over M}
\ess .\cr
}
\eqno {\rm I}.49 
$$
Formula I.49 a) enables us to obtain completely explicit expressions for the bivariate
Frobenius characteristics of the modules $\BA_{ij}^x$. 
\sas

\heading{\bol Theorem I.2}

{\ita Letting $l$ and $a$ be the leg and arm of $(i,j)$ and assuming
I.48, with the above conventions, we have
$$
A_{ij}^x/q^{a}\ses A_{ij}^y/t^{l}\ses \biggl(\ssp
 \prod_{s=1}^{m-1}\ssp \Bigl(1-\nabla {u_s^{ij}\over T_{\mu/ij}}\ssp
 \Bigr) \biggr)\ssp \phi_S^{(m)}\ess .
\eqno {\rm I}.50 
$$
}

This result has a truly surprising consequence. For a moment let $\Pred(\mu)$
 be as in I.30 and let the weight of $\mu/\nu^{(i)}$ be $\ssp
t^{l'_i}q^{a'_i}$. For any pair $i,j\in [1,m]$ set
$$
R_{i,j}\ses
\{\ssp s\in \mu\ssp :\ssp a_{i-1}'<a'(s)\leq a_i'\ssp\ess ;\ess\ssp l_{j+1}'<l'(s)\leq l_{j}'\ess \}
\ess ,
\eqno {\rm I}.51
$$
where for convenience we set $a_0'=l_{m+1}'=-1$. In words
$R_{i_0,j_0}$ is the subrectangle of $\mu$ consisting of the cells
which have in their shadow only the corner cells
$$
(l_i'\scs a_i')\ess\ess\ess {\rm for}\ess\ess\ess i_0\leq i\leq
j_0\ess .
$$
This given, from I.50 we immediately deduce the following.
\sas

\heading{\bol Theorem I.3}

{\ita The bigraded modules $\BA_{i'j'}^x$ and $\BA_{i'j'}^y\,$, up to a
bidegree shift, remain isomorphic as the cell $(i',j')$ varies in a
rectangle $R_{i,j}$.}
\sap

This paper is divided into five sections. In Section 1 we prove all
the propositions and theorems stated in the Introduction. Some of
these proofs rely on material presented in [1]. The reader will be
well advised to have a copy of that paper at hand in reading the
present work. The main goal in Section 2 is to give a
representation theoretical interpretation of the ``crucial identity''
(I.20).  The basic tool there is an algorithm for constructing bases
for all our modules $\BM_{\mu/ij}$. Since this algorithm is based on
the heuristics proposed in [1], its validity depends on the validity
of those heuristics, which at the present time are still
conjectural. Nevertheless it will be seen that the symmetric function
identities implied by the validity of the algorithm are in complete
agreement with massive computational evidence provided by the theory
of Macdonald polynomials.  In Section 3 we treat in full detail the
case when $\mu$ is a ``hook'' shape and show that all our conjectures
are indeed correct in this case to the finest detail.  In Section 4 we
give a combinatorial argument proving that for all $\mu\part n$, each
of the modules $\BM_{\mu/ij}$ has dimension bounded above by $n!$
times the number of cells in the shadow of $(i,j)$.  Finally, in
Section 5 we show that some of the modules whose existence was
conjectured in [8] have a natural setting within the theory of
``atoms'' we have developed in the present work. In particular we are
able to explain the origin of some puzzling identities derived in [8].
\sas

\heading{\bol 1. Basic properties of our lattice modules.} 
\sas

This section is dedicated to proving all the propositions and theorems
we stated in the introduction. 
\sas

\heading{\bol Proof of Proposition I.1}

For $L=\{  (p_1,q_1)\scs  (p_2,q_2)\scs \ldots \scs  (p_n,q_n)\ssp \}$ we can write
$$
\DD_L(x;y)\ses {1\over p!q!}\sum_{\sig \in S_n}\ssp
\sign(\sig)\ssp  x_{\sig_1}^{p_1} y_{\sig_1}^{q_1}\ssp
  x_{\sig_2}^{p_2} y_{\sig_2}^{q_2}\ssp \cdots \ssp  x_{\sig_n}^{p_n} y_{\sig_n}^{q_n}\ess .
\eqno 1.1
$$
Thus using the diagonal symmetry of the operator $D_{hk}$ we have
$$
D_{hk}\ssp \DD_L \ses \sum_{i=1}^n\ssp {1\over p!q!}
\sum_{\sig \in S_n}\ssp
\sign(\sig)\ssp  \del_{x_{\sig_i}}^h  \del_{y_{\sig_i}}^k\ssp
\bigl(\ssp x_{\sig_1}^{p_1} y_{\sig_1}^{q_1}\ssp  
x_{\sig_2}^{p_2} y_{\sig_2}^{q_2}\ssp \cdots \ssp  x_{\sig_n}^{p_n} y_{\sig_n}^{q_n}\ssp \bigr) \ess .
\eqno 1.2
$$
Now, 
$$\del_{x_{\sig_i}}^h \del_{y_{\sig_i}}^k 
\bigl(  x_{\sig_1}^{p_1} y_{\sig_1}^{q_1} 
 \cdots 
  x_{\sig_n}^{p_n} y_{\sig_n}^{q_n}  \bigr)
=
\cases{
(p_i)_h(q_i)_k\,
x_{\sig_1}^{p_1} y_{\sig_1}^{q_1}\,\cdots\,   
x_{\sig_i}^{p_i-h} y_{\sig_i}^{q_i-k} \,\cdots\,  x_{\sig_n}^{p_n} y_{\sig_n}^{q_n}
& if $h\leq p_i$ and $k\leq q_i\,$,\cr\cr
0 & otherwise,
} 
$$
where for two integers $h\leq p$ we set $\ssp (p)_h=p(p-1)\cdots (p-h+1)$. 
\sas

Moreover, we can easily see that the determinant
$$
\sum_{\sig \in S_n}\ssp
\sign(\sig)\ssp
x_{\sig_1}^{p_1} y_{\sig_1}^{q_1}\ssp\cdots   
x_{\sig_i}^{p_i-h} y_{\sig_i}^{q_i-k}\ssp \cdots \ssp  x_{\sig_n}^{p_n} y_{\sig_n}^{q_n} 
\eqno 1.3
$$
fails to vanish if and only if the biexponent pairs
$$
(p_1,q_1)\scs \ldots \scs (p_{i-1},q_{i-1})\scs(p_i-h,q_i-k)\scs(p_{i-1},q_{i-1})\scs\ldots \scs (p_n,q_n)
\eqno 1.4
$$
are all distinct. Putting all this together, formula I.5 follows from our conventions concerning
lattice determinants.
\sa
\heading{\bol Proof of Proposition I.2}

What we assert there is just a special case of Proposition I.1.
\sa

\heading{\bol Proof of Proposition I.3}

Note that from Proposition I.2 it immediately follows that
$$
\pm D_x^hD_y^k\ssp \DD_{\mu/ij}\ses \pm D_{hk}\ssp \DD_{\mu/ij}\ses \DD_{\mu/i+h,j+k}\ess ,
\eqno 1.5
$$
and this is easily seen to imply I.6 and I.7. To show the stated surjectivity, 
we use the nonsingularity of the $\flip$ map and write every element
$Q\in \BM_{\mu/i+h,j+k}$  in the form  $Q=P(\del_x,\del_y)\DD_{\mu/i+h,j+k}$
with $P(x,y)$ a uniquely determined element of  $\BM_{\mu/i+h,j+k}$. Now, 
we see from 1.5 that we also have
$$
Q\ses \pm D_x^hD_y^k\ssp P(\del_x,\del_y)\DD_{\mu/ij}\ses \pm  D_{hk}\ssp P(\del_x,\del_y)\DD_{\mu/ij}\ess .
$$
This shows that  both $ D_x^hD_y^k$ and $ D_{hk}$
map the subspace
$$
\{\ssp P(\del_x,\del_y)\DD_{\mu/ij}\ess :\ess P\in \BM_{\mu/i+h,j+k}\ess \}\ses 
\flip_{\DD_{\mu/i j }}\ess \BM_{\mu/i+h,j+k}\ess\con\ess \BM_{\mu/ij}
$$  
isomorphically onto $\BM_{\mu/i+h,j+k}$. This completes our proof.
\sa

\heading{\bol Remark 1.1}

We get a better picture of what is going on here if we make use of Proposition I.7.
For instance, if we let  $\BK_{ij}^{hk}$
denote the kernel of $D_{hk}$ as a map of $\BM_{ij}$ onto $\BM_{i+h,j+k}\,$,
then I.29 b), with $\DD=\DD_{\mu/ij}\ess $  and $\ess {\tilde \DD}=\DD_{\mu/i+h,j+k}\ssp $,
gives the direct sum decomposition
$$
\BM_{\mu/ij}\ses \flip_{\DD_{\mu/ij}}  \BM_{\mu/i+h,j+k}
\ess\ssp  \oplus \ess\ssp \BK_{ij}^{hk}\ess .  
\eqno 1.6
$$
Moreover, since $D_{hk}$, (up to a bidegree shift of $(-h,-k)$), gives also an isomorphism of 
bigraded $S_n$-modules of 
$\flip_{\DD_{\mu/ij}}  \BM_{\mu/i+h,j+k}$ onto $\BM_{\mu/i+h,j+k}\ess $,
we see from 1.6 that the bigraded Frobenius characteristic 
$K_{ij}^{hk}(x;q,t)$  of $\BK_{ij}^{hk}$  must be given by the formula
$$
K_{ij}^{hk}\ses C_{ij}(x;q,t)\sms t^hq^k \ssp C_{i+h,j+k}(x;q,t)\ess . 
\eqno 1.7
$$  
\sa

\heading{\bol Proof of Proposition I.4}
\sas

Our proof proceeds by induction with respect to the partial order $(i,j)\leq (i',j')$. 
We know from [10] that, up to a scalar factor, $\DD_\aaa(x;y)$ is the only alternant 
in $\BM_\aaa$. This can also be seen from the following reasoning.
Note that all the monomials
$$
x_1^{p_1}y_1^{q_1}x_2^{p_2}y_2^{q_2}\cdots x_n^{p_n}y_n^{q_n}
\eqno 1.8
$$
occurring in $\DD_\aaa(x;y)$ consist of factors $x_i^{p_i}y_i^{q_i}$
with  $(p_i,q_i)\in \aaa$. Since $\aaa$ has only $n$ cells, all
the monomials contained in 
any derivative of $\DD_\aaa(x;y)$ will have 
at least one pair of equal biexponents. This forces the vanishing of
the antisymmetrization of  every derivative of
$\DD_\aaa(x;y)$. This proves the assertion when $(i,j)$ is a corner cell of
$\mu$ and $\alpha=\mu/ij$. So let us assume that the assertion is true for any $(i',j') >(i,j)$. This given,
note that every bihomogeneous alternating polynomial $\DD(x;y)\in \BM_{\mu/ij}$ 
can be written in the form
$$
\DD(x;y)\ses P(\del_x;\del_y)\ssp \DD_{\mu/ij}(x;y)
\eqno 1.9
$$ 
with $P$ bihomogeneous and invariant under the diagonal action.
Now it is well known (see [24]) that the ideal generated by the diagonal invariant
polynomials with vanishing constant term is also generated by the polynomials
$$
\sum_{i=1}^n\ssp x_i^h\ssp y_i^k \bigsp \hbox{\ with $\ess 1\leq h+k\leq n$}\ess .
$$
Thus, if $P(x;y)$ is not a constant, we may express it in the form
$$
P(x;y)\ses \sum_{1\leq h+k\leq n}\ssp A_{hk}(x;y)\ssp \sum_{i=1}^n\ssp x_i^h\ssp y_i^k\ess .
\eqno 1.10
$$
Substituting this into 1.9 and using 1.5 gives
$$
\DD(x;y)=\hskip -.1in\sum_{1\leq h+k\leq n} \hskip -.08in  A_{hk}(\del_x;\del_y)\ssp D_{hk}\ssp \DD_{\mu/ij}(x;y)
=
\hskip -.1in
\sum_{\multi{(i+h,j+k)\ssp \in\ssp \mu\cr 1\leq h+k\leq n\cr }}  
\hskip -.08in \pm A_{hk}(\del_x;\del_y) \ssp \DD_{\mu/i+h,j+k}(x;y)\ess . 
\eqno 1.11
$$
Thus, from the induction hypothesis we derive that any bihomogeneous alternant of $\BM_{\mu/ij}\ssp $, 
with lesser total degree than $\DD_{\mu/ij}\ssp $, must be a linear combination of the $\DD_{\mu/i'j'}$
with $(i',j')>(i,j)$. This completes the induction since the only elements of $\BM_{\mu/ij}$
of the same total degree as $\DD_{\mu/ij}$ are its scalar  multiples.
\sa

\heading{\bol Proof of Proposition I.5 }

From Proposition I.1 it immediately follows that for any $\mu\part n+1$ we have
$$
\sum_{i=1}^{n+1}\ssp \del_{x_i}^h\del_{y_i}^k\ssp \DD_\mu(x;y)\ses 0\bigsp (\ssp \forall \ess\ess h+k\geq 1\ssp )
\ess .
\eqno 1.12
$$
In particular, if $D_x$ and $D_y$ are as given in I.6,  we deduce that 
$$
\eqalign{
&\del_{x_{n+1}}\ssp \DD_\mu(x;y)\ssp =\ssp - D_x\ssp \DD_\mu(x;y){\ \atop \ }\ess ,\cr 
&\del_{y_{n+1}}\ssp \DD_\mu(x;y)\ssp =\ssp - D_y\ssp \DD_\mu(x;y)\ess .\cr
}
\eqno 1.13
$$
This means that in constructing a basis for $\BM_\mu$ of the form
$$
\CB_\mu\ses \{\ssp b(\del_x;\del_y)\ssp \DD_\mu(x;y)\ssp :\ssp b\in \CC \ssp \}\ess ,
\eqno 1.14
$$
the polynomials in $\CC$ need not contain any of the variables $x_{n+1}\scs y_{n+1}$.
Now we have the following 
\ess

\heading{\bol Lemma 1.1}

{\ita If $\CC $ is a collection of polynomials in $\BQ[\xon\ssp ;\ssp \yon]$ then
the collection $\CB_\mu$ given in 1.14 is a basis for $\BM_\mu$ if and only if the collection
$$
\CB_{\mu/00}\ses  \{\ssp b(\del_x;\del_y)\ssp \DD_{\mu/00}(x;y)\ssp :\ssp b\in \CC \ssp \} 
\eqno 1.15 
$$ 
is a basis for $\BM_{\mu/00}\ess .$
}
\sas

\heading{\bol Proof}

The Laplace expansion of the determinant giving $\DD_\mu$, with respect to the last row, gives that
$$
\DD_\mu(x;y)\ses \sum_{(i,j)\in \mu} x_{n+1}^i y_{n+1}^j\ess \eee_{ij}\ssp \DD_{\mu/ij}(\xon;\yon) 
\eqno 1.16 
$$
with $\eee_{ij}=\pm 1$. Note then that for $f\in  \BQ[\xon\ssp ;\ssp \yon]$ we necessarily have
$$
a)\ess\ess f(\del_x;\del_y)\ssp \DD_\mu(x;y)\ses  0
\ess\ess\ess \longleftrightarrow \ess\ess\ess
b)\ess\ess f(\del_x;\del_y)\ssp \DD_{\mu/00}(x;y)\ses  0\ess .
$$
In fact, we see from 1.16 that b) immediately follows from a) by setting $x_{n+1}=y_{n+1}=0$.
Conversely, if b) holds true then by applying to it the operator $D_{i,j}$ we obtain that
$$
f(\del_x;\del_y)\ssp \DD_{\mu/ij}(x;y)\ses  0
$$
must hold as well for all $(i,j)\in \mu$ and then a) again follows by applying $f(\del_x;\del_y)$ to
both sides of 1.16. We thus derive that, for a given collection $\CC$, $\ess \CB_\mu$ is an independent set    
if and only if $\CB_{\mu/00}$ is. In particular, both spaces $\BM_\mu$ and $\BM_{\mu/00}$ must have
the same dimension.
\bigsp {\bf Q.E.D.}
\sap 

This given, I.11 follows by choosing $\CC$ so that both $\CB_{\mu}$ and
$\CB_{\mu/00}$ are  bihomogeneous bases and noting that (because of Lemma 1.1) the action of
$S_n$ on corresponding bihomogeneous components of $\CB_{\mu}$ and
$\CB_{\mu/00}$  are given by the same matrices. This completes the proof of 
Proposition I.5.
\sa

We should note that a useful consequence of Lemma 1.1 is the following.
\sas

\heading{\bol Proposition 1.1}

{\ita If $\CB_\mu^*(x_1,\ldots,x_n,x_{n+1};y_1,\ldots,y_n,y_{n+1})$ is a basis for
$\BM_\mu$ then $\CB_\mu^*(x_1,\ldots,x_n,0;y_1,\ldots,y_n,0)$  is
a basis for $\BM_{\mu/00}\ess .$}
\sas

\heading{\bol Proof}

Let $\CC\con {\bf Q}[\xon;\yon]$ be chosen so that both $\CB_\mu$ and  $\CB_{\mu/00}$, 
(as given by 1.14 and 1.15) are bases for
$\BM_\mu$ and  $\BM_{\mu/00}$ respectively. By assumption, for every element of $b\in \CC$
we have the expansion
$$
\eqalign{
b(\del_{x_1},\ldots ,\del_{x_n};\del_{y_1},\ldots ,\del_{y_n})\ssp \DD_\mu(x_1,&\ldots ,x_n,x_{n+1};y_1,\ldots ,y_n,y_{n+1})
  \cr
& \ses \sum_{b^*\in \CB_\mu^*} c_{b^*}\ssp b^*(x_1,\ldots ,x_n,x_{n+1};y_1,\ldots ,y_n,y_{n+1})\ess .
\cr
}
$$
However, setting $x_{n+1}=y_{n+1}=0$  here (and using 1.16) gives the identity
$$
\eqalign{
b(\del_{x_1},\ldots ,\del_{x_n};\del_{y_1},\ldots ,\del_{y_n})\ssp \DD_{\mu/00}(x_1,&\ldots ,x_n;y_1,\ldots ,y_n)
  \cr
& \ses \sum_{b^*\in \CB_\mu^*} c_{b^*}\ssp b^*(x_1,\ldots ,x_n,0;y_1,\ldots ,y_n,0)\ess .
\cr
}
$$
This shows that $\CB_\mu^*(x_1,\ldots ,x_n,0;y_1,\ldots ,y_n,0)$  spans $\BM_{\mu/00}$.
However, it must be a basis since its cardinality is no larger than the dimension of $\BM_\mu$
and the latter has the same dimension as $\BM_{\mu/00}$.  
\sa

\heading{\bol Proof of Theorem I.1 }

For a given cell $(i,j)\in \mu$ we are to determine if there are constants $x$, $y$ and $z$ such that 
$$
C_{\mu/ij}\sms x\ssp  C_{\mu/i,j+1}\sms y\ssp  C_{\mu/i+1,j} \sps z\ssp  C_{\mu/i+1,j+1}\ses 0 \ess .
\eqno 1.17 
$$
Let us begin with the generic case, that is when the shadows of the four cells 
$\ess (i,j)$, $(i,j+1)$, $(i+1,j)$, $(i+1,j+1)$
contain the same corners of $\mu$. To this end,  
let $\tau$ be the partition contained in the shadow of $(i,j)$ and $\rho$ be one of the
predecessors of $\tau$. 
Denoting by $c_\rho^{ij}$, $c_\rho^{i,j+1}$, $c_\rho^{i+1,j}$ and  $c_\rho^{i+1,j+1}$
the coefficients of $\TH_{\mu-\tau+\rho}(x;q,t)$ in $C_{\mu/ij}(x;q,t)$, 
$C_{\mu/i,j+1}(x;q,t)$, $C_{\mu/i+1,j}(x;q,t)$ and $C_{\mu/i+1,j+1}(x;q,t)$ respectively, 
it is not difficult to derive from I.16 and the definition I.14 that we must have
$$
\eqalign{
&c_\rho^{ij}\ses {t^{l_1}-q^{a_1+1}\over t^{l_1}-q^{a_1 } }\ess {q^{a_2}-t^{l_2+1}\over q^{a_2}-t^{l_2} }
\ess c_\rho^{i+1,j+1}\ess ,\cr 
&c_\rho^{i,j+1}\ses {q^{a_2}-t^{l_2+1}\over q^{a_2}-t^{l_2} } 
\ess c_\rho^{i+1,j+1}\ess ,\cr
&c_\rho^{i+1,j}\ses  {t^{l_1}-q^{a_1+1}\over t^{l_1}-q^{a_1 } }
\ess c_\rho^{i+1,j+1}\ess ,
}
$$
with
$$
l_1=l+i-l'\ess\scs\ess\ess 
a_1=a'-j\ess\scs\ess\ess 
l_2=l'-i\ess\scs\ess\ess 
a_2=a+j-a'\ess , 
\eqno 1.18 
$$
where $l$ and $a$ give the leg and arm of $(i,j)$ and 
$l'$ and $a'$ give the
coleg and coarm of the cell $\mu/\mu-\tau+\rho$.
\sa

This given, equating to zero the coefficient of  $\TH_{\mu-\tau+\rho}(x;q,t)$ in 1.17 yields the
equation 
$$
\biggl({t^{l_1}-q^{a_1+1}\over t^{l_1}-q^{a_1 } }\ess {q^{a_2}-t^{l_2+1}\over q^{a_2}-t^{l_2} }
\sms x\ess {q^{a_2}-t^{l_2+1}\over q^{a_2}-t^{l_2} } 
\sms y\ess {t^{l_1}-q^{a_1+1}\over t^{l_1}-q^{a_1 } }
\sps z\ess\biggr)\ess c_\rho^{i+1,j+1}\ses 0
$$
Since by definition the coefficients $c_{\mu\nu}$ are never zero, we see that 
1.17 will hold true if and only if we can find $x,y$ and $z$ {\bol independent} of $\rho$ 
such that
$$
(t_1 - q\ssp  q_1) (q_2 - t\ssp t_2) - x\ssp  (q_2 - t\ssp  t_2) (t_1 - q_1)
                - y\ssp  (t_1 - q\ssp  q_1) (q_2 - t_2) + z\ssp  (t_1 - q_1) (q_2 - t_2)= 0 \, ,
\eqno 1.19 
$$
where for convenience we have set
$$
t_1=t^{l_1}\scs\ess\ess t_2=t^{l_2}\scs\ess\ess q_1=q^{a_1}\scs\ess\ess q_2=q^{a_2} \ess .
$$
Setting $T=t^l$ and $Q=q^a$, from 1.18 we deduce that
 $t_2=T/t_1\scs q_2=Q/q_1$. Thus,  
making these substitutions and multiplying by $\ssp t_1\ssp q_1\ssp $, reduces 1.19  to
$$                                  
        Q\ssp (  x  +  y - z  - 1 )\, t_1^2\sms 
\bigl( x (t T +  Q) + y ( T + q Q) - z (T + Q) - (t T + q Q)\bigr)\, q_1 t_1\sps
T\ssp  (x t   + y q -z   - t   q   )\,  q_1^2 \ses 0\ess .
$$
Now this is most fortunate since the coefficients of $t_1^2$, $t_1q_1$ and $q_1^2$ 
are independent of $\rho$. 

Setting to zero these coefficients yields the system    
$$
\matrix{
\hfill   x & + & \hfill  y &-& \hfill  z & = &1 \cr
   (t T +  Q)\ssp x& +&  ( T + q Q)\ssp y& - &  (T + Q)\ssp z& = & t T + q Q \cr
\hfill  t\ssp x  & + &\hfill q\ssp y& -&\hfill z &  =& t   q 
\cr}
$$
whose unique solution
$$
x\ses {T-q\ssp Q\over T-Q }\ess \scs \ess\ess\ess 
y\ses {t \ssp T-  Q\over T-Q  } \ess \scs \ess\ess\ess 
z\ses{t\ssp T-q\ssp Q\over T-Q }
$$
establishes the identity in I.17, in this case. 

Let us deal next with the case when the leftmost corner of $\tau$
is in the shadow of $(i+1,j)$ but not in the shadow of $(i,j+1)$ and $(i+1,j+1)$. 
Let $\rho_1$ be the partition obtained by removing this corner from $\tau$.
This given we derive from  I.16 that neither $C_{\mu/i,j+1}(x;q,t)$ nor $C_{\mu/i+1,j+1}(x;q,t)$ 
will contain a term involving  $\TH_{\mu-\tau+\rho_1}(x;q,t)$ in their expansion.
So, taking the coefficient of this polynomial in 1.17 reduces it to
$$
c_{\rho_1}^{ij}\sms y\ssp c_{\rho_1}^{i+1,j}\ses  0\ess .
$$ 
Now, using again the same notation, we may write 
$$
c_{\rho_1}^{ij}\ses {q^{a_2} -t^{l_2+1}\over q^{a_2} -t^{l_2} }\ess c_{\rho_1}^{i+1,j}\ess .
$$
These two equations give that
$$
y\ses {q^{a_2} -t^{l_2+1}\over q^{a_2} -t^{l_2} } \ess .
$$
However, in this case it is easily seen that $l_2=l'-i=l$ and $a'=j$,
giving $a_2=a+j-a'=a$, and we are led again to the solution
$$
y\ses    {Q -t \ssp T  \over Q  -T }\ess . 
$$
The remaining cases can be easily checked to yield the same values of $x$ and $y$. 
This completes 
the proof of Theorem I.1 
since the other assertions are immediate consequences of I.16.
\sa

\heading{\bol Proof of Proposition I.6  }

There is very little left to do here since (see Remark 1.1)
both equations in I.19 a) are but particular cases of 1.7
and the equations in I.19 b) as well as I.20 are immediate consequences of the definitions. 
\sa

\heading{\bol Proof of Proposition I.7  }

By the properties of the map $\flip_\TDD\ssp $, a polynomial $Q(x,y)$  in $\BTM$  may be written in the form
$$
Q(x,y)= P(\del_x;\del_y)\TDD(x;y)\ess ,\ess\ess\ess\ess\ess \ess {\rm with}\ess\ess P(x;y)\in \BM_\TDD\ess .
$$
Since $\TDD(x;y)= D(\del_x;\del_y)\DD(x;y)\,$, we may also write $Q(x,y)$ in the form
$$
Q(x,y)= P(\del_x;\del_y) D(\del_x;\del_y) \DD(x;y)= D(\del_x;\del_y) P(\del_x;\del_y)\DD(x;y)
$$
with
$$
 P(\del_x;\del_y)\DD(x;y)\in \BM_\DD\ess .
$$
This establishes surjectivity and the containment $\BTM\subseteq\BM_\DD$. In fact,
this argument shows that $D(\del_x;\del_y)$ maps the space
$$
\flip_\DD\ess \BM_\TDD\ses \bigl\{\ess P(\del_x;\del_y)\ssp \DD\ess :\ess P\in \BM_\TDD\ess \bigr\}
$$
surjectively onto $\BM_\TDD$.

Now we establish I.28, the description of the kernel. To this end note that
the polynomial $f=P(\del_x;\del_y)\DD(x;y)=\flip_\Delta P$  is
in $\BK$ if and only if
$0=D(\del_x;\del_y)f(x;y)=P(\del_x;\del_y)\TDD(x;y)$, or equivalently,
$P\in\BTMperp$. Thus we may write
$$
\eqalign{
\BK & \ses \bigl\{\ssp f=P(\del_x;\del_y)\DD(x;y)\ssp :
\ssp P\in \BM_\DD\ess\ess \&\ess\ess P(\del_x;\del_y)\TDD(x;y)=0\ssp \bigr\}  \cr
&\ses \flip_\DD\ssp \bigl\{\ssp  P\in \BM_\DD \ssp :\ssp P(\del_x;\del_y)\TDD(x;y)=0\ssp \bigr\}  \cr
&\ses \flip_\DD\ssp  \BM_\DD   \cap  \BTMperp \ess , \cr
}
$$
and I.28 follows by an application of $\flip_\DD^{-1}$ to both sides of this relation.
This shows that the orthogonal decomposition 
$$
\BM_\DD\ses \BM_\TDD\ess \oplus_\perp \ess \BM_\DD   \cap  \BTMperp 
$$
in this case can be written in the form
$$
\BM_\DD\ses \BM_\TDD\ess \oplus_\perp\ess \flip_\DD^{-1}\ess \BK\ess , 
$$
establishing I.29 a). Applying $\flip_\DD$ to both sides gives I.29 b), completing our proof. 
\sa

\heading{\bol Proof of Proposition I.8  }

Our point of departure is formula I.16. So let $\tau$ be the partition in the shadow of $(i,j)$
and let $\ssp x_0^{ij},\ldots ,x_m^{ij}\ess ;\ess u_0^{ij},\ldots ,u_m^{ij}$ be the corner weights of $\tau$.
Let  $\rho^{(1)},\rho^{(2)},\ldots ,\rho^{(m)}$ be the predecessors of $\tau$ ordered from
left to right so that $\ssp x_1^{ij},\ldots ,x_m^{ij}$ are the respective weights of the
cells $\tau/\rho^{(1)},\ldots ,\tau/\rho^{(m)}$.
This given, using formula I.44 with $\mu$ replaced by $\tau$ and $\nu^{(i)}$ replaced by $\rho^{(s)}$,
formula I.16 becomes
$$
C_{\mu/ij}(x;q,t)\ses {1\over M}\ssp \sum_{s=1}^m\ssp {1\over x_s^{ij}}
\ssp {\prod_{r=0}^m \bigl(  x_s^{ij}- u_r^{ij}\bigr)\over \prod_{r=1\,;\,r\neq s}^m \bigl(   x_s^{ij}- x_r^{ij}\bigr) }
\ess \TH_{\mu-\tau+\rho^{(s)}}\ess .
\eqno 1.20 
$$
 For convenience set $ \mu-\tau+\rho^{(s)}=\aaa^{(s)} $,
so that as in I.48 we have
$$ 
S\ses \Pred_{ij}(\mu)\ses \bigl\{\aaa^{(1)},\aaa^{(2)},\ldots ,\aaa^{(m)} \bigr\}\ess .
$$

\noindent
Now, formula I.40 for $\aaa=\aaa^{(s)}$ may be written as
$$
\TH_{\aaa^{(s)}}(x;q,t)\ses 
\prod_{r=1\,;\,r\neq s}^m  \Bigl( 1- {\nabla\over T_{\aaa^{(r)}} }\Bigr)\ssp \phi_S^{(m)}\ess .
\eqno 1.21 
$$
Note next that from the definition of $\;\mu/ij\;$ it follows that $T_\mu=t^iq^jT_{\mu/ij}\,$, and
since $t^iq^j x_s^{ij}$ is the weight of the cell $\mu/\mu-\tau+\rho^{(s)}=\mu/\aaa^{(s)}\,$,
we also have $t^iq^j x_s^{ij}T_{\aaa^{(s)}}=T_\mu$. In conclusion we see  that
$$
{1\over T_{\aaa^{(s)}} }\ses { x_s^{ij}\over T_{\mu/ij} }\ess .
\eqno 1.22 
$$
Using this in 1.21 and substituting the resulting expression in 1.20 we finally obtain  
$$
C_{\mu/ij}(x;q,t)\ses {1\over M}\ssp \sum_{s=1}^m\ssp {1\over x_s^{ij}}
\ssp {\prod_{r=0}^m \bigl(  x_s^{ij}- u_r^{ij}\bigr)\over \prod_{r=1\,;\,r\neq s}^m \bigl(   x_s^{ij}- x_r^{ij}\bigr) }
\ess \prod_{r=1\,;\,r\neq s}^m 
 \Bigl( 1- \nabla { x_r^{ij}\over T_{\mu/ij} }\Bigr)\ssp \phi_S^{(m)}\ess . 
\eqno 1.23 
$$ 
Now it develops that we have the following identity.
\sas

\heading{\bol Lemma 1.2}

{\ita If $x_0,x_1,\ldots ,x_m$ and $u_0,u_1,\ldots ,u_m$ are any quantities such that 
$$
x_0  x_1 \cdots  x_m=u_0  u_1 \cdots  u_m\ess ,
\eqno 1.24 
$$
then for all $z$ we have} 
$$
\sum_{s=0}^m{1\over x_s}
\ssp  
{
\prod_{r=0}^m (x_s-u_r) 
\over
\prod_{r=1\,;\,r\neq s}^m (x_s-x_r) 
}
\ssp \prod_{r=1\,;\,r\neq s}^m (1-z x_r)
\ses
{1\over z}\ssp \biggl(\prod_{s=0}^m \bigl(1-z u_s)\sms \prod_{s=0}^m \bigl(1-z x_s\bigr)\biggr)\ess .
\eqno 1.25 
$$
{\bol Proof}

Note that because of 1.24 the expression on the right hand side of 1.25 evaluates to a polynomial of degree
at most $m-1$. We can thus apply the Lagrange interpolation formula at the points 
$$
1/x_1\scs 1/x_2\scs \ldots \scs 1/x_m\scs
$$
and obtain that
$$
{1\over z}\ssp \biggl(\prod_{s=0}^m \bigl(1-z u_s\bigr)\sms \prod_{s=0}^m \bigl(1-z x_s\bigr)\biggr)\ses 
\sum_{s=0}^m{  x_s}
\ssp 
\prod_{r=0}^m \bigl(1-{ u_r\over x_s}\bigr) 
\ssp \prod_{r=1\,;\,r\neq s}^m {(z-{1\over  x_r})\over({1\over  x_s}-{1\over  x_r}) }\ess .
$$
Clearly this is just another way of writing 1.25.
\sa

This given, evaluating 1.25 at
$$
x_s=x_s^{ij}\ess \scs\ess\ess u_s=u_s^{ij}\ess {\rm and} \ess\ess z={\nabla\over T_{\mu/ij}}\ess ,
$$ 
applying both sides to $\phi_S^{(m)}$ and using 1.23, we finally obtain that
$$
C_{\mu/ij}\ses {1\over M}{T_{\mu/ij}\over \nabla}\ess 
\biggl[\ssp \prod_{s=0}^m\ssp \Bigl(1-\nabla {u_s^{ij}\over T_{\mu/ij}}\ssp \Bigr)
-
\ssp \prod_{s=0}^m\ssp \Bigl(1-\nabla {x_s^{ij}\over T_{\mu/ij}}\ssp \Bigr)
 \biggr]\ssp \phi_S^{(m)}\ess  .
\eqno 1.26 
$$
We claim that this formula contains both I.49 a) and b).  In fact,
expanding the products, collecting terms according to powers of $\nabla$, 
and using the identity
$$
u_0^{ij}u_1^{ij}\cdots u_m^{ij}\ses x_0^{ij}x_1^{ij}\cdots x_m^{ij}
\eqno 1.27 
$$
gives
$$
C_{\mu/ij}\ses
\sum_{k=1}^m\ssp {(-\nabla)^{m-k}\phi_S^{(m)}\over T_{\mu/ij}^{m-k}}\ssp 
{ e_{m+1-k} [x_0^{ij}+\cdots +x_d^{ij}]\sms e_{m+1-k} [u_0^{ij}+\cdots +u_d^{ij}]\over M}
\ess .
$$
In view of I.35, we see that this just another  way of writing I.49 b).
Note next that, using 1.22, we may write 
$$
\prod_{s=0}^m\ssp \Bigl(1-\nabla {x_s^{ij}\over T_{\mu/ij}}\ssp \Bigr)
\ses\Bigl(1-\nabla {x_0^{ij}\over T_{\mu/ij}}\ssp \Bigr)
\prod_{s=1}^m\ssp \Bigl(1- {\nabla\over T_{\aaa^{(s)}}}\ssp \Bigr)\ess .
$$ 
However, using I.36, we derive that
$$
\prod_{s=1}^m\ssp \Bigl(1- {\nabla \over T_{\aaa^{(s)}}}\ssp \Bigr)\phi_S^{(m)}
\ses
\sum_{\aaa\in S}\ssp 
\Bigl(\prod_{\bbb\in S/\{\aaa\}}{1\over 1-T_\aaa/T_\bbb}\ssp \Bigr)
\prod_{s=1}^m\ssp \Bigl(1- {T_\aaa\over T_{\aaa^{(s)}}}\ssp \Bigr)\ess\TH_\aaa
\ses  0\ess  .
\eqno 1.28
$$
Thus the second product in 1.26 is entirely superfluous and we see that
1.26 is also another way of writing I.49 a). This
completes the proof of Proposition I.8.  

\heading{\bol Proof of Theorem I.2}

We see from the recurrences in I.19 that we may write
$$
\eqalign{
&a)\ess\ess A_{ij}^x\ses C_{\mu/ij}\sms t\ssp C_{\mu/i+1,j}\sms C_{\mu/i,j+1}
\sps t\ssp C_{\mu/i+1,j+1}\ess ,{\ \atop\  }\cr 
&b)\ess\ess A_{ij}^y\ses C_{\mu/ij}\sms q\ssp C_{\mu/i,j+1}\sms  C_{\mu/i+1,j}
\sps q\ssp C_{\mu/i+1,j+1}\ess .\cr
}
\eqno 1.29
$$
Let us, for a moment, assume (as in the proof of Theorem I.1) that
the shadows of the four cells $(i,j)$, $(i+1,j)$, $(i,j+1)$, $(i+1,j+1)$
contain the same corners of $\mu$. This means that if
$$
\Pred_{i,j}(\mu)\ses  S\ses\{\aaa^{(1)},\aaa^{(2)},\ldots ,\aaa^{(m)}\}
$$
then we also have
$$
S\ses \Pred_{i+1,j}(\mu)\ses\Pred_{i,j+1}(\mu)\ses\Pred_{i+1,j+1}(\mu)\ess .
$$
Thus formula I.49 a) applied to each of the four cells gives
$$
\eqalign{
C_{\mu/ij}& \ses {1\over M}{T_{\mu/ij}\over \nabla}\ess 
\biggl(\ssp \prod_{s=0}^m\ssp \Bigl(1-\nabla {u_s^{ij}\over T_{\mu/ij}}\ssp \Bigr)
 \biggr)\ssp \phi_S^{(m)}\ssp ,\cr
C_{\mu/i+1,j}& \ses {1\over M}{T_{\mu/i+1,j}\over \nabla}\ess 
\biggl(\ssp \prod_{s=0}^m\ssp \Bigl(1-\nabla {u_s^{i+1,j}\over T_{\mu/i+1,j}}\ssp \Bigr)
 \biggr)\ssp \phi_S^{(m)}\ssp ,\cr
C_{\mu/i,j+1}& \ses {1\over M}{T_{\mu/i,j+1}\over \nabla}\ess 
\biggl(\ssp \prod_{s=0}^m\ssp \Bigl(1-\nabla {u_s^{i,j+1}\over T_{\mu/i,j+1}}\ssp \Bigr)
 \biggr)\ssp \phi_S^{(m)}\ssp ,\cr
C_{\mu/i+1,j+1}& \ses {1\over M}{T_{\mu/i+1,j+1}\over \nabla}\ess 
\biggl(\ssp \prod_{s=0}^m\ssp \Bigl(1-\nabla {u_s^{i+1,j+1}\over T_{\mu/i+1,j+1}}\ssp \Bigr)
 \biggr)\ssp \phi_S^{(m)}\ssp .\cr
}
\eqno 1.30
$$
In the figure below we have depicted  the generic situation we are dealing with.
$$
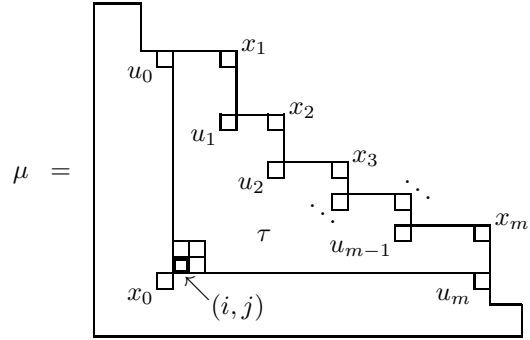
\hbox{%
	\myput(0,200){\llap{$\mu\;\;=\;\;\;$}}
	\hseg(0,0;540)					
	\vseg(540,0;40)
	\hseg(500,40;40)
	\corn(500,80)\labLL(500,80;u_m)
	\drawL(400,500,40,140;x_m,u_{m-1})
	\drawL(320,400,140,180;,)
		\myput(390,170){$\ddots$}
		\myput(270,140){$\ddots$}
	\drawL(240,320,180,220;x_3,u_2)
	\drawL(180,240,220,280;x_2,u_1)
	\drawL(100,180,280,360;x_1,u_0)
	\corn(100,80)\labLL(100,80;x_0)
	\hseg(60,360;40)
	\vseg(60,360;60)
	\hseg(0,420;60)
	\vseg(0,0;420)
%
	\hseg(100,80;400)
	\vseg(100,80;280)
	%
	\bcorn(120,100)\corn(140,100)\corn(120,120)\corn(140,120)
	\labUR(200,120;\tau)
	\myput(115,75){\vbox to0pt{\hbox{$\raise.5\baselineskip\hbox{$\nwarrow$}(i,j)$}\vss}}
}
\hskip 540\mylength
$$

We have the partition $\tau$ that is in the shadow of $(i,j)$, its
corner cells, the corresponding corner weights, the cell $(i,j)$ and the adjacent
cells $(i+1,j)$, $(i,j+1)$, $(i+1,j+1)$.
Now a look at the figure should reveal that in this case we have the identities
$$
{u_s^{ij}\over T_{\mu/ij}}\ses 
{u_s^{i+1,j}\over T_{\mu/i+1,j}}\ses 
{u_s^{i,j+1}\over T_{\mu/i,j+1}}\ses 
{u_s^{i+1,j+1}\over T_{\mu/i+1,j+1}}
\bigsp \hbox{( for  $\ess  1\leq s\leq m-1$ )} \ess .
$$ 
This implies that $C_{\mu/ij}$, $C_{\mu/i+1,j}$, $C_{\mu/i,j+1}$ and $C_{\mu/i+1,j+1}$
have the common factor
$$
CF\ses {1\over M \ssp \nabla}\ess \prod_{s=1}^{m-1}\ssp 
\Bigl(1- {u_s^{ij}\over T_{\mu/ij}}\nabla\ssp \Bigr)\ess .
\eqno 1.31
$$
Note further that, from the figure and the definition of the corner weights,
we see that we must also have 
$$
\matrix
{
 u_0^{ij} &=& t\ssp u_0^{i+1,j}    &=&    u_0^{i,j+1}  &=&  t\ssp u_0^{i+1,j+1} \ess ,\cr
&&&&\cr
 u_m^{ij} &=&    u_m^{i+1,j} &=&   q\ssp u_m^{i,j+1}   &=&  q\ssp u_m^{i+1,j+1}\ess ,\cr
&&&&\cr
T_{\mu/ij} &=&   t\ssp T_{\mu/i+1,j} &=&  q\ssp T_{\mu/i,j+1}   &=& t\ssp q\ssp T_{\mu/i+1,j+1}\ess . \cr
} 
$$
Thus, setting for convenience
$$
z_0^{ij}\ses {u_0^{ij}\over T_{\mu/ij}}\ess \nabla\ess\ess \ess\ess 
\hbox{ and}\ess\ess\ess\ess 
 z_m^{ij}\ses {u_m^{ij}\over T_{\mu/ij}}\ess \nabla\ess ,
$$
we can  rewrite the identities in 1.30 in the form
$$
\eqalign{
C_{\mu/ij}& \ses CF\ess {T_{\mu/ij}}\ess 
  \Bigl(1-{z_0^{ij}}\ssp \Bigr)\Bigl(1-{z_m^{ij}}\ssp \Bigr)
\ssp \phi_S^{(m)}\ssp ,\cr
C_{\mu/i+1,j}& \ses  CF\ess {T_{\mu/ij}\over t}\ess 
  \Bigl(1-{z_0^{ij}}\ssp \Bigr)\Bigl(1-t\ssp {z_m^{ij}}\ssp \Bigr)
\ssp \phi_S^{(m)}\ssp ,\cr
C_{\mu/i,j+1}& \ses  CF\ess {T_{\mu/ij}\over q}\ess 
  \Bigl(1-q\ssp {z_0^{ij}}\ssp \Bigr)\Bigl(1- {z_m^{ij}}\ssp \Bigr)
\ssp \phi_S^{(m)}\ssp ,\cr
C_{\mu/i+1,j+1}& \ses  CF\ess {T_{\mu/ij}\over tq}\ess 
  \Bigl(1-q\ssp {z_0^{ij}}\ssp \Bigr)\Bigl(1-  t\ssp {z_m^{ij}}\ssp \Bigr)
\ssp \phi_S^{(m)}\ssp .\cr
}
\eqno 1.32
$$
Substituting these expressions in 1.29 a) and grouping terms we get
$$
\eqalign{
A_{ij}^x&\ses CF\ssp\cdot\ssp T_{\mu/ij}\ssp
\Bigl[ (1-z_0^{ij})\bigl(\ssp 1-z_m^{ij}-1+t\ssp z_m^{ij}\ssp\bigr)
\sms
{1\over q}\ssp(1-q\ssp z_0^{ij})\bigl(\ssp 1-z_m^{ij}-1+t\ssp z_m^{ij}\ssp\bigr)\ssp \Bigr] 
\ssp \phi_S^{(m)} \cr
&\ses CF\ssp\cdot\ssp T_{\mu/ij}\ssp
\Bigl[ (1-z_0^{ij}) \ssp z_m^{ij}\ssp (t-1)
\sms
{1\over q}\ssp(1-q\ssp z_0^{ij})  \ssp z_m^{ij}\ssp (t-1)\ssp \Bigr] 
\ssp \phi_S^{(m)} \cr
&\ses CF\ssp\cdot\ssp T_{\mu/ij}\ssp z_m^{ij}\ssp (t-1)\ssp
\bigl(1-z_0^{ij}   
\sms
{1\over q} +  z_0^{ij} \bigr) \ses CF\ssp\cdot\ssp q^a\ssp \nabla\ssp (1-1/t)\ssp
(1- 1/q )\ssp \phi_S^{(m)} \ess , 
\cr
}
$$
where the last equality is due to the fact that we have $u_m^{ij}=q^a/t$ with $a$ the
arm of the cell $(i,j)$. Using 1.31 yields our desired formula
$$
A_{ij}^x\ses q^a \ess \prod_{s=1}^{m-1}\ssp \Big(1-{u_s^{ij}\over T_{\mu/ij}}\ssp\nabla\Big)
\ssp \phi_S^{(m)}\ess .
\eqno 1.33
$$
Similarly, starting from 1.29 b) we derive that
$$
\eqalign{
A_{ij}^y&\ses CF\ssp\cdot\ssp T_{\mu/ij}\ssp
\Bigl[ (1-z_m^{ij})\bigl(\ssp 1-z_0^{ij}-1+q\ssp z_0^{ij}\ssp\bigr)
\sms
{1\over t}\ssp(1-t\ssp z_m^{ij})\bigl(\ssp 1-z_0^{ij}-1+q\ssp z_0^{ij}\ssp\bigr)\ssp \Bigr] 
\ssp \phi_S^{(m)} \cr
&\ses CF\ssp\cdot\ssp T_{\mu/ij}\ssp
\Bigl[1-z_m^{ij}-{1\over t}+z_m^{ij}\ssp\Bigr] \ssp z_0^{ij}\ssp (q-1)
\ssp \phi_S^{(m)} \cr
&\ses CF\ssp\cdot\ssp T_{\mu/ij}\ssp(1-1/t)\ess
 z_0^{ij}\ssp (q-1)\ssp 
 \ses CF\ssp\cdot\ssp t^l\ssp \nabla\ssp (1-1/t)\ssp
(1- 1/q )\ssp \phi_S^{(m)} \ess ,  
\cr
}
$$
where we have set $u_0^{ij}=t^l/q$ with $l$ the leg of $(i,j)$.
This gives
$$
{A_{ij}^x\over q^a}\ses {A_{ij}^y\over t^l}
\eqno 1.34
$$
as desired. 
\sas

Let us assume next that the shadows of $(i,j)$ and $(i+1,j)$ contain the same corners of
$\mu$ with 
$$
\Pred_{i,j}(\mu)\ses 
\Pred_{i+1,j}(\mu)\ses S\ses\{\aaa^{(1)},\aaa^{(2)},\ldots ,\aaa^{(m)}\}
$$
but (see figure below) the shadows of $(i,j+1)$ and $(i+1,j+1)$ miss the corner
$\mu/\aaa^{(1)}$. Thus
$$
\Pred_{i,j+1}(\mu)\ses 
\Pred_{i+1,j+1}(\mu)\ses S^*\ses\{\aaa^{(2)},\ldots ,\aaa^{(m)}\}\ses
\ess S\ssp /\ssp\{\aaa_1\}\ess .
\eqno 1.35
$$
$$
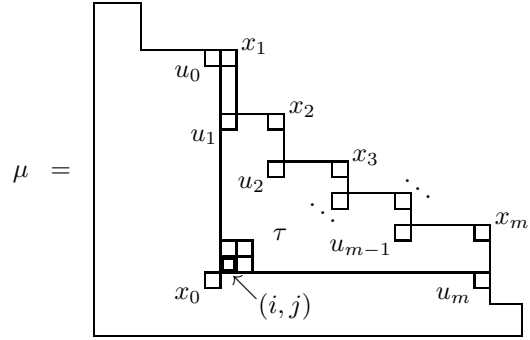
\hbox{%
	\myput(0,200){\llap{$\mu\;\;=\;\;\;$}}
	\hseg(0,0;540)					
	\vseg(540,0;40)
	\hseg(500,40;40)
	\corn(500,80)\labLL(500,80;u_m)
	\drawL(400,500,40,140;x_m,u_{m-1})
	\drawL(320,400,140,180;,)
		\myput(390,170){$\ddots$}
		\myput(270,140){$\ddots$}
	\drawL(240,320,180,220;x_3,u_2)
	\drawL(180,240,220,280;x_2,u_1)
	\drawL(160,180,280,360;x_1,u_0)
	\corn(160,80)\labLL(160,80;x_0)
	\hseg(60,360;100)
	\vseg(60,360;60)
	\hseg(0,420;60)
	\vseg(0,0;420)
%
	\hseg(160,80;340)
	\vseg(160,80;280)
	%
	\bcorn(180,100)\corn(200,100)\corn(180,120)\corn(200,120)
	\labUR(220,120;\tau)
	\myput(175,75){\vbox to0pt{\hbox{$\raise.5\baselineskip\hbox{$\nwarrow$}(i,j)$}\vss}}
}
\hskip 540\mylength
$$
Remarkably it develops that all the relations in 1.32 do hold true also in this
case so that the final conclusions in 1.33 and 1.34 still  remain unchanged.
To see how this comes about note first that,
since the situation is the same as before  as far as $(i,j)$ and
$(i+1,j)$ are concerned, there is no change in the first two 
equations of 1.30 and 1.32. On the other hand, in this case, 
the remaining two equations in 1.30  become
$$
\eqalign{a)\ess\ess\ess
C_{\mu/i,j+1}& \ses {1\over M}{T_{\mu/i,j+1}\over \nabla}\ess 
\biggl(\ssp \prod_{s=0}^{m-1}\ssp \Bigl(1-\nabla {u_s^{i,j+1}\over T_{\mu/i,j+1}}\ssp \Bigr)
 \biggr)\ssp \phi_{S^*}^{(m-1)}\ssp ,\cr
b)\ess\ssp C_{\mu/i+1,j+1}& \ses {1\over M}{T_{\mu/i+1,j+1}\over \nabla}\ess 
\biggl(\ssp \prod_{s=0}^{m-1}\ssp \Bigl(1-\nabla {u_s^{i+1,j+1}\over T_{\mu/i+1,j+1}}\ssp \Bigr)
 \biggr)\ssp \phi_{S^*}^{(m-1)}\ssp .\cr
}
\eqno 1.36
$$
Now, using 1.35, from I.37 we get
$$
\phi_{S^*}^{(m-1)}\ses \Bigl(  1-{\nabla\over T_{\aaa^{(1)} }} \Bigr)
\ssp \phi_{S}^{(m)}\ess .
$$
However, since $\ssp T_{\aaa^{(1)}}\ssp x_1^{ij}= T_{\mu/ij}\ssp $ and in this case
$x_1^{ij}=q\ssp u_0^{ij} $, this may be rewritten in the form
$$
\phi_{S^*}^{(m-1)}\ses \Bigl(  1-{q\ssp u_0^{ij}\nabla\over T_{\mu/ij} } \Bigr)
\ssp \phi_{S}^{(m)}\ess .
\eqno 1.37
$$
Note further that we also have 
$$
u_s^{i,j+1}\ses {1\over q}\ess u_{s+1}^{ij}\bigsp \hbox{ (for $s=0,\ldots ,m-1$)  }\ess .
\eqno 1.38
$$
Since $T_{\mu/i,j+1}=T_{\mu/ij}/q $, substituting 1.37 and 1.38 in 1.36 a) gives
$$
C_{\mu/i,j+1}  \ses {1\over M}{T_{\mu/ij}\over q\nabla}\ess 
\biggl(\ssp \prod_{s=1}^{m}\ssp \Bigl(1-\nabla {u_s^{ij}\over T_{\mu/ij}}\ssp \Bigr)
 \biggr)\ssp \Bigl(  1-{q\ssp u_0^{ij}\nabla\over T_{\mu/ij} } \Bigr)\ssp \phi_{S}^{(m)}\ssp .
\eqno 1.39
$$
Setting again $z_0^{ij}=u_0^{ij}\nabla/T_{\mu/ij}$ and $z_m^{ij}=u_m^{ij}\nabla/T_{\mu/ij}$,
we see that 1.39 may be written as
$$
C_{\mu/i,j+1}  \ses CF\ess {T_{\mu/ij}\over q}\ssp \bigl(1-q\ssp z_0^{ij}\bigr)\bigl(1- z_m^{ij}\ssp \bigr) 
\ssp \phi_{S}^{(m)}\ssp ,
$$
which is in perfect agreement with 1.32.
\sas

Similarly,  using 1.37 and the identities
$$
u_s^{i+1,j+1}= {1\over qt}\ssp u_{s+1}^{ij}\ess\ess\ess  \hbox{(for $s=0,\ldots ,m-2$)}
\ess\scs \ess\ess\ess u_{m-1}^{i+1,j+1}={1\over q}\ssp u_m^{ij}\ess\ess\ess
 \hbox{  and   } \ess\ess T_{\mu/i+1,j+1}=T_{\mu/ij}/qt
$$
we may write 1.36 b) as
$$
C_{\mu/i+1,j+1}  \ses {1\over M}{T_{\mu/ij}\over qt\nabla}\ess 
\biggl(\ssp \prod_{s=1}^{m-1}\ssp \Bigl(1-\nabla {u_s^{ij}\over T_{\mu/ij}}\ssp \Bigr)
 \biggr)
\ssp \Bigl(  1-{t\ssp u_m^{ij}\nabla\over T_{\mu/ij} } \Bigr)
\ssp \Bigl(  1-{q\ssp u_0^{ij}\nabla\over T_{\mu/ij} } \Bigr)
\ssp \phi_{S}^{(m)}\ssp ,
$$
which is easily seen to be again in perfect agreement with 1.32.
The case we have just considered should be sufficient evidence that 
we have an underlying mechanism here that forces the same final answer to
come out in all the possible cases, completing our proof of Theorem I.1.
\sas

An immediate consequence of I.50 is the following remarkable fact
\sas

\heading{\bol Theorem 1.1}

{\ita Under the SF-heuristic and the $n!$ conjecture, all the modules $\BA_{ij}^x$, $\BA_{ij}^y$, for $\mu\part n+1$, are 
bigraded versions of the left regular representation of $S_n$.
}

\heading{\bol Proof}

Note that formula I.50 may also be written in the form
$$
A_{ij}^x/q^a\ses A_{ij}^y/t^l\ses \sum_{k=1}^{m}\ssp 
{(-\nabla)^{m-k}\phi_S^{(m)}\over  T_{\mu/ij}^{m-k}  }  
\ssp e_{m-k}\bigl[\ssp u_1^{ij}+u_2^{ij}+\cdots +u_{m-1}^{ij}\ssp \bigr] \ess ;
$$
thus, using I.35 we get
$$
A_{ij}^x/q^a\ses A_{ij}^y/t^l\ses \sum_{k=1}^{m}\ssp {\phi_S^{(k)}\over  T_{\mu/ij}^{m-k} }  
\ssp e_{m-k}\bigl[\ssp u_1^{ij}+u_2^{ij}+\cdots +u_{m-1}^{ij}\ssp \bigr] \ess .
\eqno 1.40
$$
To compare this formula with I.41 we should set there, for every $\aaa_i\in S$, 
$$
T_{\aaa_i}= T_{\mu}/x_i
$$ 
and obtain
$$
\TH_{\aaa_i}\ses \sum_{k=1}^{m}\ssp {\phi_S^{(k)}\over  T_{\mu }^{m-k} }  
\ssp e_{m-k}\bigl[\ssp x_1 +x_2 +\cdots +x_{m}\sms x_i\ssp \bigr] \ess .
\eqno 1.41
$$
Since $\aaa_i\RA  \mu$ gives that $\aaa_i\part n$, we have the expansion
$$
\TH_{\aaa_i}\ses \sum_{\la\part n}\ssp S_\la(x) \ssp \TK_{\la,\aaa_i}(q,t)
$$
which, together with the Macdonald result $\TK_{\la,\aaa_i}(1,1)=f_\la$, yields that
$$ 
\TH_{\aaa_i}(x;1,1)\ses \sum_{\la\part n}\ssp  S_\la(x) \ssp f_\la \ses  h_1^n(x)\ess .
$$
Thus, placing $q=t=1$ in 1.40 and 1.41 we obtain that
$$
A_{ij}^x(x;1,1)\ses A_{ij}^y(x;1,1)\ses
\sum_{k=1}^{m}\ssp {\phi_S^{(k)}(x;1,1)}  
\ssp { m-1 \choose k-1} \ses \TH_{\aaa_i}(x;1,1)\ses h_1^n(x)\ess ,
$$
which proves that both $A_{ij}^x$ and $ A_{ij}^y$ are Frobenius characteristics of
bigraded left regular representations. 
\sas

Another interesting identity relating the characteristics $A_{ij}^x$ and $ A_{ij}^y$ may be derived
from the theory of Macdonald polynomials as well as the present heuristics:
\sas

\heading{\bol Theorem 1.2}
$$
T_{\mu/ij} \DA   A_{ij}^x\ses A_{ij}^y\ess .
\eqno 1.42
$$

\heading{\bol Proof}

Due to the fact that $C_{\mu/ij}(x;q,t)$ is the bivariate Frobenius characteristic 
of $\BM_{\mu/ij}=\CL_\del[\DD_{\mu/ij}]$,
from I.25 with $\DD=\DD_{\mu/ij}$ we get that
$$
T_{\mu/ij} \DA C_{\mu/ij}\ses C_{\mu/ij}\ess .
\eqno 1.43
$$
Now  I.19 a) and b) give
$$
\eqalign{
a)\ess\ess A_{ij}^x &\ses C_{\mu/ij}\sms t\ssp C_{\mu/i+1,j}\sms  \ssp C_{\mu/i,j+1}\sps  t\ssp C_{\mu/i+1,j+1}\ess ,\cr
b)\ess\ess  A_{ij}^y &\ses C_{\mu/ij}\sms  q  \ssp C_{\mu/i,j+1}\sms  \ssp C_{\mu/i+1,j}\sps  q\ssp C_{\mu/i+1,j+1}\ess .\cr
}
\eqno 1.44
$$
Using 1.43 on each term of the expansion in 1.44 a) yields
$$
 \DA   A_{ij}^x\ses 
{1 \over T_{\mu/ij}}\ssp 
C_{\mu/ij}-  
{1 \over t\ssp  T_{\mu/i+1,j}}\ssp 
C_{\mu/i+1,j}-  
{1 \over  \ssp  T_{\mu/i ,j+1}}
\ssp C_{\mu/i,j+1}+  
{1 \over t\ssp  T_{\mu/i+1,j+1}}
\ssp C_{\mu/i+1,j+1}\ess .
$$
Multiplying both sides by $T_{\mu/ij}$ and using the identities 
$$
T_{\mu/ij}\ses t\ssp T_{\mu/i+1,j}\ses q\ssp T_{\mu/i,j+1}\ses tq\ssp T_{\mu/i+1,j+1}
$$
we finally obtain
$$
T_{\mu/ij} \DA   A_{ij}^x\ses 
C_{\mu/ij}\sms C_{\mu/i,j+1}\sms q\ssp C_{\mu/i,j+1} \sps q \ssp C_{\mu/i+1,j+1}
$$
whose right-hand side is seen to be a rearrangement of the right-hand side of 1.44 b).
This proves 1.42.
\sas

\heading{\bol Remark 1.2}

We should emphasize at this point that each symmetric function identity  
we write down here may be studied 
from two different viewpoints. On one hand it can be viewed as an identity 
involving Macdonald polynomials and may be verified using purely symmetric function manipulations.
On the other hand, if we view it as an identity relating two bigraded Frobenius characteristics,
we may try to give it a representation theoretical proof. It develops that 1.42, 
which here and after we shall refer to as the ``{\ita flip identity,}'' may also be shown
in this manner. It is significant that the ``crucial identity,'' which on the surface
appears quite similar, nevertheless turns out so much more difficult to prove.
\sas

Our point of departure is the introduction of a bilinear form in each  of the spaces
$\BM_{\mu/ij}\ssp ,$ which is defined by setting
$$
\LLL  P\scs Q \RRR \ses \LL \ssp  \flip_{ij}^{-1} P\scs Q\ssp \RR 
	\ess ,
\eqno 1.45
$$
where for convenience,  we set  
$$
\flip_{ij}^{-1} \ses \flip_{\DD_{\mu/ij}}^{-1} \ess .
$$
In other words, for any polynomial $(x;y),\ess $  $\flip_{ij}^{-1} P$ denotes the unique polynomial
$P_1\in \BM_{\mu/ij}$ such that $P=P_1(\del)\DD_{\mu/ij}$. In particular, we see that
if $P_1=\flip_{ij}^{-1} P$ and $Q_1=\flip_{ij}^{-1} Q$ then  1.45 may also be rewritten as
$$
\LLL  P\scs Q \RRR \ses P_1(\del)Q_1(\del)\ssp \DD_{\mu/ij}\ess \big |_{x=y=0}
	\ess ,
$$
yielding that $\LLL\scs \RRR$ is a symmetric bilinear form.
Now it develops that this form may be used to construct a nondegenerate pairing of $\BA_{ij}^x$ 
with $\BA_{ij}^y$ that forces the identity in 1.42.
More precisely, we have the following remarkable result. 
\sas

\heading{\bol Proposition 1.2}

{\ita The two spaces $\BA_{ij}^x$ and $\BA_{ij}^y$ have the same dimension and we
 can construct two bihomogeneous bases $\CB_{ij}^x=\{f_1^x,f_2^x,\ldots ,f_N^x\}$ and
$\CB_{ij}^y=\{f_1^y,f_2^y,\ldots ,f_N^y\}$ for $\BA_{ij}^x$ and $\BA_{ij}^y$ respectively such that
$$
\LLL f_r^x\scs f_s^y \RRR\ses 
\cases {
1 & if $r=s$\ssp ,\cr\cr
0 & if $r\neq s$\ssp .\cr
}
\eqno 1.46
$$
In particular, we must have
$$
\weight(f_r^x)\ess \times \ssp  \weight(f_r^y)\ses T_{\mu/ij}
	\ess ,
\eqno 1.47
$$
where for convenience for a bihomogeneous polynomial $P$ of bidegree $(h,k)$ we set
$$
\weight(P)=t^hq^k\ess .
$$ 
Moreover, if for all $\sig\in S_n$ we have
$$
\sig f_s^x\ses \sum_{r=1}^N\ssp f_r^x\ess a_{r,s}(\sig)
\eqno 1.48
$$
then 
$$
\sig f_s^y\ses \sign(\sig)\ess \sum_{r=1}^N\ssp f_r^y\ess a_{s,r}(\sig^{-1})\ess .
\eqno 1.49
$$
}

\heading{\bol Proof}

Note that I.29 gives the orthogonal decompositions 
$$
\BM_{\mu/ij}= \BM_{\mu/i,j+1}
\ssp \oplus_\perp  \ess \flip_{ij}^{-1}\ssp \BK_{ij}^y
\ess \scs \ess\ess\ess
\BM_{\mu/ij}= \BM_{\mu/i+1,j}
\ssp \oplus_\perp  \ess \flip_{ij}^{-1}\ssp \BK_{ij}^x \ess .
\eqno 1.50
$$
In particular this means that, $\ssp$ for $P\in \BM_{\mu/ij}\ssp ,$  we have
$$
\LL\ssp P\scs Q \ssp \RR \ses 0 
$$
for all $Q\in \flip_{ij}^{-1}\ssp \BK_{ij}^y$ if and only if 
$$
P\in \BM_{\mu/i,j+1}\ess .
$$
We thus deduce the equivalences
$$
\ess P\in \BK_{ij}^x \ess
	  \hbox{~~and~~}
\ess P\in(\flip_{ij}^{-1})^\perp\ssp \BK_{ij}^y\ess
\ess\ess \Longleftrightarrow \ess\ess \ess  P\in \BK_{i,j+1}^x\ess .
\eqno 1.51
$$
Similarly we derive that
$$
\ess P\in \BK_{ij}^y \ess
	\hbox{~~and~~}
\ess P\in(\flip_{ij}^{-1})^\perp \ssp \BK_{ij}^x\ess
\ess\ess \Longleftrightarrow \ess\ess \ess  P\in \BK_{i+1,j}^y\ess .
\eqno 1.52
$$
In view of our definition 1.45 of the form $\LLL\scs \RRR$ we deduce from 1.51 and 1.52
that if $P_1,P_2\in \BK_{ij}^x$ and $Q_1,Q_2\in \BK_{ij}^y$, with 
$P_1-P_2\in \BK_{i,j+1}^x$ and $Q_1-Q_2\in \BK_{i+1,j}^y\ssp$,
then 
$$
 \LLL P_1-P_2\scs Q_1 \RRR =0\ess\ess  {\rm and}\ess\ess   \LLL P_2\scs Q_1-Q_2 \RRR= 0\ess .
$$
Thus
$$
\LLL P_1 \scs Q_1 \RRR= \LLL P_1-P_2 \scs Q_1 \RRR + \LLL P_2 \scs Q_1 \RRR 
\ses \LLL P_2 \scs Q_1 \RRR 
$$
and similarly
$$
\LLL P_2 \scs Q_1 \RRR= \LLL P_2 \scs Q_1-Q_2 \RRR + \LLL P_2 \scs Q_2 \RRR 
\ses \LLL P_2 \scs Q_2 \RRR \ess ,
$$
yielding
$$
\LLL P_1 \scs Q_1 \RRR\ses \LLL P_2 \scs Q_2 \RRR \ess .
\eqno 1.53
$$
This shows that $\LLL \scs \RRR$ is a well-defined pairing of $\BA_{ij}^x=\BK_{ij}^x/\BK_{i,j+1}^x$
with $\BA_{ij}^y=\BK_{ij}^y/\BK_{i+1,j}^x\ess $. We are left to show that it is nondegenerate.
To this end suppose that for some representative element $P\in \BK_{ij}^x$ of $\BK_{ij}^x/\BK_{i,j+1}^x$ we have
$$
\LLL P\scs Q\RRR =\LL\ssp  \flip^{-1} _{ij}Q\scs P\ssp \RR\ses 0
\eqno 1.54
$$ 
for all representatives $Q\in \BK_{ij}^y$ of $\BK_{ij}^y/\BK_{i+1,j}^x$. In view of
1.53, the relation in 1.54 must hold true for all $Q\in \BK_{ij}^y$, but then 
the first equation in 1.50 yields that $P\in \BM_{\mu/i,j+1}$ and this, together with $P\in \BK_{ij}^x\,$,
forces $P\in \BK_{i,j+1}^x$. In other words, $P$ is equal to zero in the quotient  $\BK_{ij}^x/\BK_{i,j+1}^x$.
Similarly we show that 1.54  for all $P\in \BA_{ij}^x$ can only hold for $Q=0$ in $\BA_{ij}^y$.
Thus $\LLL \scs \RRR$ is nondegenerate. 

Now let $\{f_1,f_2,\ldots, f_N\}$ and $\{g_1,g_2,\ldots, g_M\}$
be  bihomogeneous bases for $\BA_{ij}^x$ and $\BA_{ij}^y$, and set
$$
R_{i,j}\ses \LLL f_i ,g_j\RRR\ess .
$$ 
Note that we cannot have $N<M$ for otherwise we would be able to construct a nontrivial solution 
$c_1,c_2,\ldots ,c_M$ of the homogeneous system of equations 
$$
c_1 R_{i,1}+c_2 R_{i,2}+\cdots +c_M R_{i,M}\ses 0 \bigsp (\ess {\rm for}\ess\ess i=1,2,\ldots ,N\ess ) 
$$
and this would contradict the nondegeneracy of $\ess \LLL\scs \RRR$. For the same reason
we can't have $N>M$ nor $M=N$ with $R=  \| R_{i,j} \|_{i=1}^ N$ a singular matrix. 
Thus $\BA_{ij}^x$ and $\BA_{ij}^y$ have the same dimension  
and $R$ must be invertible. This given, the two bases $\{f_1^x,f_2^x,\ldots ,f_N^x\}$ 
with the asserted properties may be obtained by setting
$$
\{f_1^x,f_2^x,\ldots ,f_N^x\}\ses \{f_1 ,f_2 ,\ldots ,f_N \} 
$$
and
$$
\{f_1^y,f_2^y,\ldots ,f_N^y\}\ses \{f_1^x,f_2^x,\ldots ,f_N^x\}\times R^{-1}\ess .
$$
With this choice, 1.46 is immediate and then 1.47 follows from the fact that if for two bihomogeneous 
polynomials $P ,Q$ we have 
$$
P(\del )Q(\del)\DD_{\mu/ij}=1
$$
then necessarily their bidegrees must add up to the bidegree of $\DD_{\mu/ij}$.
Finally, to show that 1.49 follows from 1.48 note first that from 1.46 we derive
that the expansion of any element $Q\in \BA_{ij}^y$ in terms
of the basis $\{f_1^y,f_2^y,\ldots ,f_N^y\}$ may be written in the form
$$
Q\ses \sum_{r=1}^N\ssp f_r^y \ess  \LLL Q,f_r^x\RRR \ess .
$$
Thus we may write
$$
\sig   f_s^y \ses \sum_{r=1}^N\ess f_r^y\ess \LLL  \sig   f_s^y\scs f_r^x\RRR \ess .
\eqno 1.55
$$
However, we see that
$$
\eqalign{
\LLL  \sig   f_s^y\scs f_r^x\RRR &= \LL\flip_{ij}^{-1}\sig f_s^y \scs f_r^x\RR
= \sign(\sig)\ssp \LL\sig \, \flip_{ij}^{-1}  f_s^y \scs f_r^x\RR \cr 
&= \sign(\sig)\ssp \LL \flip_{ij}^{-1} f_s^y \scs \sig ^{-1} f_r^x\RR= 
\sign(\sig)\ssp \LLL f_s^y\scs \sig ^{-1} f_r^x\RRR\ess .\cr
}
$$
Substituting this in 1.55 gives
$$
\sig   f_s^y \ses \sign(\sig) \ess \sum_{r=1}^N\ess f_r^y\ess \LLL  f_s^y\scs \sig ^{-1}f_r^x\RRR \ess .
\eqno 1.56
$$
Now, from 1.48 for $\sig^{-1}$ and $r,s$ interchanged we derive that
$$
\sig^{-1}   f_r^x \ses  \ess \sum_{s=1}^N\ess f_s^x\ess a_{s,r}(\sig ^{-1})  
$$
and 1.46 then gives that
$$ 
\LLL  f_s^y\scs \sig ^{-1}f_r^x\RRR\ses  a_{s,r}(\sig ^{-1}) \ess .
$$
Substituting in 1.56 yields 1.49 as desired, completing our proof.
\sa

\heading{\bol Remark 1.3}

We should note that the fact that $\BA_{ij}^x$ and $\BA_{ij}^y$ have the same dimension is 
also an immediate consequence of 1.44 a) and b). In fact, setting $q=t$ there yields the stronger result
that these two modules (graded by total degree) have identical Frobenius characteristics.
We should also emphasize that this argument as well as the proof of Proposition 1.2 
makes no use of any of our conjectures nor any identification of the polynomials 
$C_{\mu/ij}$ with expressions (such as in 1.20) involving the Macdonald polynomials $\TH_\mu$. 
\sa

\heading{\bol Remark 1.4}

Proposition 1.2 leads to an alternate proof of Theorem 1.2 and 
a direct representation theoretical interpretation of the  identity in 1.42.
To see this note that 1.48 yields that the bigraded  characters  of $\BA_{ij}^x$
and $\BA_{ij}^x$ ,
are respectively given by the expressions
$$
\bigl(\ch \BA_{ij}^x\bigr)(\sig;q,t) = \sum_{s=1}^N\ssp \weight(f_r^x)\ssp a_{r,r}(\sig)
\ssp   ,
\ess \ess 
\bigl( \ch \BA_{ij}^y\bigr)(\sig;q,t) = \sign(\sig)\ssp 
\sum_{s=1}^N\ssp \weight(f_r^y)\ssp a_{r,r}(\sig^{-1})\ess .
\eqno 1.57
$$ 
Now,  1.46 gives 
$$
\weight(f_r^y)\ses T_{\mu/ij} \ssp /\weight(f_r^y)\ess ,
$$
and from 1.57 we derive that 
$$
\bigl( \ch \BA_{ij}^y\bigr)(\sig;q,t)= T_{\mu/ij}\ssp \sign(\sig)\ssp 
\sum_{s=1}^N  a_{r,r}(\sig^{-1})/ \weight(f_r^x)
=
T_{\mu/ij}\ssp \sign(\sig)\ssp \bigl( \ch \BA_{ij}^x\bigr)(\sig;1/q,1/t)
\ess .
$$
Equating the Frobenius images of both sides yields 1.42.
\sa

We shall terminate this section by showing that 
a proof of the ``crucial identity'' 
would in one stroke establish Conjecture I.2 as well as all the conjectured 
expansions in I.16. This implication is based on a result of M.\ Haiman in [14]
which asserts that a proof of the $n!$ conjecture for a given $\mu$
yields that the bigraded Frobenius characteristic of $\BM_\mu$ for that
same $\mu$ must be given by the polynomial $\TH_\mu(x;q,t)$. Since the 
$n!$ conjecture has been verified by computer for all $|\mu|\leq 8$, the
argument can proceed by induction on $|\mu|$. So let us assume that
for a given $\mu\part n$ we have $\dim \BM_\nu=(n-1)!$ for 
all $\nu\RA \mu$. The Haiman result then yields that for all $\nu\RA \mu$
the bigraded Frobenius characteristic of  $\BM_\nu$ is $\TH_\nu$. 
Since I.20 is just another way of writing the four term recursion in I.17,
its validity implies (by Theorem I.1) that I.16 must hold true as well. 
Now, as we have seen, I.16, for $(i,j)=(0,0)$ states (via the Macdonald identity
in I.13) that 
$$
C_{\mu/00}(x;q,t)\ses \del_{p_1} \TH_\mu(x;q,t)\ess .
$$ 
Combining this with Proposition I.5 gives
$$
\del_{p_1} C_{\mu}(x;q,t)\ses \del_{p_1} \TH_\mu(x;q,t)\ess .
$$ 
In particular, applying $\del_{p_1}^{n-1}$ to both sides we get that the bigraded Hilbert series of
the module $\BM_\mu$ is given by the polynomial
$$
F_\mu(q,t)\ses \del_{p_1}^{n}\TH_\mu(x;q,t)\ses \sum_{\la\part n}\ssp f_\la \ssp \TK_{\la\mu}(q,t) \ess .
$$
Here, the last equality follows from I.8. But now, the Macdonald result
that $\TK_{\la\mu}(1,1)=f_\la$ yields that
$$
\dim \BM_\mu\ses F_\mu(1,1)\ses \sum_{\la\part n}\ssp f_\la^2\ses n!\ess ,
$$
completing the induction. Then of course we can combine this with Haiman's result
and obtain that the $\TK_{\la\mu}(q,t)$ are polynomials with positive integer coefficients. 
To show that I.20 implies Conjecture I.2, we use I.16 with the $c_{\tau\rho}(q,t)$
given by I.14 and obtain
$$
C_{\mu/ij}(x;q,t)\ses \sum_{\rho\RA \tau}
\prod_{s\in \CR_{\tau/\rho}}\hskip -.07in
{ t^{l_\tau (s)}-q^{a_\tau (s)+1}\over { t^{l_\tau (s)}-q^{a_\tau (s)}} }  
\hskip -.05in
\prod_{s\in \CC_{\tau/\rho}} 
\hskip -.07in
{
q^{a_\tau (s)}- t^{l_\tau (s)+1}
\over 
{q^{a_\tau (s)}-t^{l_\tau (s)}} 
}
\ess
\TH_{\mu-\tau+\rho}(x;q,t)\ess .  
$$ 
Now this identity, for $t=1/q$,  may be rewritten as
$$
\eqalign{
C_{\mu/ij}(x;q,1/q) &\ses \sum_{\rho\RA \tau}
\Bigl(\prod_{s\in R_{\tau/\rho}\cup \ssp C_{\tau/\rho} }\hskip -.1in{1-q^{h_\tau(s)}\over 1-q^{h_\rho(s)}}
\ess \Bigr) {1\over q^{|C_{\tau/\rho}|} }\ess  
\TH_{\mu-\tau+\rho}(x;q,1/q) \cr
&\ses \sum_{\rho\RA \tau}
{1\over 1-q}{  \prod_{s\in \tau} (1-q^{h_\tau(s) }) \over \prod_{s\in \rho} (1-q^{h_\rho(s)}) }
\ess  {1\over q^{|C_{\tau/\rho}|} }\ess  
\TH_{\mu-\tau+\rho}(x;q,1/q) \ess . \cr
}  
\eqno 1.58
$$ 
Here, the symbols $h_\tau(s),h_\rho(s)$ denote the hook lengths of the cell $s$
with respect to the two partitions $\tau$ and $\rho$.
Using the fact that $\TH_{\mu-\tau+\rho}(x;1,1)=h_1^{n-1}$, we see that
letting $q\RA 1$ reduces 1.58 to
$$
C_{\mu/ij}(x;1,1) \ses \Bigl(\sum_{\rho\RA \tau}
{  \prod_{s\in \tau}  h_\tau(s)   \over \prod_{s\in \rho}  h_\rho(s)  }
\ess
\Bigr)  h_1^{n-1}\ess .
\eqno 1.59
$$
Now the classical recursion for the number of standard tableaux gives
$$
\sum_{\rho\RA \tau}
{  \prod_{s\in \tau}  h_\tau(s)   \over \prod_{s\in \rho}  h_\rho(s)  }
\ses {|\tau|\over f_\tau} \sum_{\rho\RA \tau}\ssp  f_\rho \ses |\tau|\ess .
$$
Thus, 1.59 may be rewritten as
$$
C_{\mu/ij}(x;1,1) \ses |\tau|\ssp  h_1^{n-1}\ssp , 
$$
which establishes that $\BM_{\mu/ij}$ consists of $|\tau|$
occurrences of the left regular representation of $S_{n-1}\ssp $, 
precisely as asserted  by Conjecture I.2.  

\vfill\supereject

\def\st{\,:\,}
\def\set#1{\left\{#1\right\}}
\def\setof#1#2{\set{\,{#1}\st{#2}\,}}

\def\Mepsd{\BM^{\epsilon_1\cdots\epsilon_d}}


\heading{\bol 2. Conjectural Bases and the ``crucial identity''. }
\sa
As we have seen in the previous section, the proof of Conjecture I.3
is reduced to establishing the ``crucial  identity''
$$
t^l\ssp A_{ij}^x\ses q^a \ssp A_{ij}^y\ess .
\eqno 2.1
$$
Although at the moment we are unable to prove this identity except in some
special cases (see the next section),
we can nevertheless search for
the underlying representation theoretical mechanism that produces it.  
Our main goal in this section is to provide such a mechanism. This will be obtained
by means of an algorithm for constructing bases for all of our spaces $\BM_{\mu/ij}$
which is an extension of an algorithm given in [1]. All our constructs here,
as in [1],
are heavily dependent on the SF-heuristic, and as such they are still conjectural.
Nevertheless, the validity of our arguments is strongly supported by amply
verifiable theoretical and numerical consequences.   
Remarkably, these heuristics not only yield 2.1 but reveal that both 2.1 and  
the ``flip'' identity 
$$
A_{ij}^x\ses T_{\mu/ij}\ssp \da A_{ij}^y 
\eqno 2.2
$$
are consequences of considerably more refined versions. Before we can state these results 
we need to introduce some notation. Given that 
$$
\Pred (\mu)\ses \bigl\{\ssp \nu^{(1)},\nu^{(2)},\ldots ,\nu^{(d)}\ssp \bigr\}\ess ,
$$
It will be convenient here to use the shorter symbol $S_{ij}$ to represent the subset
$\Pred_{ij}(\mu)$ define in I.31. That is, we are setting
$$
S_{ij}\ses \bigl\{\ssp \nu^{(i)}\ssp :\ssp \mu/\nu^{(i)}\ess\hbox{is in the shadow of $(i,j)$}\ess \bigr\} 
\eqno 2.3
$$
Given that $S_{ij}=\{\nu^{(i_1)},\nu^{(i_2)},\ldots ,\nu^{(i_m)}\}$,    with   $i_1<i_2<\cdots <i_m$,
here and after we shall identify a subset $T$ of $S_{ij}$ with its corresponding 
$0,1$-word $\eee(T)=\eee_1 \cdots  \eee_m$
defined by setting $\eee_s=1\ssp{\rm or} \ssp 0\ess $ according as $\nu^{(i_s)}\in T$ or $\nu^{(i_s)}\not\in T$. 
Conversely, given such a word $\eee$, we shall set
$$
T(\eee)\ses \{\ssp \nu^{(i_s)}\ssp :\ssp \eee_s=1\ssp \}\ess .
$$
This given, recalling the definition in I.33, we shall also set (when $|S_{ij}|=m$) 
$$
\BM_{ij}^{\eee_1\cdots\eee_m}\ses \BM_{S_{ij}}^{T(\eee)}\ess .
\eqno 2.4
$$ 
Assuming that the corners of $\mu$ in the shadow of $(i,j)$ have weights
$$
x_r^{ij}\ses  t^{l_r'}q^{a_r'}\bigsp (\ssp {\rm for}\ess\ess r=1\ldots m\ssp ) 
$$
we shall set
$$
w_r^{ij}=  a_r'- a_{r-1}'\ess\ess\ess {\rm and }\ess\ess\ess v_r^{ij}= l_r'- l_{r+1}'\ess .
$$
Of course when dealing with a fixed pair $(i,j)$ we shall drop the superscripts $ij$ and simply write
$x_r, w_r, v_r$. In the figure below we have illustrated the geometric meaning of the parameters $w_r$
and $v_r$ as representing the exposed ``{\ita width}'' of corner $r$ and the vertical ``{\ita drop}''
immediately after it.
$$
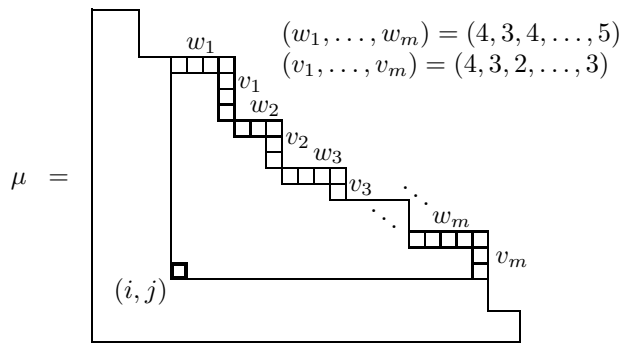
\hbox{%
	\myput(0,200){\llap{$\mu\;\;=\;\;\;$}}
	\hseg(0,0;540)					
	\vseg(540,0;40)
	\hseg(500,40;40)
	\myput(510,100){$v_m$}
	\myput(430,150){$w_m$}
	\vseg(500,40;80)
	\multiputline(500,140)(0,-20){3}{\corn(0,0)}
	\multiputline(500,140)(-20,0){5}{\corn(0,0)}
	\vseg(400,140;40)\hseg(320,180;80)
		\myput(390,170){$\ddots$}
		\myput(350,140){$\ddots$}
	\myput(325,190){$v_3$}
	\myput(280,230){$w_3$}
	\multiputline(320,220)(0,-20){2}{\corn(0,0)}
	\multiputline(320,220)(-20,0){4}{\corn(0,0)}
	\myput(245,250){$v_2$}
	\myput(200,290){$w_2$}
	\multiputline(240,280)(0,-20){3}{\corn(0,0)}
	\multiputline(240,280)(-20,0){3}{\corn(0,0)}
	\myput(185,320){$v_1$}
	\myput(120,370){$w_1$}
	\multiputline(180,360)(0,-20){4}{\corn(0,0)}
	\multiputline(180,360)(-20,0){4}{\corn(0,0)}
	\hseg(60,360;40)
	\vseg(60,360;60)
	\hseg(0,420;60)
	\vseg(0,0;420)
%
	\hseg(100,80;400)
	\vseg(100,80;280)
	%
	\bcorn(120,100)
	\labLL(120,100;{(i,j)})
	\myput(240,380){$(w_1,\ldots,w_m)=(4,3,4,\ldots,5)$}
	\myput(240,340){$(v_1,\ldots,v_m)=(4,3,2,\ldots,3)$}
}
\hskip 540\mylength
$$
Given a subset $T=T(\eee)\con S_{ij}$ we  shall let $D_{ij}(T)$ denote the 
subdiagram of $\mu$ obtained by the following construction:
 \item{}{\ita Divide the shadow of $(i,j)$ in $\mu$
into $m$ rectangles, of widths $\ssp w_1,\ldots,w_m\ssp $
from left to right, by dropping vertical lines from each
of its corners.  Then delete the $r^{th}$ rectangle if $\epsilon_r=0$, and
slide the remaining rectangles horizontally left to fill the
gaps, setting the southwest corner of the leftmost rectangle at $(i,j)$.  This done, 
the  cells covered by the  resulting rectangles form $D_{ij}(T)$.}
\sas

\noindent
In the figure below we illustrate this construction when $\ess \mu=(15,15,11,11,6,6,6,6,3,3,2,2)\ess $
$m=5$, $\ess i=j=0\ess $ and $T=\{2,3,5\}$ or $\eee=01101\,$.
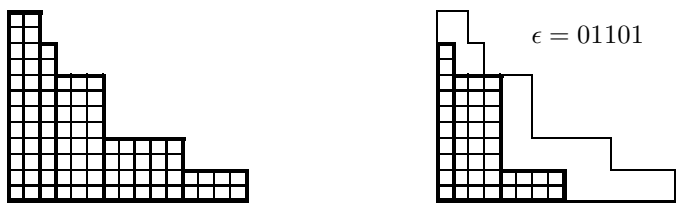
 $$
\hbox{%
	\multiputgrid(40,240)(-20,0){2}(0,-20){12}{\corn(0,0)}
	\multiputgrid(60,200)(-20,0){1}(0,-20){10}{\corn(0,0)}
	\multiputgrid(120,160)(-20,0){3}(0,-20){8}{\corn(0,0)}
	\multiputgrid(220,80)(-20,0){5}(0,-20){4}{\corn(0,0)}
	\multiputgrid(300,40)(-20,0){4}(0,-20){2}{\corn(0,0)}
%
	\brectangle(40,240)(40,240)
	\brectangle(60,200)(20,200)
	\brectangle(120,160)(60,160)
	\brectangle(220,80)(100,80)
	\brectangle(300,40)(80,40)
}\hskip 300\mylength
\hskip 1in
%
\hbox{%
	\hseg(0,0;300)
	\vseg(300,0;40)
	\hseg(220,40;80)
	\vseg(220,40;40)
	\hseg(120,80;100)
	\vseg(120,80;80)
	\hseg(60,160;60)
	\vseg(60,160;40)
	\hseg(40,200;20)
	\vseg(40,200;40)
	\hseg(0,240;40)
	\vseg(0,0;240)
	\multiputgrid(20,200)(-20,0){1}(0,-20){10}{\corn(0,0)}
	\multiputgrid(80,160)(-20,0){3}(0,-20){8}{\corn(0,0)}
	\multiputgrid(160,40)(-20,0){4}(0,-20){2}{\corn(0,0)}
%
	\brectangle(20,200)(20,200)
	\brectangle(80,160)(60,160)
	\brectangle(160,40)(80,40)
	\myput(120,200){$\epsilon=01101$}
}\hskip 300\mylength
$$

We need one further convention before we can present our basic construction.
In some of the formulas that follow it will be more illuminating to use the symbol
``$\BM_1(\del)\DD$'' rather than ``$\flip_\DD\ssp \BM_1$'' to denote the image
of $\BM_1$ by $\flip_\DD$. In other words, we are setting
$$
\BM_1(\del)\DD\ses \bigl\{\ssp P(\del)\DD\ssp :\ssp P\in \BM_1\ssp \bigr\}\ess .
\eqno 2.5
$$
This given, extensive numerical and theoretical evidence strongly suggests that 
\sa

\heading{\bol Conjecture 2.1}

{\ita The space $\BM_{\mu/ij}$ affords the following direct sum decomposition:}
$$
\BM_{\mu/ij}\ses \bigoplus_{\emptyset\neq T\con S_{ij}}\ess 
\bigoplus_{(i',j')\in D_{ij}(T)}\ess \BM_{S_{ij}}^T(\del)\DD_{\mu/i'j'}\ess .
\eqno 2.6
$$
The constructions underlying this identity are of course heavily dependent
on the ``Science Fiction  Conjecture'' (see [1]) which states that the modules 
$\BM_{\nu_1},\BM_{\nu_2},\dots  ,\BM_{\nu_d} $ generate a distributive lattice under 
span and intersection. Under this assumption, 2.6 yields an algorithm for constructing
bihomogeneous bases for the modules $\BM_{\mu/ij}$. To be explicit, this algorithm
consists in preconstructing bihomogeneous bases $\CB_S^T$ for all the subspaces $\BM_S^T$ 
given in I.33 and for all pairs 
$$
\bigl\{\ssp (T,S)\ssp :\ssp \emptyset\neq T\con S\con \{1,2,\ldots ,d\}\ess \bigr\}\ess .
$$ 
This done, for all $(i,j)\in \mu$ a basis for $\BM_{\mu/ij}$ should be given by the collection
$$
\CB_{\mu/ij}\ses \sum_{\emptyset\neq T\con S_{ij}}\ess 
\sum_{(i',j')\in D_{ij}(T)}\ess \flip_{i'j'}\ssp \CB_{S_{ij}}^T 
\eqno 2.7
$$
where for convenience we have set $\ess\flip_{i'j'}= \ssp \flip_{\DD_{\mu/i'j'}}\ssp $.
\sa

Before we proceed any further it will be good to see what 2.6 yields in at least one
concrete example. We shall illustrate it in the case $\mu=(3,2,1)$ and $(i,j)=(0,0)$.
To this end, we begin by noting that under the SF hypotheses, in any three-corner case, the module
$\ess \BM_{\nu_1}+\BM_{\nu_2}+\BM_{\nu_3}\ess  $ decomposes into the direct sum
of the submodules $\BM^{\eee_1\eee_2\eee_3}$.
For $\mu=(3,2,1)$ we have $\nu_1=(3,2)\scs \nu_2=(3,1,1)\scs \nu_3=(2,2,1)$. Accordingly we have the 
direct sum decompositions 
$$
\eqalign{
\BM_{32}&\ses \BM^{100}\oplus\BM^{110}\oplus\BM^{101}\oplus\BM^{111} \cr \cr 
\BM_{311}&\ses \BM^{010}\oplus\BM^{110}\oplus\BM^{011}\oplus\BM^{111} \cr \cr 
\BM_{221}&\ses \BM^{001}\oplus\BM^{101}\oplus\BM^{011}\oplus\BM^{111} \ess .\cr 
}
\eqno 2.9
$$
After constructing bases $\CB^{\eee_1\eee_2\eee_3}$ for each of the submodules $\BM^{\eee_1\eee_2\eee_3}$
appearing above, the result of applying the recipe in 2.7 with $(i,j)=(0,0)$ and $S_{00}=\{1,2,3\}$
may be described by the following diagrams:
$$
\eee(T)={100} \ssp\longrightarrow   \tableau{
\CB^{100}\cr
\CB^{100}& 0\cr
\CB^{100}&0&0 \cr}
\bigsp
\eee(T)={101} \ssp\longrightarrow  \tableau{
\CB^{101}\cr
\CB^{101}& 0\cr
\CB^{101}&\CB^{101}&0 \cr}
 $$ $$
\eee(T)={111} \ssp\longrightarrow   \tableau{
\CB^{111}\cr
\CB^{111}&\CB^{111}\cr
\CB^{111}&\CB^{111}&\CB^{111} \cr}
\bigsp
\eee(T)={110} \ssp\longrightarrow  \tableau{
\CB^{110}\cr
\CB^{110}& \CB^{110}\cr
\CB^{110}&\CB^{101}&0 \cr}
 $$ $$
\eee(T)={010} \ssp\longrightarrow   \tableau{
0\cr
\CB^{010}&0 \cr
\CB^{010}&0&0 \cr}
\bigsp
\eee(T)={011} \ssp\longrightarrow   \tableau{
0\cr
\CB^{011}&0\cr
\CB^{011}&\CB^{011}&0 \cr}
 $$ $$
\eee(T)={001} \ssp\longrightarrow   \tableau{
0\cr
0&0 \cr
\CB^{001}&0&0 \cr}
$$ 
Placing a basis $\CB^{\eee_1\eee_2\eee_3}$ in cell $(i',j')$ means 
that in the construction of the basis for the module $\BM_{321/00}$ we are to 
apply each element of $\CB^{\eee_1\eee_2\eee_3}$ as a differential operator on the polynomial
$\DD_{321/i'j'}$. 
We should note  
that a basis  for the module $\BM_{321/00}$ is given by the following collection:
$$
\eqalign{
\CB_{321/00}\ses 
&\flip_{20}\ess\CB(\BM_{32})\sps
\flip_{10}\ess\CB(\BM_{32}+ \BM_{311})\sps
\flip_{11}\ess\CB(\BM_{32}\cap \BM_{311})\cr
&
\hskip .3in
\sps\flip_{00}\ess\CB(\BM_{32}+ \BM_{311}+ \BM_{221})\sps
\flip_{01}\Bigl(\,\CB(\BM_{32}\cap \BM_{311} )+ \CB^{101}+\CB^{011}\, \Bigr)\cr 
&
\hskip .6in
\sps \flip_{02}\ess\CB(\BM_{32}\cap \BM_{311}\cap \BM_{221})\ess .\cr
}
$$
Here we have used the symbol $\CB(\BM)$ to denote a bihomogeneous basis
for a module $\BM$. If we compare this result with the developments in Section 2 of [1] we see
that although the algorithm described there involved a recursive process rather than
an assignment of bases to cells, the results are identical. This is in fact a  theorem
that we shall soon establish. But before we get into that, it will be good to 
see how the present construction leads to a representation theoretical 
explanation of the crucial identity.  
\sas

We should note that Conjecture 2.1 may also be stated in a manner which interchanges the roles of
$x$ and $y$. This ``dual'' version requires that the subdiagram $D_{ij}(T)$ be replaced
by one obtained by dividing the shadow of $(i,j)$
into $m$ rectangles, of heights  $v_1,v_2,\ldots , v_m$ from bottom to top,
 by drawing horizontal lines from each of the corners of $\mu$,
then deleting the $r^{th}$ rectangle if $\eee_r=0$ and vertically dropping the remaining rectangles to fill
the gaps. This given, any of the constructions and proofs  that follow have dual versions
which can be routinely derived from their counterparts. We leave it to the reader to fill the gaps 
that result from our systematically dealing with only one of the versions.
With this proviso our basic result here may be stated as follows.
\sa

\heading{\bol Theorem 2.1}

{\ita Let $(i,j)\in \mu$ and $|S_{ij}|=m$. Then on the validity of Conjecture 2.1, the following
are isomorphic as bigraded $S_n$-modules to $\BK_{ij}^x \ssp $, $ \BK_{ij}^y\ssp $, $\BA_{ij}^x$
and $\BA_{ij}^y$ respectively:
$$
\BTK^x_{ij} \ses
\bigoplus_{\multi{\epsilon_1\cdots\epsilon_m  \cr \epsilon_m=1}}
\bigoplus_{\multi{(i,j')\in\mu \cr j' \;\ge\; j}}
\ssp\BM_{S_{ij}}^{T(\epsilon)}(\del)
\DD_{i+\epsilon_cv_c+\cdots+\epsilon_mv_m-1,j'}
\eqno 2.11
$$
$$
\BTK^y_{ij} \ses
\bigoplus_{\multi{\epsilon_1\cdots\epsilon_m \cr \epsilon_1=1}}
\bigoplus_{\multi{(i',j)\in\mu \cr i' \;\ge\; i}}
\BM_{S_{ij}}^{T(\epsilon)}(\del)
\DD_{i',j+\epsilon_1w_1+\cdots+\epsilon_rw_r-1}
\eqno 2.12
$$
$$
\BTA^x_{ij} \ses
\bigoplus_{\epsilon_1\cdots\epsilon_m \st \epsilon_m=1}
\BM_{S_{ij}}^{T(\epsilon)}(\del)
\DD_{i+\epsilon_1v_1+\cdots+\epsilon_mv_m-1,j}
\eqno 2.13
$$
$$
\BTA^y_{ij} \ses
\bigoplus_{\epsilon_1\cdots\epsilon_m \st \epsilon_1=1}
 \BM_{S_{ij}}^{T(\epsilon)}(\del)
\DD_{i,j+\epsilon_1w_1+\cdots+\epsilon_mw_m-1}
\eqno 2.14
$$
where $r$ in 2.12 is determined so that within each term, the lowest corner weakly above $(i',j')$
is the $r ^{th}$  (that is,  $l'_{r+1}<i'-i\le l'_r$), and in 2.11, the leftmost corner
weakly right of $(i',j')$ is the $c^{th}$ (that is, $a'_{c-1}<j'-j\le a'_c$).
}

\heading{\bol Proof}

We shall prove the relations for $\BTK^y_{ij}$ and $\BTA^y_{ij}$.
The relations for the other two are proved similarly by means of the dual version
of 2.6.

The kernel $\BK^y_{ij}$ is isomorphic to 
$$
\BTK_{ij}^y\ses\BM_{ij}/D_y^{-1}(\BM_{i,j+1})
\eqno 2.15
$$
where $D_y^{-1}(\BM_{i,j+1})$ denotes any submodule of $\BM_{ij}$ whatsoever
that is in one-to-one correspondence with $\BM_{i,j+1}$ via $D_y$.
We shall choose a preimage obtained by shifting each contribution to 2.6
one cell to the left, noting that 
$$
D_y\BM_{S_{ij}}^T(\del)\DD_{\mu/i'j'}\ses \BM_{S_{ij}}^T(\del)\DD_{\mu/i',j'+1}
\ess .
$$
This given we may set 
$$
D_y^{-1}(\BM_{i,j+1})  
\ses \bigoplus_{\emptyset\neq T\con S_{i,j+1}}\ess 
\bigoplus_{(i',j')\in D_{i,j+1}(T)}\ess \BM_{S_{i,j+1}}^T(\del)\DD_{\mu/i',j'-1}\ess .
\eqno 2.16
$$
For simplicity we shall only deal with the case when the shadows of $(i,j)$ and $(i,j+1)$
contain the same corners of $\mu$. In this case we may set $S_{i,j+1}=S_{ij}$ in 2.16 and obtain
$$
D_y^{-1}(\BM_{i,j+1})  
\ses \bigoplus_{\emptyset\neq T\con S_{ij }}\ess 
\bigoplus_{(i',j')\in D_{i,j+1 }(T)}\ess \BM_{S_{ij }}^T(\del)\DD_{\mu/i',j'-1}\ess .
$$
This may also be rewritten in the form
$$
D_y^{-1}(\BM_{i,j+1})  
\ses \bigoplus_{\emptyset\neq T\con S_{ij }}\ess 
\bigoplus_{(i',j')\in  D_{i,j+1}^\leftarrow(T)}\ess \BM_{S_{ij }}^T(\del)\DD_{\mu/i'j'}
	\ess ,
\eqno 2.17
$$
where the symbol ``$ D_{i,j+1}^\leftarrow(T)$'' is to represent
the subdiagram of $\mu$ obtained by shifting all cells of  $D_{i,j+1}(T)$
one unit to the left. Now note that when $\eee_1(T)=0$, the diagram $D_{i,j+1}(T)$ 
is identical in shape with $D_{ij}(T)$ but shifted one column to the right.
Thus in this case $ D_{i,j+1}^ \leftarrow(T) =D_{i,j}(T)$. On the other hand,
when $\eee_1(T)=1$ then $D_{i,j+1}(T)$ is $D_{i,j}(T)$ with the leftmost column removed
and then the difference $D_{i,j}(T)- D_{i,j+1}^\leftarrow(T) $ is obtained by picking 
the rightmost cell from each of the rows of $D_{i,j}(T)$. In any case, in view of 2.15, we may write
$$
\BTK_{ij}^y\ses 
\bigoplus_{\emptyset\neq T\con S_{ij }}\ess 
\bigoplus_{(i',j')\in D_{ij}(T)- D_{i,j+1}^\leftarrow(T)}\ess \BM_{S_{ij }}^T(\del)\DD_{\mu/i'j'}\ess ,
\eqno 2.18
$$
which is easily seen to be another way of writing 2.12.

To  prove 2.14 we shall assume that the shadows of $(i,j),(i+1,j),(i,j+1)$ and $(i+1,j+1)$
contain the same corners of $\mu$, so that in this case we can also write
$$
\BTK_{i+1,j}^y\ses 
\bigoplus_{\emptyset\neq T\con S_{ij }}\ess 
\bigoplus_{(i',j')\in D_{i+1,j}(T)- D_{i+1,j+1}^\leftarrow(T)}\ess \BM_{S_{ij }}^T(\del)\DD_{\mu/i'j'}\ess .
$$
Moreover, we see that under these assumptions the diagrams $D_{i+1,j}(T)$ and $ D_{i+1,j+1}^\leftarrow(T) $
are simply $D_{i,j}(T)$ and $ D_{i,j+1}^\leftarrow(T) $ with the bottom row removed, thus the difference
$$
SED_{ij}(T)=\big(D_{i ,j}(T)- D_{i ,j+1}^\leftarrow(T) \big)\sms \big(D_{i+1,j}(T)- D_{i+1,j+1}^\leftarrow(T) \big)
$$
reduces to the southeast corner cell of $D_{i ,j}(T)$ when $\eee_1=1$ and is otherwise empty
when $\eee_1=0$. In any case we may write
$$
\BTA_{ij}^y\ses \BTK_{1j}^y/\BTK_{i+1,j}^y\ses 
\bigoplus_{\emptyset\neq T\con S_{ij }}\ess 
\bigoplus_{(i',j')\in SED_{ij}(T)}\ess \BM_{S_{ij }}^T(\del)\DD_{\mu/i'j'}\ess .
\eqno 2.19
$$
Since when $\eee_1=1$, we have $SED_{ij}(T) \ses \bigl\{(i,\eee_1w_1+\cdots +\eee_mw_m-1)\bigr\}\ssp$,
we see that 2.19 is just another way of writing 2.14.
\sas

The cases we have omitted here are a bit more tedious to deal with if we stick with the convention of
making the set $S_{ij}$ vary with $(i,j)$. A way to deal with all cases at the same time is to 
fix $S=\{\aaa^{(1)},\ldots, \aaa^{(m)}\}$ to be a set of predecessors of $\mu$ obtained by removing some consecutive 
corners from left to right. Suppose that corners $\mu/\aaa^{(b)},\ldots ,\mu/\aaa^{(c)}$ are the ones
in the shadow of $(i,j)$. In 2.6 we would have $\BM_{S_{ij}}^{T(\eee)}$ with (up to renumbering)
$\eee=\eee_b\cdots\eee_c$; however, this decomposes further into the sum of $\BM_S^{T(\eee_1\cdots \eee_m)}$
where $\eee_1\cdots \eee_{b-1}$ and $\eee_{c+1}\cdots \eee_{m}$ vary freely. Setting
$w_s=0$ and $v_s=0$ for each corner $s$ where $\mu/\aaa^{(s)}$ is not in the shadow of $(i,j)$,
the only dependence on $T(\eee)$ in our construction is on $\eee_b,\ldots,\eee_c$. 
In particular, if we use the same set $S$ in our decompositions of $\BM_{ij}$, $\BM_{i+1,j}$, $\BM_{i,j+1}$, 
and $\BM_{i+1,j+1}$, the the above reasoning works even in  the omitted cases. 
\sas

In the figure below we have illustrated 2.12 and 2.14 in the case
$$
\mu=(27^2,25^5,20^2,16^2,12^3,9^4,3^3)
\ess \scs \ess\ess
(i,j)=(4,5) 
\ess \scs \ess\ess
\epsilon=10011
\ess \scs \ess\ess
m=5\ess .
$$ 
Here the vertical rectangles in bold lines give $D_{ij}(T(\eee))$, and the drawn individual cells along the righthand edge
give the contribution of $T(\eee)$ to $\BK_{ij}^y$ with the lowest giving the contribution
to $\BTA_{ij}^y$. 
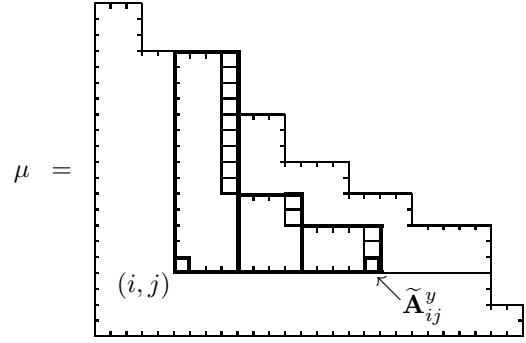
$$
\hbox{%
        \myput(0,200){\llap{$\mu\;\;=\;\;\;$}}
        \hseg(0,0;540)                                  
        \vseg(540,0;40)
        \hseg(500,40;40)
        \vseg(500,40;100)\hseg(400,140;100)
        \vseg(500,40;80)
        \vseg(400,140;40)\hseg(320,180;80)
        \vseg(320,180;40)\hseg(240,220;80)
        \vseg(240,220;60)\hseg(180,280;60)
        \vseg(180,280;80)\hseg(100,360;80)
        \hseg(60,360;40)
        \vseg(60,360;60)
        \hseg(0,420;60)
        \vseg(0,0;420)

        \hseg(100,80;400)
        \vseg(100,80;280)

        \bcorn(120,100)
        \labLL(120,100;{(i,j)})

        \brectangle(180,360)(80,280)
        \brectangle(260,180)(80,100)
        \brectangle(360,140)(100,60)

        \multiputline(360,140)(0,-20){3}{\corn(0,0)}
        \multiputline(260,180)(0,-20){2}{\corn(0,0)}
        \multiputline(180,360)(0,-20){9}{\corn(0,0)}
        \bcorn(360,100)


        \myput(355,75){\vbox to0pt{\hbox{$\raise.5\baselineskip\hbox{$\nwarrow$}
{\bf\widetilde A}^y_{ij}$}\vss}}

        \lhair(540,0){3}
        \dhair(500,40){2}
        \lhair(500,40){6}
        \dhair(420,140){4}      
        \lhair(400,140){3}
        \dhair(320,180){4}
        \lhair(320,180){3}
        \dhair(240,220){4}
        \lhair(240,220){4}
        \dhair(180,280){4}
        \lhair(180,280){4}
        \dhair(60,360){7}
        \lhair(60,360){3}
        \dhair(0,420){3}
        \rhair(0,0){21}
        \uhair(0,0){27}

        \lhair(360,80){3}
        \dhair(260,137){5}
        \lhair(260,140){2}
        \dhair(180,177){4}
        \lhair(180,180){9}
        \dhair(100,357){4}
        \rhair(103,80){15}
        \uhair(100,83){13}

}
\hskip 540\mylength
$$
\sa

To proceed we need the following identity satisfied by the characteristics defined by I.35.
\sas

\heading{\bol Proposition 2.1}
$$
\da\ssp \phi_S^{(k)}\ses {\phi_S^{(m+1-k)}\over \prod_{\bbb\in S} T_\aaa  }
$$ 
{\bol Proof}

Combining I.35 and I.36 we derive that
$$
\phi^{(k)}_S= \sum_{\aaa\in S}\ssp \Bigl(\prod_{\bbb\in S/\{\aaa\}}{1\over 1-T_\aaa/T_\bbb}\ssp \Bigr)\ssp 
(-\nabla)^{m-k}\TH_\aaa
= \sum_{\aaa\in S}\ssp \Bigl(\prod_{\bbb\in S/\{\aaa\}}{1\over 1-T_\aaa/T_\bbb}\ssp \Bigr)\ssp 
(-T_{\aaa})^{m-k}
\TH_\aaa\ess .
$$
Thus (since $\ess \da\TH_\aaa=\TH_\aaa/T_\aaa$):
$$
\eqalign{
\da\phi^{(k)}_S &= \sum_{\aaa\in S}\ssp \Bigl(\prod_{\bbb\in S/\{\aaa\}}{T_\aaa\over T_\aaa-T_\bbb}\ssp \Bigr)\ssp 
\Big({-1\over T_{\aaa}}\Big)^{m-k}\ssp \TH_\aaa\cr
&= \sum_{\aaa\in S}\ssp \Bigl(\prod_{\bbb\in S/\{\aaa\}}{1\over 1-T_\aaa/T_\bbb}\ssp \Bigr)\ssp 
{({-T_{\aaa}})^{k-1}\over\prod_{\bbb\in S}T_\bbb }\ssp \TH_\aaa\ses {\phi^{(m+1-k)}_S\over\prod_{\bbb\in S}T_\bbb } 
\ess .\cr
}
$$
The last equality results from I.35 with $k$ replaced by $\ess m+1-k$.
\bigsp {\bf Q.E.D.}
\sap

We are now ready to show that both 2.1 and 2.2  
may be derived from  geometric properties of lattice diagrams.
 To see how this comes about, for given $0,1$-words $\eee=\eee_1\cdots \eee_m$
and $\eta=\eta_1\cdots \eta_m$ set
$$
\BTA^x_{ij}(\eee) \ses
\BM_{S_{ij}}^{T(\epsilon)}(\del)
\DD_{\mu/i+\epsilon_1v_1+\cdots+\epsilon_mv_m-1,j}
\eqno 2.20
$$
and
$$
\BTA^y_{ij} (\eta)\ses
 \BM_{S_{ij}}^{T(\eta)}(\del)
\DD_{\mu/i,j+\eta_1w_1+\cdots+\eta_mw_m-1}\ess .
\eqno 2.21
$$
We see from 2.13 and 2.14 that
$$
\BTA^x_{ij} =
\bigoplus_{\eee_1\cdots\eee_m \st \epsilon_m=1}
 \BTA^x_{ij}(\eee)
\ess\ess\ess\ess {\rm and }\ess\ess\ess \ess
\BTA^y_{ij} =
\bigoplus_{\eta_1\cdots\eta_m \st \eta_1=1}
\BTA^y_{ij} (\eta)
	\ess .
$$
This given we have the following
refinements of the crucial and flip 
identities. 
\sa

\heading{\bol Theorem 2.2}

{
\ita
For $\eee=(\eee_1,\ldots,\eee_{m-1},1)$ and $\eta=(1,\eee_1,\ldots ,\eee_{m-1})$  we have
(with the same $l$ and $a$ as in 2.1)
$$
t^l\ssp \Fch \BTA^x_{ij}(\eee)\ses q^a\ssp \Fch \BTA^y_{ij}(\eta) 
	\ess ,
\eqno 2.22
$$
while for $\eta=(1,1-\eee_1,\ldots ,1-\eee_{m-1})$ we have
$$
\Fch \BTA^x_{ij}(\eee)\ses T_{\mu/ij}\ssp \da \Fch \BTA^y_{ij}(\eta)
	\ess .
\eqno 2.23
$$
}
\heading{\bol Proof}

We first determine the Frobenius characteristic of $\BTA^x_{ij}(\eee)$,
and then that of $\BTA^y_{ij}(\eta)$.
Let $\eta=(1,\eee_1,\ldots ,\eee_{m-1})$.
Set $\eee_1+\cdots +\eee_m=k$ and $V(\eee)=\epsilon_1v_1+\cdots+\epsilon_mv_m-1$. Then 
from 2.20 we get that
$$
\Fch \BTA^x_{ij}(\eee)\ses {T_\mu\over  t^iq^j t^{V(\eee) }}\da\ssp \Fch \BM_{S_{ij}}^{T(\eee)}
\ess .
\eqno 2.24
$$
Setting $\nu^{(i_s)}=\aaa^{(s)}$, that is 
$$
S_{ij}\ses \bigl\{\ssp \aaa^{(1)},\aaa^{(2)},\ldots ,\aaa^{(m)} \bigr\} \ess,
$$
the definition I.34 gives  
$$
\Fch \BM_{S_{ij}}^{T(\eee)}\ses {\phi^{(k)}_{S_{ij}}\over \prod_{s=1}^mT_{\aaa^{(s)}}^{1-\eee_s} }\ess ,
$$
and  Proposition 2.1 yields
$$
\da  \Fch \BM_{S_{ij}}^{T(\eee)}\ses
{\phi^{(m+1-k)}_{S_{ij}}\over \prod_{s=1}^mT_{\aaa^{(s)}}^{\eee_s } }\ess .
$$
This reduces 2.24 to
$$
\Fch \BTA^x_{ij}(\eee)\ses {T_\mu\over  t^iq^j t^{V(\eee) }} \ssp
{\phi^{(m+1-k)}_{S_{ij}}\over \prod_{s=1}^mT_{\aaa^{(s)}}^{\eee_s } }\ess .
$$
Recalling the definition of $V(\eee)$, we may write
$$
\Fch \BTA^x_{ij}(\eee)\ses {T_\mu\over  t^{i-1}q^j}   \ssp
{\phi^{(m+1-k)}_{S_{ij}}\over \prod_{s=1}^m\bigl(T_{\aaa^{(s)}}t^{v_s}\bigr)^{\eee_s } }\ess .
\eqno 2.25
$$
In a similar way, for $\eta_1+\cdots +\eta_m=k$,  we derive that    
$$
\Fch \BTA^y_{ij}(\eta)\ses {T_\mu\over  t^{i }q^{j-1}}   \ssp
{\phi^{(m+1-k)}_{S_{ij}}\over \prod_{s=1}^m\bigl(T_{\aaa^{(s)}}q^{w_s}\bigr)^{\eta_s } }\ess .
\eqno 2.26
$$
In conclusion, 2.25 and 2.26 yield the identity
$$
\Bigl({1\over t}\ssp \prod_{s=1}^m\bigl(T_{\aaa^{(s)}}t^{v_s}\bigr)^{\eee_s }\Bigr)\ssp \Fch \BTA^x_{ij}(\eee)\ses
\Bigl({1\over q}\ssp \prod_{s=1}^m\bigl(T_{\aaa^{(s)}}t^{w_s}\bigr)^{\eta_s }\Bigr) \ssp \Fch \BTA^y_{ij}(\eta)
\eqno 2.27
$$
To see that this is 2.22, note that the definition of the coefficients $v_s$ and $w_s$ gives
$$
u^{ij}_s= x^{ij}_s/t^{v_s}=x^{ij}_{s+1}/q^{w_s}
	\ess ,
\eqno 2.28
$$ 
and since $T_{\aaa^{(s)}}=T_{\mu/ij}/x_s^{ij}$ we may write (using $\eee_m=1$)
$$
{1\over t}\ssp \prod_{s=1}^m\bigl(T_{\aaa^{(s)}}t^{v_s}\bigr)^{\eee_s }
\ses 
{1\over t}\ssp \prod_{s=1}^m\left({T_{\mu/ij}\over x_s^{ij}}\ssp t^{v_s}\right)^{\eee_s }
\ses
{T_{\mu/ij}^k\over t\ssp u_m^{ij}}\ess \prod_{s=1}^{m-1}\left({1 \over  u_s^{ij}}\right)^{\eee_s }\ess .
\eqno 2.29
$$
Similarly when $\eta_1=1$ we obtain that
$$
{1\over q}\ssp \prod_{s=1}^m\left(T_{\aaa^{(s)}}q^{w_s}\right)^{\eee_s }
\ses 
{T_{\mu/ij}^k\over q\ssp u_0^{ij}}\ess \prod_{s=2}^{m }\left({1 \over  u_s^{ij}}\right)^{\eta_s }\ess .
$$
Thus, taking account of 2.28, we see that when $ \eta_{s+1}=\eee_s$ for  $1\leq s\leq m-1 $,  this last 
expression may be written in the form
$$
{T_{\mu/ij}^k\over q\ssp u_0^{ij}}\ess  \prod_{s=1}^{m-1}\left({1 \over  u_s^{ij}}\right)^{\eee_s }\ess .
$$
Comparing with 2.29 we finally derive that, after  making the approriate cancellations,
2.27 reduces to
$$
{1\over t\,u_m^{ij}}\ess   \Fch \BTA^x_{ij}(\eee)\ses {1\over q\,u_0^{ij}}\ess  \Fch \BTA^y_{ij}(\eta)
$$
which is another way of writing 2.22 since $t\,u_m^{ij}=q^a$ and  $q\,u_0^{ij}=t^l$.
\sas
Next let us assume that $\eta=(1,1-\eee_1,\ldots,1-\eee_{m-1})$.                            
 Since $\eee_1+\cdots +\eee_{m-1}=k-1$,
this choice gives $\ess \sum_{s=1}^m\eta_s=m+1-k$ and in this case $A_{ij}^y(\eta)$ is given by
2.26 with $k$ replaced by $m+1-k$. Thus we may write
$$
A_{ij}^y(\eta)\ses q\ssp T_{\mu/ij}\ess  
{ \phi_{S_{ij}}^{(k)}\over \prod_{s=1}^m\bigl( T_{\aaa^{(s)}}q^{w_s}\bigr)^{\eta_s} }\ess .
$$
Now using 2.20 again we obtain
$$
\eqalign{
\da A_{ij}^y(\eta) 
&\ses {1\over q\ssp T_{\mu/ij}} {\phi_{S_{ij}}^{(m+1-k)}  \over \prod_{s=1}^m T_{\aaa^{(s)}}}
\ssp \prod_{s=1}^m\bigl( T_{\aaa^{(s)}}q^{w_s}\bigr)^{\eta_s}\cr
&\ses
{1\over q\ssp T_{\mu/ij}} {\phi_{S_{ij}}^{(m+1-k)}
\ess \prod_{s=1}^m q^{w_s}  \over \prod_{s=1}^m\bigl( T_{\aaa^{(s)}}q^{w_s}\bigr) }
\ssp \prod_{s=1}^m\bigl( T_{\aaa^{(s)}}q^{w_s}\bigr)^{\eta_s} \cr
&\ses
{1\over q\ssp T_{\mu/ij}} {\phi_{S_{ij}}^{(m+1-k)}\ess \prod_{s=1}^m q^{w_s}  
\over \prod_{s=1}^m\bigl( T_{\aaa^{(s)}}q^{w_s}\bigr) ^{1-\eta_s}}
  \cr
}
\eqno 2.30 
$$
Taking account of 2.28 and recalling that here $\eta=(1,1-\eee_1,\ldots ,1-\eee_{m-1})$, we see
that we have
$$
T_{\mu/ij}\prod_{s=1}^m\bigl( T_{\aaa^{(s)}}q^{w_s}\bigr) ^{1-\eta_s}
\ses
T_{\mu/ij}^k \prod_{s=2}^m \left({  q^{w_s}\over x_s^{ij} }\right)^{\eee_{s-1}}
=\ess
T_{\mu/ij}^k \prod_{s=1}^{m-1} \left({  1\over u_{s }^{ij} }\right)^{\eee_{s }}\ess .
$$
Thus since $\sum_{s=1}^mw_s=a+1$, the identity in 2.30 reduces to
$$
\da A_{ij}^y(\eta)\ses {q^a\ssp \phi_{S_{ij}}^{(m+1-k)}\over  \ssp T_{\mu/ij}^k} \ess 
\prod_{s=1}^{m-1} \bigl( u_{s }^{ij}  \bigr)^{\eee_{s }}\ess .
\eqno 2.31
$$
On the other hand, using 2.28 again, we  can also rewrite 2.25 in the form 
$$
 A_{ij}^x(\eee)\ses {t\ssp T_{\mu/ij}\ssp \phi_{S_{ij}}^{(m+1-k)}\ssp u_m^{ij}\over T_{\mu/ij}^k} \ess 
\prod_{s=1}^{m-1} \bigl( u_{s }^{ij}  \bigr)^{\eee_{s }}\ess .
$$
Comparing with 2.31, we see that 
$$
A_{ij}^x(\eee)\ses {t\ssp u_m^{ij} \over q^a}\ess T_{\mu/ij}\da A_{ij}^y(\eta)
	\ess ,
$$
and this is 2.23 since, as we have seen, $t   u_{s }^{ij}=q^a$.
This completes our proof.
\sa

The refined crucial identity 2.22 and flip identity 2.23 each relate a term
2.20 in the direct sum decomposition 2.13 of the $x$-atom to a term 2.21
in the direct sum decomposition 2.14 of the $y$-atom.  We illustrate this
in the figure that follows in the case
$$\mu=(24^2,22^4,19^3,17^2,15^2,11^4,8^2,6^2,2^2)
	\quad,\quad
	m=7
	\quad,\quad
	(i,j)=(3,4)
	\quad,\quad
	\epsilon=(0,1,0,1,0,0,{\bf1}).
$$

Draw 6 copies of the diagram of
$\mu$ with the shadow of $(i,j)$ marked off.  Put three diagrams on the
right and three on the left, labelled D1--D6, as shown.
In diagram D3, write
${\bf 1},\epsilon_1,\ldots,\epsilon_{m-1},{\bf1}$ just northeast of the
inner corner cells $u^{ij}_0,\ldots,u^{ij}_m$.
Drop vertical lines from each corner to form $m$ vertical rectangles,
and shade the rectangles underneath
$1$'s.  The $\bf 1$ at the bottom right
does not contribute a rectangle since there is
nothing beneath it.  Slide the shaded rectangles to the left to fill in the
gaps, forming the shaded region $D_{ij}(T(\epsilon))$ in D5.
The rightmost cell $(i,j')$
on the bottom row of this region is drawn in, and
gives a term 2.21 of the direct sum 2.14:
${\bf \tilde A}^y_{ij}(\eta)=
	{\bf M}^{T(\eta)}_{S_{ij}}(\partial) \Delta_{\mu/i,j'}$,
where $\eta=(1,\epsilon_1,\ldots,\epsilon_{m-1})$.

Via the refined crucial identity 2.22, this piece of the $y$-atom
corresponds to a piece of the $x$-atom that we locate as follows.
Extend horizontal lines to the left from each corner in D3,
forming $m$ horizontal rectangles.  Shade the rectangles
that are left of
$1$'s.  The $\bf 1$ at the top left
does not contribute a rectangle since there is
nothing to its left.
  Slide the shaded rectangles down to fill in the gaps,
forming the shaded region in D1.
  The topmost
cell $(i',j)$ in the left column of this region is drawn in, and
gives a term 2.20 of the direct sum 2.13:
${\bf \tilde A}^x_{ij}(\epsilon)=
	{\bf M}^{T(\epsilon)}_{S_{ij}}(\partial) \Delta_{\mu/i',j}$.
This term is related to the term from D5 via 2.22.

The three diagrams of $\mu$ on the right side of the figure
illustrate what happens when we apply
${\bf flip}_{\mu/ij}$ to the modules constructed on the left side
of the figure.
Let $\tilde\epsilon_i=1-\epsilon_i$.
In D4, write
${\bf 1},\tilde\epsilon_1,\ldots,\tilde\epsilon_{m-1},{\bf1}$ just northeast
of the inner corner cells, and then shade vertical and horizontal
rectangles according to whether they have a $1$ along the edge at their
end.  This has the effect of complementing which
rectangles are shaded in or not shaded in, except that the vertical
rectangle on the left and the horizontal rectangle on the bottom are
fixed.  Now slide all vertical rectangles left to fill in the gaps,
and place the result in D2.
Its bottom rightmost cell gives a term
${\bf\tilde A}^y_{ij}(\eta)$ of 2.14
for which the refined flip identity 2.23 holds with
$\eta=(1,\tilde\epsilon_1,\ldots,\tilde\epsilon_{m-1}).$
Finally, slide down the horizontal rectangles in D4
to form D6.
Its top left cell gives a term of 2.13 corresponding to the one
in 2.14 from D2 via the crucial identity 2.22 and to the one
in D5 via the flip identity 2.23.

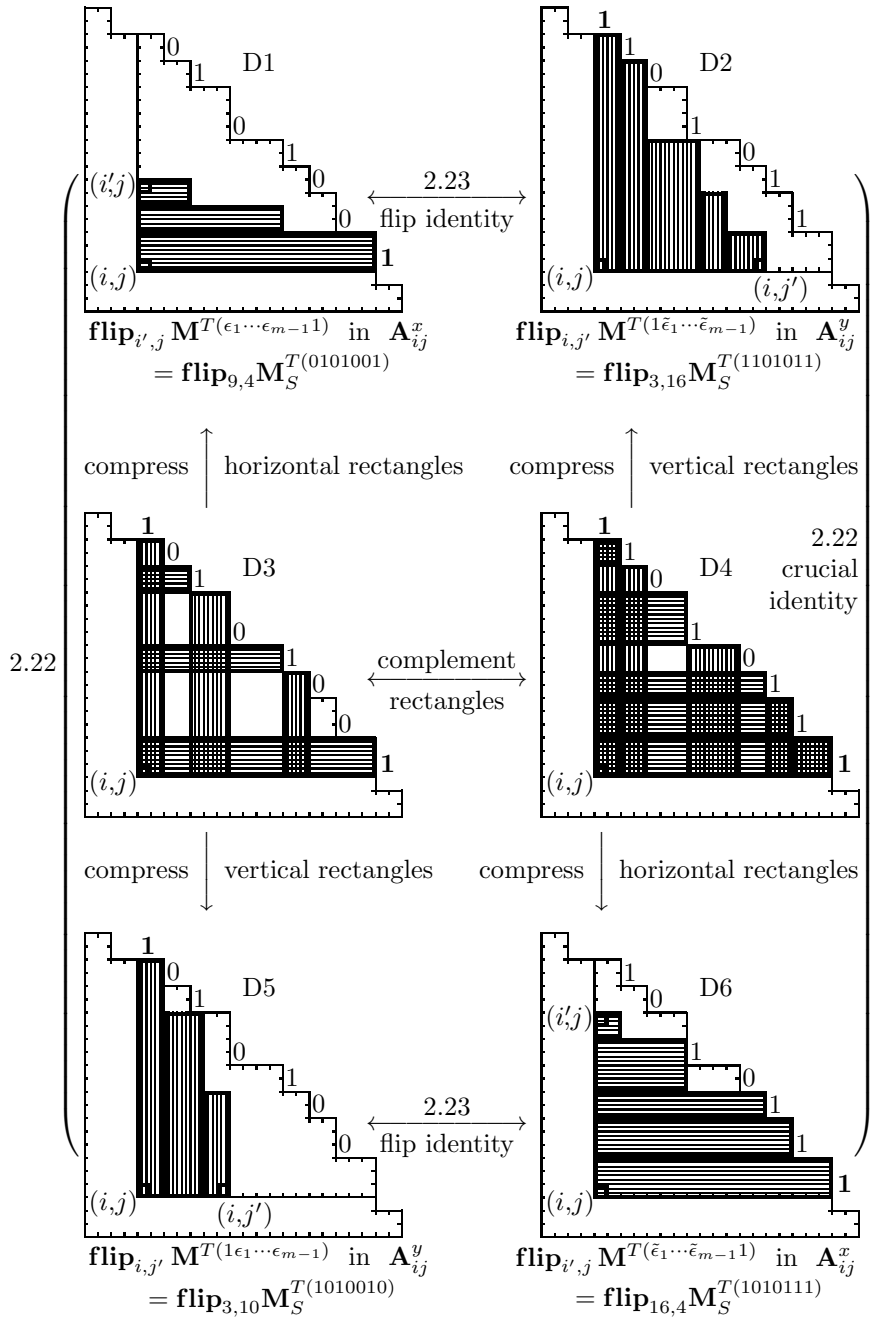
\vfill\eject
\def\bitA{\myput(3,3){0}}
\def\bitB{\myput(3,3){1}}
\def\bitBB{\myput(3,3){\bf1}}

\long\def\AfigI{
\hbox{%
	\vseg(80,60;360)
	\rhair(80,60){18}
	\hseg(80,60;360)
	\uhair(80,60){18}
	\hseg(80,420;40)
	\vseg(120,380;40)
	\dhair(80,420){2}
	\lhair(120,380){2}
	\hseg(120,380;40)
	\vseg(160,340;40)
	\myput(122,385){\bitA}
	\dhair(120,380){2}
	\lhair(160,340){2}
	\hseg(160,340;60)
	\vseg(220,260;80)
	\myput(162,345){\bitB}
	\dhair(160,340){3}
	\lhair(220,260){4}
	\hseg(220,260;80)
	\vseg(300,220;40)
	\myput(222,265){\bitA}
	\dhair(220,260){4}
	\lhair(300,220){2}
	\hseg(300,220;40)
	\vseg(340,180;40)
	\myput(302,225){\bitB}
	\dhair(300,220){2}
	\lhair(340,180){2}
	\hseg(340,180;40)
	\vseg(380,120;60)
	\myput(342,185){\bitA}
	\dhair(340,180){2}
	\lhair(380,120){3}
	\hseg(380,120;60)
	\vseg(440,60;60)
	\myput(382,125){\bitA}
	\dhair(380,120){3}
	\lhair(440,60){3}
	\myput(445,65){\bitBB}
	\hrectangle(440,120)(360,60)
	\Brectangle(440,120)(360,60)
	\hrectangle(300,160)(220,40)
	\Brectangle(300,160)(220,40)
	\hrectangle(160,200)(80,40)
	\Brectangle(160,200)(80,40)
	\bcorn(100,200)
	\vseg(0,0;460)
	\rhair(0,0){23}
	\hseg(0,0;480)
	\uhair(0,0){24}
	\hseg(0,460;40)
	\vseg(40,420;40)
	\dhair(0,460){2}
	\lhair(40,420){2}
	\hseg(40,420;80)
	\vseg(120,380;40)
	\dhair(40,420){4}
	\lhair(120,380){2}
	\hseg(120,380;40)
	\vseg(160,340;40)
	\dhair(120,380){2}
	\lhair(160,340){2}
	\hseg(160,340;60)
	\vseg(220,260;80)
	\dhair(160,340){3}
	\lhair(220,260){4}
	\hseg(220,260;80)
	\vseg(300,220;40)
	\dhair(220,260){4}
	\lhair(300,220){2}
	\hseg(300,220;40)
	\vseg(340,180;40)
	\dhair(300,220){2}
	\lhair(340,180){2}
	\hseg(340,180;40)
	\vseg(380,120;60)
	\dhair(340,180){2}
	\lhair(380,120){3}
	\hseg(380,120;60)
	\vseg(440,40;80)
	\dhair(380,120){3}
	\lhair(440,40){4}
	\hseg(440,40;40)
	\vseg(480,0;40)
	\dhair(440,40){2}
	\lhair(480,0){2}
	\myput(0,-75.0){$\matrix{\flip_{i',j}{\bf M}^{T(\epsilon_1\cdots\epsilon_{m-1}1)}\hbox{~~in~~}{\bf A}^{x}_{ij}\cr\quad={\flip_{9,4}}{\bf M}^{T(0101001)}_S}$}	\myput(7.0,180){$(i'\!\!,\!j)$}
	\myput(7.0,40){$(i,\!j)$}
	\bcorn(100,80)
}\hskip 480\mylength
}

\long\def\AfigII{
\hbox{%
	\vseg(80,60;360)
	\rhair(80,60){18}
	\hseg(80,60;360)
	\uhair(80,60){18}
	\hseg(80,420;40)
	\vseg(120,380;40)
	\myput(82,425){\bitBB}
	\dhair(80,420){2}
	\lhair(120,380){2}
	\hseg(120,380;40)
	\vseg(160,340;40)
	\myput(122,385){\bitB}
	\dhair(120,380){2}
	\lhair(160,340){2}
	\hseg(160,340;60)
	\vseg(220,260;80)
	\myput(162,345){\bitA}
	\dhair(160,340){3}
	\lhair(220,260){4}
	\hseg(220,260;80)
	\vseg(300,220;40)
	\myput(222,265){\bitB}
	\dhair(220,260){4}
	\lhair(300,220){2}
	\hseg(300,220;40)
	\vseg(340,180;40)
	\myput(302,225){\bitA}
	\dhair(300,220){2}
	\lhair(340,180){2}
	\hseg(340,180;40)
	\vseg(380,120;60)
	\myput(342,185){\bitB}
	\dhair(340,180){2}
	\lhair(380,120){3}
	\hseg(380,120;60)
	\vseg(440,60;60)
	\myput(382,125){\bitB}
	\dhair(380,120){3}
	\lhair(440,60){3}
	\vrectangle(120,420)(40,360)
	\Brectangle(120,420)(40,360)
	\vrectangle(160,380)(40,320)
	\Brectangle(160,380)(40,320)
	\vrectangle(240,260)(80,200)
	\Brectangle(240,260)(80,200)
	\vrectangle(280,180)(40,120)
	\Brectangle(280,180)(40,120)
	\vrectangle(340,120)(60,60)
	\Brectangle(340,120)(60,60)
	\bcorn(340,80)
	\vseg(0,0;460)
	\rhair(0,0){23}
	\hseg(0,0;480)
	\uhair(0,0){24}
	\hseg(0,460;40)
	\vseg(40,420;40)
	\dhair(0,460){2}
	\lhair(40,420){2}
	\hseg(40,420;80)
	\vseg(120,380;40)
	\dhair(40,420){4}
	\lhair(120,380){2}
	\hseg(120,380;40)
	\vseg(160,340;40)
	\dhair(120,380){2}
	\lhair(160,340){2}
	\hseg(160,340;60)
	\vseg(220,260;80)
	\dhair(160,340){3}
	\lhair(220,260){4}
	\hseg(220,260;80)
	\vseg(300,220;40)
	\dhair(220,260){4}
	\lhair(300,220){2}
	\hseg(300,220;40)
	\vseg(340,180;40)
	\dhair(300,220){2}
	\lhair(340,180){2}
	\hseg(340,180;40)
	\vseg(380,120;60)
	\dhair(340,180){2}
	\lhair(380,120){3}
	\hseg(380,120;60)
	\vseg(440,40;80)
	\dhair(380,120){3}
	\lhair(440,40){4}
	\hseg(440,40;40)
	\vseg(480,0;40)
	\dhair(440,40){2}
	\lhair(480,0){2}
	\myput(-45,-75.0){$\matrix{\flip_{i,j'}{\bf M}^{T(1\tilde\epsilon_1\cdots\tilde\epsilon_{m-1})}\hbox{~~in~~}{\bf A}^{y}_{ij}\cr\quad={\flip_{3,16}}{\bf M}^{T(1101011)}_S}$}	\myput(320,22.0){$(i,\!j')$}
	\myput(7.0,40){$(i,\!j)$}
	\bcorn(100,80)
}\hskip 480\mylength
}

\long\def\AfigIII{
\hbox{%
	\vseg(80,60;360)
	\rhair(80,60){18}
	\hseg(80,60;360)
	\uhair(80,60){18}
	\hseg(80,420;40)
	\vseg(120,380;40)
	\myput(82,425){\bitBB}
	\dhair(80,420){2}
	\lhair(120,380){2}
	\hseg(120,380;40)
	\vseg(160,340;40)
	\myput(122,385){\bitA}
	\dhair(120,380){2}
	\lhair(160,340){2}
	\hseg(160,340;60)
	\vseg(220,260;80)
	\myput(162,345){\bitB}
	\dhair(160,340){3}
	\lhair(220,260){4}
	\hseg(220,260;80)
	\vseg(300,220;40)
	\myput(222,265){\bitA}
	\dhair(220,260){4}
	\lhair(300,220){2}
	\hseg(300,220;40)
	\vseg(340,180;40)
	\myput(302,225){\bitB}
	\dhair(300,220){2}
	\lhair(340,180){2}
	\hseg(340,180;40)
	\vseg(380,120;60)
	\myput(342,185){\bitA}
	\dhair(340,180){2}
	\lhair(380,120){3}
	\hseg(380,120;60)
	\vseg(440,60;60)
	\myput(382,125){\bitA}
	\dhair(380,120){3}
	\lhair(440,60){3}
	\myput(445,65){\bitBB}
	\hrectangle(440,120)(360,60)
	\Brectangle(440,120)(360,60)
	\hrectangle(300,260)(220,40)
	\Brectangle(300,260)(220,40)
	\hrectangle(160,380)(80,40)
	\Brectangle(160,380)(80,40)
	\vrectangle(120,420)(40,360)
	\Brectangle(120,420)(40,360)
	\vrectangle(220,340)(60,280)
	\Brectangle(220,340)(60,280)
	\vrectangle(340,220)(40,160)
	\Brectangle(340,220)(40,160)
	\vseg(0,0;460)
	\rhair(0,0){23}
	\hseg(0,0;480)
	\uhair(0,0){24}
	\hseg(0,460;40)
	\vseg(40,420;40)
	\dhair(0,460){2}
	\lhair(40,420){2}
	\hseg(40,420;80)
	\vseg(120,380;40)
	\dhair(40,420){4}
	\lhair(120,380){2}
	\hseg(120,380;40)
	\vseg(160,340;40)
	\dhair(120,380){2}
	\lhair(160,340){2}
	\hseg(160,340;60)
	\vseg(220,260;80)
	\dhair(160,340){3}
	\lhair(220,260){4}
	\hseg(220,260;80)
	\vseg(300,220;40)
	\dhair(220,260){4}
	\lhair(300,220){2}
	\hseg(300,220;40)
	\vseg(340,180;40)
	\dhair(300,220){2}
	\lhair(340,180){2}
	\hseg(340,180;40)
	\vseg(380,120;60)
	\dhair(340,180){2}
	\lhair(380,120){3}
	\hseg(380,120;60)
	\vseg(440,40;80)
	\dhair(380,120){3}
	\lhair(440,40){4}
	\hseg(440,40;40)
	\vseg(480,0;40)
	\dhair(440,40){2}
	\lhair(480,0){2}
	\myput(7.0,40){$(i,\!j)$}
	\bcorn(100,80)
}\hskip 480\mylength
}

\long\def\AfigIV{
\hbox{%
	\vseg(80,60;360)
	\rhair(80,60){18}
	\hseg(80,60;360)
	\uhair(80,60){18}
	\hseg(80,420;40)
	\vseg(120,380;40)
	\myput(82,425){\bitBB}
	\dhair(80,420){2}
	\lhair(120,380){2}
	\hseg(120,380;40)
	\vseg(160,340;40)
	\myput(122,385){\bitB}
	\dhair(120,380){2}
	\lhair(160,340){2}
	\hseg(160,340;60)
	\vseg(220,260;80)
	\myput(162,345){\bitA}
	\dhair(160,340){3}
	\lhair(220,260){4}
	\hseg(220,260;80)
	\vseg(300,220;40)
	\myput(222,265){\bitB}
	\dhair(220,260){4}
	\lhair(300,220){2}
	\hseg(300,220;40)
	\vseg(340,180;40)
	\myput(302,225){\bitA}
	\dhair(300,220){2}
	\lhair(340,180){2}
	\hseg(340,180;40)
	\vseg(380,120;60)
	\myput(342,185){\bitB}
	\dhair(340,180){2}
	\lhair(380,120){3}
	\hseg(380,120;60)
	\vseg(440,60;60)
	\myput(382,125){\bitB}
	\dhair(380,120){3}
	\lhair(440,60){3}
	\myput(445,65){\bitBB}
	\hrectangle(440,120)(360,60)
	\Brectangle(440,120)(360,60)
	\hrectangle(380,180)(300,60)
	\Brectangle(380,180)(300,60)
	\hrectangle(340,220)(260,40)
	\Brectangle(340,220)(260,40)
	\hrectangle(220,340)(140,80)
	\Brectangle(220,340)(140,80)
	\hrectangle(120,420)(40,40)
	\Brectangle(120,420)(40,40)
	\vrectangle(120,420)(40,360)
	\Brectangle(120,420)(40,360)
	\vrectangle(160,380)(40,320)
	\Brectangle(160,380)(40,320)
	\vrectangle(300,260)(80,200)
	\Brectangle(300,260)(80,200)
	\vrectangle(380,180)(40,120)
	\Brectangle(380,180)(40,120)
	\vrectangle(440,120)(60,60)
	\Brectangle(440,120)(60,60)
	\vseg(0,0;460)
	\rhair(0,0){23}
	\hseg(0,0;480)
	\uhair(0,0){24}
	\hseg(0,460;40)
	\vseg(40,420;40)
	\dhair(0,460){2}
	\lhair(40,420){2}
	\hseg(40,420;80)
	\vseg(120,380;40)
	\dhair(40,420){4}
	\lhair(120,380){2}
	\hseg(120,380;40)
	\vseg(160,340;40)
	\dhair(120,380){2}
	\lhair(160,340){2}
	\hseg(160,340;60)
	\vseg(220,260;80)
	\dhair(160,340){3}
	\lhair(220,260){4}
	\hseg(220,260;80)
	\vseg(300,220;40)
	\dhair(220,260){4}
	\lhair(300,220){2}
	\hseg(300,220;40)
	\vseg(340,180;40)
	\dhair(300,220){2}
	\lhair(340,180){2}
	\hseg(340,180;40)
	\vseg(380,120;60)
	\dhair(340,180){2}
	\lhair(380,120){3}
	\hseg(380,120;60)
	\vseg(440,40;80)
	\dhair(380,120){3}
	\lhair(440,40){4}
	\hseg(440,40;40)
	\vseg(480,0;40)
	\dhair(440,40){2}
	\lhair(480,0){2}
	\myput(7.0,40){$(i,\!j)$}
	\bcorn(100,80)
}\hskip 480\mylength
}

\long\def\AfigV{
\hbox{%
	\vseg(80,60;360)
	\rhair(80,60){18}
	\hseg(80,60;360)
	\uhair(80,60){18}
	\hseg(80,420;40)
	\vseg(120,380;40)
	\myput(82,425){\bitBB}
	\dhair(80,420){2}
	\lhair(120,380){2}
	\hseg(120,380;40)
	\vseg(160,340;40)
	\myput(122,385){\bitA}
	\dhair(120,380){2}
	\lhair(160,340){2}
	\hseg(160,340;60)
	\vseg(220,260;80)
	\myput(162,345){\bitB}
	\dhair(160,340){3}
	\lhair(220,260){4}
	\hseg(220,260;80)
	\vseg(300,220;40)
	\myput(222,265){\bitA}
	\dhair(220,260){4}
	\lhair(300,220){2}
	\hseg(300,220;40)
	\vseg(340,180;40)
	\myput(302,225){\bitB}
	\dhair(300,220){2}
	\lhair(340,180){2}
	\hseg(340,180;40)
	\vseg(380,120;60)
	\myput(342,185){\bitA}
	\dhair(340,180){2}
	\lhair(380,120){3}
	\hseg(380,120;60)
	\vseg(440,60;60)
	\myput(382,125){\bitA}
	\dhair(380,120){3}
	\lhair(440,60){3}
	\vrectangle(120,420)(40,360)
	\Brectangle(120,420)(40,360)
	\vrectangle(180,340)(60,280)
	\Brectangle(180,340)(60,280)
	\vrectangle(220,220)(40,160)
	\Brectangle(220,220)(40,160)
	\bcorn(220,80)
	\vseg(0,0;460)
	\rhair(0,0){23}
	\hseg(0,0;480)
	\uhair(0,0){24}
	\hseg(0,460;40)
	\vseg(40,420;40)
	\dhair(0,460){2}
	\lhair(40,420){2}
	\hseg(40,420;80)
	\vseg(120,380;40)
	\dhair(40,420){4}
	\lhair(120,380){2}
	\hseg(120,380;40)
	\vseg(160,340;40)
	\dhair(120,380){2}
	\lhair(160,340){2}
	\hseg(160,340;60)
	\vseg(220,260;80)
	\dhair(160,340){3}
	\lhair(220,260){4}
	\hseg(220,260;80)
	\vseg(300,220;40)
	\dhair(220,260){4}
	\lhair(300,220){2}
	\hseg(300,220;40)
	\vseg(340,180;40)
	\dhair(300,220){2}
	\lhair(340,180){2}
	\hseg(340,180;40)
	\vseg(380,120;60)
	\dhair(340,180){2}
	\lhair(380,120){3}
	\hseg(380,120;60)
	\vseg(440,40;80)
	\dhair(380,120){3}
	\lhair(440,40){4}
	\hseg(440,40;40)
	\vseg(480,0;40)
	\dhair(440,40){2}
	\lhair(480,0){2}
	\myput(0,-75.0){$\matrix{\flip_{i,j'}{\bf M}^{T(1\epsilon_1\cdots\epsilon_{m-1})}\hbox{~~in~~}{\bf A}^{y}_{ij}\cr\quad={\flip_{3,10}}{\bf M}^{T(1010010)}_S}$}	\myput(200,22.0){$(i,\!j')$}
	\myput(7.0,40){$(i,\!j)$}
	\bcorn(100,80)
}\hskip 480\mylength
}

\long\def\AfigVI{
\hbox{%
	\vseg(80,60;360)
	\rhair(80,60){18}
	\hseg(80,60;360)
	\uhair(80,60){18}
	\hseg(80,420;40)
	\vseg(120,380;40)
	\dhair(80,420){2}
	\lhair(120,380){2}
	\hseg(120,380;40)
	\vseg(160,340;40)
	\myput(122,385){\bitB}
	\dhair(120,380){2}
	\lhair(160,340){2}
	\hseg(160,340;60)
	\vseg(220,260;80)
	\myput(162,345){\bitA}
	\dhair(160,340){3}
	\lhair(220,260){4}
	\hseg(220,260;80)
	\vseg(300,220;40)
	\myput(222,265){\bitB}
	\dhair(220,260){4}
	\lhair(300,220){2}
	\hseg(300,220;40)
	\vseg(340,180;40)
	\myput(302,225){\bitA}
	\dhair(300,220){2}
	\lhair(340,180){2}
	\hseg(340,180;40)
	\vseg(380,120;60)
	\myput(342,185){\bitB}
	\dhair(340,180){2}
	\lhair(380,120){3}
	\hseg(380,120;60)
	\vseg(440,60;60)
	\myput(382,125){\bitB}
	\dhair(380,120){3}
	\lhair(440,60){3}
	\myput(445,65){\bitBB}
	\hrectangle(440,120)(360,60)
	\Brectangle(440,120)(360,60)
	\hrectangle(380,180)(300,60)
	\Brectangle(380,180)(300,60)
	\hrectangle(340,220)(260,40)
	\Brectangle(340,220)(260,40)
	\hrectangle(220,300)(140,80)
	\Brectangle(220,300)(140,80)
	\hrectangle(120,340)(40,40)
	\Brectangle(120,340)(40,40)
	\bcorn(100,340)
	\vseg(0,0;460)
	\rhair(0,0){23}
	\hseg(0,0;480)
	\uhair(0,0){24}
	\hseg(0,460;40)
	\vseg(40,420;40)
	\dhair(0,460){2}
	\lhair(40,420){2}
	\hseg(40,420;80)
	\vseg(120,380;40)
	\dhair(40,420){4}
	\lhair(120,380){2}
	\hseg(120,380;40)
	\vseg(160,340;40)
	\dhair(120,380){2}
	\lhair(160,340){2}
	\hseg(160,340;60)
	\vseg(220,260;80)
	\dhair(160,340){3}
	\lhair(220,260){4}
	\hseg(220,260;80)
	\vseg(300,220;40)
	\dhair(220,260){4}
	\lhair(300,220){2}
	\hseg(300,220;40)
	\vseg(340,180;40)
	\dhair(300,220){2}
	\lhair(340,180){2}
	\hseg(340,180;40)
	\vseg(380,120;60)
	\dhair(340,180){2}
	\lhair(380,120){3}
	\hseg(380,120;60)
	\vseg(440,40;80)
	\dhair(380,120){3}
	\lhair(440,40){4}
	\hseg(440,40;40)
	\vseg(480,0;40)
	\dhair(440,40){2}
	\lhair(480,0){2}
	\myput(-45,-75.0){$\matrix{\flip_{i',j}{\bf M}^{T(\tilde\epsilon_1\cdots\tilde\epsilon_{m-1}1)}\hbox{~~in~~}{\bf A}^{x}_{ij}\cr\quad={\flip_{16,4}}{\bf M}^{T(1010111)}_S}$}	\myput(7.0,320){$(i'\!\!,\!j)$}
	\myput(7.0,40){$(i,\!j)$}
	\bcorn(100,80)
}\hskip 480\mylength
}

\def\vc#1#2{{\hsize=480\mylength\vcenter{#1}}}
\def\vc#1#2{{\hsize=480\mylength\vcenter{\smash{\hbox{\myput(240,-100){#2}}}\vbox{#1}}}}
\def\svc#1#2{\smash{\vc{#1}{#2}}}
\def\longbiarrow{\mathop{\leftarrow\joinrel\relbar\joinrel\relbar\joinrel\relbar\joinrel\relbar\joinrel\relbar\joinrel\relbar\joinrel\relbar\joinrel\rightarrow}\limits}

\hbox{}\vskip1in
$$\baselineskip 12pt\mylength=.25pt
\hbox{2.22}\!
\left(\matrix{
	\svc{\AfigI}{D1}
	&\!\!\!\!\!\!\!\!\!\!\!\!\!\!\longbiarrow^{\hbox{2.23}}_{\hbox{flip identity}}\!\!\!
	&\svc{\AfigII}{D2}
\cr\noalign{\vskip 1in}
	\rlap{compress $\Bigg\uparrow$ horizontal rectangles} \hfill &&
	\hfill \llap{compress $\Bigg\uparrow$ vertical rectangles}
\cr
	\vc{\AfigIII}{D3}
	&\!\!\!\!\!\!\!\!\!\!\!\!\!\!\longbiarrow^{\hbox{complement}}_{\hbox{rectangles}}\!\!\!
	&\vc{\AfigIV}{D4}
\cr\noalign{\vskip -.1in}
	\rlap{compress $\Bigg\downarrow$ vertical rectangles} \hfill &&
	\hfill \llap{compress $\Bigg\downarrow$ horizontal rectangles}
\cr\noalign{\vskip .95in}
	\svc{\AfigV}{D5}
	&\!\!\!\!\!\!\!\!\!\!\!\!\!\!\longbiarrow^{\hbox{2.23}}_{\hbox{flip identity}}\!\!\!
	&\svc{\AfigVI}{D6}
}\!
\llap{\vbox{\llap{2.22}\llap{crucial}\llap{identity}\vskip.3in}}\!
\right)
$$

\eject

\heading{\bol Remark 2.1}

We should point out that since 2.22 is equivalent to the four term recursion
and the latter in turn implies the expansion in I.16, it follows from 2.6 that
the Frobenius characteristic $C_{\mu/ij}$ of $\BM_{\mu/ij}$ is given by the formula
$$
C_{\mu/ij}\ses {1\over M}{T_{\mu/ij}\over \nabla}\ess 
\biggl(\ssp \prod_{s=0}^m\ssp \Bigl(1-\nabla {u_s^{ij}\over T_{\mu/ij}}\ssp \Bigr)
 \biggr)\ssp \phi_S^{(m)}\ess .
$$
The reader may find it challenging to derive
this identity directly from 2.6 making only use of the fact that
when $S_{ij}=\bigl\{\ssp \aaa^{(1)},\aaa^{(2)},\ldots ,\aaa^{(m)},\ssp \bigr\}\ssp $
and $\ssp \sum_{i=1}^m\eee_i=k $ we have
$$
\Fch \BM_{S_{ij}}^{T(\eee)}\ses {\phi_{S_{ij}}^{(k)}\over \prod_{s=1}^m T_{\aaa^{(s)} }^{1-\eee_s}}
\ses \biggl(\prod_{s=1}^m \Bigl(-{\nabla \over T_{\aaa^{(s)} } }\Bigr)^{1-\eee_s} \biggr)\phi_{S_{ij}}^{(m)}\ess .
$$
\sa

\sas

We terminate this section with a proof that  
the bases for $\BM_\mu$ constructed in [1] by the recursive algorithm
of Bergeron-Haiman, may directly be obtained by the same 
module assignment process we used in 2.6.
We begin with a compact summary of this algorithm; then we
give a direct formula for the final result of the recursion.
As before we set
$$
\Pred (\mu)\ses \bigl\{\nu^{(1)},\nu^{(2)},\ldots ,\nu^{(d)}\bigr\}, 
$$
with the corner cells $\mu/\nu^{(1)},\mu/\nu^{(2)},\ldots ,\mu/\nu^{(d)}  $
ordered from left to right and respective weights
$$
x_1=t^{l_1'}q^{a_1'}\scs x_2=t^{l_2'}q^{a_2'}\scs\ldots\scs  x_d=t^{l_d'}q^{a_d'}\ess . 
$$
The Algorithm is conjectured to
produce a basis for $\BM_\mu$ from bases of
$\BM_{\nu^{(1)}},\ldots,\BM_{\nu^{(d)}}$.  We abbreviate $\BM_{\nu^{(r)}}$ as
$\BM_r$, and we work with the ``Science Fiction Conjecture'' that
$\BM_1,\ldots,\BM_d$ generate a distributive lattice under span and
intersection.
\sas

In [1] the algorithm assigns a module $\BB_{ij}$ to each cell $(i,j)$ of $\mu$ by 
a process that starts from the top row then proceeds down one row at the time
ending at first row. For notational convenience  we shall also assign modules 
here to cells left or right of $\mu$ in the strip
$0\le i<l(\mu)$,\ess $-\infty<j<\infty$, according to the following recipe:
$$
\BB_{ij} \ses \cases{\sum_{s=1}^d\BM_s & if $j<0$\ssp ,\cr
\cr
\{0\} \ess\ess\ess & if $  j\geq \mu_{i+1}\ssp$.
\cr}
\eqno 2.32
$$
The algorithm starts with setting
$$
\BB_{ij}\ses \BM_1 \bigsp \forall \ess\ess\ess (i,j)\ess\ess \hbox {in the top row of $\mu$}\ess ,
\eqno 2.33
$$
this given, for all lower rows, the assignment is
$$
\BB_{ij} \ses \cases{
\BB_{i+1,j} + (\BB_{i+1,j-w}\cap \BM_r)
& if    row $\ssp i+1\ssp$ contains  the $\ssp r^{th}\ssp $ corner of $\mu$\cr
&
 and  $\ess \mu_{i+1}-\mu_{i+2}=w\ess $,  \cr 
    \cr
\BB_{i+1,j}			& if $\ess \mu_{i+1}=\mu_{i+2}\ess $.  \cr
\cr
}
\eqno 2.34
$$
It is conjectured in [1] that,  for any $\mu\part n \ssp $, the module  $\BM_\mu=\CL_\del[\DD_\mu]$ , 
decomposes as the direct sum
$$
	\BM_\mu \ses \bigoplus_{(i,j)\in\mu} \BB_{ij}(\del)\del_{x_{n}}^i\del_{y_{n}}^j\DD_\mu\ssp .
\eqno 2.35
$$
If $\CB_{ij}$ is a basis  of $\BB_{ij}$, then a basis for $\BM_\mu$ should  
be given by the collection

$$
	\CB_\mu \ses \bigcup_{(i,j)\in\mu} \CB_{ij}(\del)\del_{x_{n}}^i\del_{y_{n}}^j\DD_\mu\ess .
\eqno 2.36
$$
Since the distributivity conjecture 
 assures that each $\BB_{ij}$ decomposes
into a direct sum of various components
$\Mepsd =\BM^{T(\epsilon)}_S$
where $S=\Pred(\mu)$ and
$T(\epsilon)=\setof{\nu^{(r)}}{\epsilon_r=1}$, 
we must have direct sum decompositions of the form
$$
\BB_{ij}\ses \bigoplus_{\eee\in{\cal E}_{ij}}\Mepsd
$$
for suitable subsets ${\cal E}_{ij}$. It develops that 
these subsets can be given explicitly by a formula 
which is essentially 2.6 for $\BM_{\mu/00}$.
In point of fact we have put together this formula 
by simply discovering how to place the components $\Mepsd$ directly into
the Young diagram of $\mu\ssp $, bypassing the recursive
process defined by 2.32, 2.33 and 2.34.
To be precise we have

\sas
\heading{\bol Proposition 2.2}

{\ita

Let $w_1=a'_1+1$ and $\ssp w_s=a'_s-a'_{s-1}$ for $\ssp s=2,\ldots,d$.
Assuming the Science Fiction Conjecture,
the Bergeron-Haiman recursion is equivalent to placing $\Mepsd$
in cells $(i,j)$ with $j<\epsilon_1w_1+\cdots+\epsilon_{r }w_{r } $,
where $r$ is the number of corners of $\mu$ that are above 
row $\ssp i+1\ssp $.   In symbols,
$$
\BB_{ij} \ses \bigoplus_{\multi{
	\epsilon_1\cdots\epsilon_d	\cr
	j \;<\; \epsilon_1w_1+\,\cdots\,+\epsilon_rw_r}}
		\Mepsd 
\eqno  2.37
$$
where $\epsilon_1,\cdots,\epsilon_d$ independently run over
$\{0,1\}$ in all ways with at least one of them being nonzero.
}
\sas

\heading{\bol Proof }

For convenience, we define
$$
W_r(\epsilon) \ses \epsilon_1 w_1 + \epsilon_2 w_2 + \cdots + \epsilon_r w_r.
\eqno 2.38
$$
We shall work our way from the top row of the partition down, to establish
that the $\BB_{ij}$ as given by 2.37 satisfy the Bergeron-Haiman recursion.
We start by  checking  the definition
of $\BB_{ij}$ for cells external to $\mu$.
Noting that $0\le W_r(\epsilon)\le w_1+\cdots+w_r=\mu_{i+1}$,
2.37 states that $\BB_{ij}$ is the span of all $\Mepsd$'s when
$j<0$ and is $\{0\}$ when $j\ge\mu_{i+1}$, in agreement with 2.32.

On the top row of the partition, we have $r=1$, $w_1=\mu_{i+1}$, and
$$W_r(\epsilon) \ses
\cases{	\mu_{i+1}	& if $\epsilon_1=1$	\cr
	0		& if $\epsilon_1=0$,}
$$
so that in 2.37, $\BB_{ij}$ is the span of all $\Mepsd$'s for which
$\epsilon_1=1$; and this is just $\BM_1$, agreeing with 2.33.

On any subsequent row that does not contain a corner, we have
$\mu_{i+1}=\mu_{i+2}$, and 2.34 gives
$\BB_{ij}=\BB_{i+1,j}$; at the same time in this case  we must use the same $r$ in 2.37 for
rows $i$ and $i+1$, and this gives $\BB_{ij}=\BB_{i+1,j}$, as desired.

Finally, consider the row containing the $r^{th}$ corner, that is,
$i<l(\mu)-1$ with $a'_r+1=\mu_{i+1}>\mu_{i+2}=a'_{r-1}+1$
and thus $\mu_{i+1}-\mu_{i+2}=w_r$. In other words we must take $w=w_r$
in the first case of 2.34.
Now according to 2.37 we have $\Mepsd\con \BB_{ij}\ssp $ if and only if
$\ssp j<W_r(\eee)=W_{r-1}(\eee)+\eee_rw_r$. On the other hand if we assume 
inductively, that both $\BB_{i+1,j}$ and $\BB_{i+1,j-w}$
are given by 2.37, then we have 
$$ 
 \Mepsd\con \BB_{i+1,j} + (\BB_{i+1,j-w}\cap \BM_r)
$$
 if and only if either

{\item{(a)}
$j<W_{r-1}(\epsilon)$, or
\item{(b)}
$j-w_r<W_{r-1}(\epsilon)$ and $\epsilon_r=1$.
}

\noindent
When $\epsilon_r=0$, (b) is false, while (a)
is equivalent to $j<W_r(\epsilon)$ because
$W_r(\epsilon)=W_{r-1}(\epsilon)+\epsilon_rw_r=W_{r-1}(\epsilon)+0$.
When $\epsilon_r=1$, (a) is equivalent to
$j<W_{r-1}(\epsilon)$, and (b) to $j<W_{r-1}(\epsilon)+w_r=W_r(\epsilon)$,
so when (a) holds, so does (b).  In total, (a) or (b) holds when
$j<W_r(\epsilon)$. This assures the equality 
$$
\BB_{ij} \ses \BB_{i+1,j} + (\BB_{i+1,j-w}\cap \BM_r)
$$
in this case and completes our proof that the assignment in 
2.37 satisfies the Bergeron-Haiman recursion.

\vfill\supereject

\heading{\bol 3. Some examples. }

In this section we shall illustrate the theory we have developed by applying it to 
the ``hook'' case $\mu=(n-k,1^k)$ . We shall see that in 
this case all  our predictions go through in the finest detail.
More significantly, this example  gives us a glimpse of the additional ingredients
that are needed to carry out our program in the general case.
\sas

We shall begin with the special case when $\ssp \mu\ssp $ reduces to a column ($\mu=(1^n)$)
or a row ($\mu=(n)$). In either case, bases for $\BM_\mu$ are well known 
(see [1], [6]). For instance in the case $\mu=(1^n)$, $\DD_\mu$ reduces to the Vandermonde determinant
in $\xon$
$$
\DD_{1^n}\ses \det\|\ssp  x_j^{i-1}\ssp \|_{i,j=1}^n\ses \DD_n(x_1,x_2,\ldots, x_n)\ess .
$$
The basic observation here is that we have  
$$
\DD_n(x_1,x_2,\ldots, x_n)\ses x_1^0x_2^1x_3^2\cdots x_n^{n-1}\sps <\ssp \cdots 
$$
where the symbol ``$<\cdots $'' is to mean that the monomials in the omitted terms are all  
greater than the preceding one in the lexicographic order.
This given, it is easily seen that, when  $  \eee_i\leq i-1$,  we also have
$$
\del_{x_1}^{\eee_1}\del_{x_2}^{\eee_2}\cdots \del_{x_n}^{\eee_n} 
\DD_n(x_1,x_2,\ldots, x_n)
\ses c(\eee)\ssp x_1^{0-\eee_1}x_2^{1-\eee_2}x_3^{2-\eee_3}\cdots x_n^{n-1-\eee_n}\sps <\ssp \cdots 
$$
with $c(\eee)$ a nonvanishing constant. This shows that the Vandermonde $\DD_n(\xon)$
has at least $n!$ independent derivatives. Since we know ([7], [10]) that $\dim \BM_\mu\leq n!$
for $\mu\part n$, it follows that the collection
$$
\CB_n(x_1,x_2,\ldots, x_n)\ses \Bigl\{\ssp \del_{x_1}^{\eee_1}\del_{x_2}^{\eee_2}\cdots \del_{x_n}^{\eee_n}
\DD_n(x_1,x_2,\ldots , x_n)\ssp :\ssp 0\leq \eee_i\leq i-1 \ssp \Bigr\}
\eqno 3.1
$$
is a basis for $\BM_{1^n}$.

Of course, we have an analogous result  in the ``row'' case $\mu=(n)$ . In fact, then  we have
$$
\DD_\mu\ses \DD_n(y_1,y_2,\ldots ,y_n)
$$
and thus a basis for $\BM_n$ is given by the collection
$$
\CB_n(y_1,y_2,\ldots ,y_n,)\ses \Bigl\{\ssp \del_{y_1}^{\eee_1}\del_{y_2}^{\eee_2}\cdots \del_{y_n}^{\eee_n}
\DD_n(y_1,y_2,\ldots , y_n)\ssp :\ssp 0\leq \eee_i\leq i-1 \ssp \Bigr\}\ess .
\eqno 3.2
$$
These classical results translate into the following basic facts concerning the modules $\BM_{1^{n+1}/i,0}$
and $\BM_{n+1/0,j}$:
\sa

\heading{\bol Theorem 3.1}

{\ita For each $1\leq i,j\leq n$ we have the following direct sum decompositions:
$$
\BM_{1^{n+1}/i,0}\ses 
\BM_{1^n}(\del)\DD_{1^{n+1}/i,0}
\ssp \oplus\ssp
\BM_{1^n}(\del)\DD_{1^{n+1}/i+1,0}
\ssp \oplus\ssp
\cdots
\ssp \oplus\ssp
\BM_{1^n}(\del)\DD_{1^{n+1}/n+1,0}
\eqno 3.3
$$
$$
\BM_{{n+1}/0,j}\ses 
\BM_{ n }(\del)\DD_{ n+1/0,j}
\ssp \oplus\ssp
\BM_{ n }(\del)\DD_{ n+1/0,j+1}
\ssp \oplus\ssp
\cdots
\ssp \oplus\ssp
\BM_{ n }(\del)\DD_{n+1/ 0,n+1}
	\ess .
\eqno 3.4
$$
Moreover, we can represent their respective atoms by the following modules.
\sas

\itemitem {(1)} For $\mu=(1^{n+1})$,
$$
a)\ess  \BA_{i,0}^x\ses \BM_{1^n}
\ess \scs\ess\ess\ess\ess 
b)\ess  \BA_{i,0}^y\ses \BM_{1^n}(\del)\DD_{1^{n+1}/i,0}\ess .
\eqno 3.5
$$

\itemitem {(2)} For $\mu=(n+1)$,
$$
a)\ess \BA_{i,0}^y\ses \BM_{n}
\ess \scs\ess\ess\ess\ess 
a)\ess  \BA_{0,j}^x\ses \BM_{n}(\del)\DD_{1^{n+1}/0,j}\ess .
\eqno 3.6
$$
}

\heading{\bol Proof}

By Proposition I.1 we have
$$
\sum_{i=1}^{n+1}\ssp \del_{x_i}\ssp \DD_{1^{n+1}}\ses 0\ess .
$$
It thus follows that if we set $\ess D_x=\sum_{i=1}^{n}\ssp \del_{x_i}$, $\ess D_y=\sum_{i=1}^{n}\ssp \del_{y_i}$
then
$$
\del_{x_{n+1}}\ssp \DD_{1^{n+1}}\ses -D_x\ssp \DD_{1^{n+1}}
\ess\scs\ess\ess
D_y\ssp \DD_{1^{n+1}}\ses 0\ess .
\eqno 3.7
$$
Thus the fact that $\CB_{n+1}^x$, as given by 3.1, is a basis for $\BM_{1^{n+1}}$ yields that
the  collection
$$
\bigcup_{k=0}^{n}\ssp 
\bigl\{\ssp \del_{x_1}^{\eee_1}\del_{x_2}^{\eee_2}\cdots \del_{x_n}^{\eee_n}D_x^k\ssp \DD_{1^{n+1}}
\ssp :\ssp 0\leq \eee_i\leq i-1\ssp \bigr\}
$$
is also a basis. Since  
$$
\del_{x_n}^{\eee_n}D_x^k\ssp \DD_{1^{n+1}}\ssp \big|_{x_{n+1}=0 }\ses \DD_{1^{n+1}/k,0}\ess ,
$$
Proposition 1.1 yields that the collection
$$
\CB_{1^{n+1}/00}\ses 
\bigcup_{k=0}^{n}\ssp 
\bigl\{\ssp \del_{x_1}^{\eee_1}\del_{x_2}^{\eee_2}\cdots \del_{x_n}^{\eee_n}
\DD_{1^{n+1}/k,0}
\ssp :\ssp 0\leq \eee_i\leq i-1\ssp \bigr\}\ess .
$$
is also a basis for $\BM_{1^{n+1}}$.
Taking account of 3.1 we may also write this in the form
$$
\CB_{1^{n+1}/00}\ses 
\bigcup_{k=0}^{n}\ssp \CB_n(\del)\ssp \DD_{1^{n+1}/k,0}\ess .
$$
From this we immediately derive the direct sum decomposition
$$
\BM_{1^{n+1}/00}\ses \bigoplus_{k=0}^n\ess  \CL\bigl[\ssp \CB_n(\del)\ssp \DD_{1^{n+1}/k,0}\ssp \bigr]\ess .
\eqno 3.8
$$
Moreover, since $\ssp \BM_{1^{n+1}/i,0}=D_x^i\ssp \BM_{1^{n+1}/00}$ and $D_x^i\ssp\DD_{1^{n+1}/k,0}=0$
for $i+k>n$, by applying $D_x^i$ to both sides of 3.8 we also get that
$$
\BM_{1^{n+1}/i,0}\ses \bigoplus_{k=i}^n\ess  \CL\bigl[\ssp \CB_n(\del)\ssp \DD_{1^{n+1}/k,0}\ssp \bigr]\ess .
\eqno 3.9
$$
In particular we deduce that 
$$
\dim \BM_{1^{n+1}/i,0}\ses (n+1-i)\times n!\ess .
$$ 
This given, since each of the summands in 3.3 has dimension $n!$ and there are $n+1-i$ 
of them, to show 3.3 we need only verify that they are independent. 
To this end, assume  that for some elements $a_i,a_{i+1},\ldots ,a_{n+1}\in \BM_{1^n}$
we have
$$
a_i(\del)\DD_{1^{n+1}/i,0}\sps a_{i+1}(\del)\DD_{1^{n+1}/i+1 ,  0}\sps
\cdots \sps  a_{n+1}(\del)\DD_{1^{n+1}/n+1 ,  0}
\ses 0\ess .
\eqno 3.10
$$
Note first that for $i=n+1$, this equation reduces to
$$
a_{n+1}(\del)\DD_{1^{n}}\ses 0
$$
and since by choice $a_{n+1}\in \CL_\del[\DD_{1^{n}}]$ this forces $a_{n+1}=0\ess $.
So to show that $a_i,\ldots ,a_{n+1}=0$ we can proceed by descent induction on $i$.
That is, we can assume that 3.10 for $i+1$ forces $a_{i+1},\ldots ,a_{n+1}=0$.
This given, note that applying $D_x^{n+1-i}$ to 3.10 reduces it to
$$
a_i(\del)\ssp \DD_{1^{n+1}/n+1,0}\ses 0
$$
or, equivalently,
$$
a_i(\del)\ssp \DD_{1^n}\ses 0\ess .
$$
But this, as we have seen, forces $a_i=0$, yielding  that we must have
$$
 a_{i+1}(\del)\DD_{1^{n+1}/1+1,0}\sps
\cdots \sps  a_{n+1}(\del)\DD_{1^{n+1}/n+1,0}
\ses 0\ess ,
$$
and the induction hypothesis yields $a_{i+1},\ldots ,a_{n+1}=0\ssp $,
completing the induction and proving 3.3. It goes without saying that
3.4 may be proved in exactly the same way.

To show 3.5 a) we simply note that  
$$
D_x\ssp \BM_{1^n}(\del)\DD_{1^{n+1}/k,0}\ses 
\cases
{ \BM_{1^n}(\del)\DD_{1^{n+1}/k+1,0} & for $k\leq n$\ess , \cr\cr
0 & otherwise\ess .\cr
}
$$
Thus from 3.3 it follows that
$$
\BK_{i,0}^x\ses \BM_{1^n}(\del)\DD_{1^n}  =\ess  \BM_{1^n}
$$
and since in this case  $\BK_{i,1}^x=\{0\}$ we  must have $\BK_{i,0}^x=\BA_{i,0}^x$.

On the other hand, 3.5  b) follows from the fact that in this case for all $i$ we have
$$
\BK_{i,0}^y\ses \BM_{1^{n+1}/i,0} 
$$
Thus, using 3.3 again we get
$$
\BA_{i,0}^y\ses \BK_{i,0}^y/\BK_{i+1,0}^y\ses \BM_{1^{n+1}/i,0}/\BM_{1^{n+1}/i+1,0}
\ses \BM_{1^n}(\del)\ssp \DD_{1^{n+1}/i,0}\ess .
$$
This completes our proof since 3.6 a) and b) can be derived from 3.4 in precisely the same way.
\sa 

\heading{\bol Remark 3.1}

We should note that the direct sum expansions in 3.3 and 3.4 bring to evidence that the modules
$\BM_{1^{n+1}/i,0}$ and $\BM_{1^{n+1}/0,i}$ afford exactly $n+1-i$ copies
of the left regular representation of $S_n$ in complete agreement with our
Conjecture I.2.
\sa
 
Before we treat the general hook case $\mu=(n+1-k,1^k)$, it will be good 
to start by working with $\mu=(5,1,1,1)$. In this case we set
$$
\eqalign{\BM^{11}&\ses \BM_{511}\cap\BM_{4111}\ess  \scs \cr
\BM^{10}&\ses \BM_{511}\cap \bigl( \BM_{511}\cap\BM_{4111}\bigr)^\perp \ssp ,\cr
\BM^{01}&\ses \BM_{4111}\cap \bigl( \BM_{511}\cap\BM_{4111}\bigr)^\perp  . 
\cr
}
\eqno 3.11
$$
and obtain the decompositions:

$$
\eqalign{
\BM_{511}&\ses \BM^{11}\ssp \oplus \ssp \BM^{10}
\cr
\BM_{4111}&\ses \BM^{11}\ssp \oplus \ssp \BM^{01}
\cr
}
$$
Using the convention that placing a module $\ssp \BM\ssp$  or a basis $\ssp\CB\ssp$ in cell $(i,j)$ of the diagram 
of $\mu$ represents applying $\BM(\del)$ or  $\CB(\del)$ to $\DD_{\mu/ij}$, formula  2.6
asserts that we must have
$$
\BM_{5111/00}= 
\tableau
{
   \BM^{11}    \cr
 \BM^{11}     \cr
   \BM^{11}   \cr
  \BM^{11}& \BM^{11} & \BM^{11} & \BM^{11} & \BM^{11}\cr
}
 \oplus \;
\tableau
{
   \BM^{10}     \cr
  \BM^{10}      \cr
  \BM^{10}     \cr
  \BM^{10}& \emptyset  & \emptyset & \emptyset & \emptyset \cr
}
  \oplus  \;
\tableau{
   \emptyset    \cr
  \emptyset      \cr
  \emptyset     \cr
   \BM^{01} & \BM^{01} & \BM^{01}  & \BM^{01} & \emptyset \cr
}$$
Taking account of 3.11 and setting $\bf A=\BM_{511}$ and $\bf B=\BM_{4111}\ssp$,
this identity may be compressed to
$$
\BM_{5111/00}= 
\vcenter{\offinterlineskip
\halign{\vrule height12pt depth6pt$\;#\;$\vrule
	&\strut$\;#\;$\vrule
	&\strut$\;#\;$\vrule
	&\strut$\;#\;$\vrule
	&\strut$\;#\;$\vrule\cr
	\multispan1\hrulefill\cr
  \BA \cr
	\multispan1\hrulefill\cr
 \BA \cr
	\multispan1\hrulefill\cr
  \BA \cr
	\multispan5\hrulefill\cr
 \BA+\BB &\BB &\BB &\BB &\BA\cap \BB \cr
	\multispan5\hrulefill\cr
}
}
$$
Letting $\ess \CB_a\scs \CB_b\scs \CB_{a+b}\scs \CB_{a\,\cap\, b} $
denote bases for $\BA\scs \BB \scs \BA+\BB \scs \BA\cap \BB\ess$
respectively, and writing $\DD_{ij}$ for $\DD_{\mu/ij}$,
 this formula  asserts that a basis for the module
$\BM_{5111/00}$ is given by the
collection
\def \ssps {\ssp +\ssp}
$$
\eqalign{
\CB_{5111/00}&= \ssp 
\CB_a(\del)\DD_{30} \;\cup\;
\CB_a(\del)\DD_{20} \;\cup\;
\CB_a(\del)\DD_{10}
\cr 
&\ess\ess\ess\ess\ess\ess\ess\ess
 \;\cup\; 
\CB_{a+b}(\del)\DD_{00} \;\cup\;
\CB_b(\del)\DD_{01} \;\cup\;
\CB_b(\del)\DD_{02} \;\cup\;
\CB_b(\del)\DD_{03} \;\cup\; \CB_{a\,\cap\, b}(\del)\DD_{04}\ssp .
\cr} 
\eqno 3.12
$$
In particular we get that this set has cardinality  
$$
4\times \dim \BA \sps 4\times \dim \BB\ess .
$$  
Thus, assuming that $\dim \BA=\dim \BB= 7!\ssp $, we deduce that
$\CB_{5111/00}$ has precisely $8!$ elements. Since it was shown in 
[7] that for $\mu\part n$ we have  $\dim \BM_\mu \leq n!\ssp $,
from  Proposition I.5 we derive that $\dim \BM_{5111/00}\leq 8!$.
Thus, to show that $\CB_{5111/00}$ is a basis we need only verify that
it is an independent set. To this end let (for $i=1,\ldots,3$)
$$
a_i\in \CL[\CB_a]\ess \scs \ess b_i\in \CL[\CB_b]\ess  \ess   
 \scs \ess u\in \CL[\CB_{a+b}]
\scs \ess v\in \CL[\CB_{a\,\cap\, b}]
$$
and suppose if possible that 
$$
\eqalign{
&
a_3(\del)\DD_{30}\ssps
a_2(\del)\DD_{20}\ssps
a_1(\del)\DD_{10}\ssps
\cr 
&\ess\ess\ess\ess\ess\ess\ess\ess
\ssps 
u(\del)\DD_{00}\ssps
b_1(\del)\DD_{01}\ssps
b_2(\del)\DD_{02}\ssps
b_3(\del)\DD_{03}\ssps v(\del)\DD_{04}\ses 0\ess.
\cr} 
\eqno 3.13
$$
\def\sdes{\ess \dot = \ess}
Let the symbol ``$\sdes$'' represent
		  equality up to a constant factor.
Proposition I.2 gives
$$
\eqalign{
& D_x^3\DD_{00}\sdes D_x^2\DD_{10}\sdes D_x \DD_{20}\sdes \DD_{30} \cr
& D_y ^4\DD_{00}\sdes D_y^3\DD_{01}\sdes D_y^2 \DD_{02}\sdes D_y\DD_{03}\sdes \DD_{04} \cr
& D_x \DD_{0,i}=0\ess\ess \&\ess \ess D_y \DD_{i,0}=0\ess\ess \ess {\rm for}\ess\ess i>0\cr
& D_x \DD_{30}=0\ess\ess \&\ess \ess D_y \DD_{04}=0\ess\ess , \cr
}
\eqno 3.14
$$
where the last of these equations results from the fact   that 
$$
\DD_{30}\ses\DD_{511}
\ess\ess {\rm and}\ess\ess \DD_{04}\ses \DD_{4111}\ess .
$$
Thus applying $D_x^3$ to 3.13 reduces it  to
$$
u(\del)\ssp \DD_{511}\ses 0\ess .
\eqno 3.15
$$
Similarly applying $D_y^4$ to 3.13 gives
$$
u(\del)\ssp \DD_{4111}\ses  0 \ess .
\eqno 3.16
$$
Since by assumption $u\in \CL_\del[\DD_{511}]+\CL_\del[\DD_{4111}]$, equations
3.15 and 3.16 force $u$ to be orthogonal to itself and therefore identically zero. 
So 3.3 becomes
$$
\eqalign{
&
a_3(\del)\DD_{30}\ssps
a_2(\del)\DD_{20}\ssps
a_1(\del)\DD_{10}
\cr 
&\ess\ess\ess\ess\ess\ess\ess\ess
\ssps 
b_1(\del)\DD_{01}\ssps
b_2(\del)\DD_{02}\ssps
b_3(\del)\DD_{03}\ssps v(\del)\DD_{04}\ses 0\ess.
\cr} 
\eqno 3.17
$$
Now, the relations in 3.14  yield that applying $D_x^2$ to 3.17 reduces it to
$$
a_1(\del)\DD_{511}\ses 0
$$ 
and since $\ess a_1\in \CL_\del[\DD_{511}]\ssp $, we derive that $a_1=0$ as well,
reducing 3.17 to
$$
\eqalign{
&
a_3(\del)\DD_{30}\ssps
a_2(\del)\DD_{20}
\cr 
&\ess\ess\ess\ess\ess\ess\ess\ess
\ssps 
b_1(\del)\DD_{01}\ssps
b_2(\del)\DD_{02}\ssps
b_3(\del)\DD_{03}\ssps v(\del)\DD_{04}\ses 0\ess.
\cr} 
\eqno 3.18
$$  
Now an application of $D_x$ yields
$$
a_2(\del)\DD_{511}\ses 0
$$
and $\ess a_2\in \CL_\del[\DD_{511}]\ssp $ yields again $a_2=0\ssp $,
reducing 3.18 to 
$$
\eqalign{
&
a_3(\del)\DD_{30}\ssps
\cr 
&\ess\ess\ess\ess\ess\ess\ess\ess
\ssps 
b_1(\del)\DD_{01}\ssps
b_2(\del)\DD_{02}\ssps
b_3(\del)\DD_{03}\ssps v(\del)\DD_{04}\ses 0\ess.
\cr} 
\eqno 3.19
$$ 
Since $D_y\DD_{30}=0$, we can now apply $D_y^3, D_y^2,  D_y^1\ssp $ in succession to 3.9
and, by a similar process, successively derive that
$$
b_1\scs b_2\scs b_3\ses   0\ess ,
$$ 
reducing 3.19 to
$$
a_3(\del)\DD_{511}\sps v(\del )\DD_{4111}\ses 0\ess .
\eqno 3.20
$$
As we let $a_3$ vary in $\CL_\del[\CB_a]$ without restriction, the term $a_3(\del)\DD_{511}$
will necessarily describe all of $\BA$. On the other hand, as $v$ varies in $\CL_\del[\CB_{a\,\cap\, b}]$
the term $v(\del )\DD_{4111}$ will describe $\flip_{4111} \bigl(\BA\cap \BB\bigr)$.
\sas

This given, to conclude from  for 3.10 that $a_3$ and $v$ must vanish we need to know that
$\BA$ and $\flip_{4111} \bigl(\BA\cap \BB\bigr)$ have no element in common other than $0$.
It is at this point that the SF hypothesis plays a role. In fact, in the particular case
of a 2-corner partition $\mu$, with two predecessors $\aaa_1,\aaa_2$, condition (iii) of
SF  asserts (see [1]  I.29) that
$$
1)\ess\ess\flip_{\aaa_1} \BM^{11}\ses \BM^{10}
\ess\ess\ess\ess {\rm and}\ess\ess\ess\ess
2)\ess\ess\flip_{\aaa_2} \BM^{11}\ses \BM^{01}\ess.
\eqno 3.21
$$
This of course guarantees that the two terms in 3.20 must separately vanish, completing
the proof that the collection $\CB_{5111/00}$ defined in 3.13 is a basis for $\BM_{5111/00}$.
\sa

\heading{\bol Remark 3.2}

We should note that although they can be verified by computer in several special 
cases, the identities in 3.21 may be too strong to be true in general. A weaker form,
which does not affect the final conclusion, is obtained by changing the definitions 
of $\BM^{10}$ and $\BM^{01}$ by dropping the condition that they be orthogonal 
complements  of $\BM^{11}$ in  $\BM_{511}$  and $\BM_{4111}$ respectively
and just require that they be simply ``complements'' constructed so that 
the relations in 3.21 are satisfied.  Another way to get around requiring
3.21 is to observe that the desired implication
$$
a_3(\del)\DD_{511}\sps v(\del )\DD_{4111}\ses 0\ess \ess\ess \Longrightarrow
\ess\ess\ess 
a_3(\del)\DD_{511}\ess \&\ess  v(\del )\DD_{4111}\ses 0
$$
immediately follows, if the collection $\CB_{a\,\cap\, b}$ is replaced by any basis
of $\ess \flip^{-1}_{4111}\ssp \BM^{01}$. This choice guarantees  that
$\CB_{5111/00}$ is an independent set. However, to conclude that $\CB_{5111/00}$ is a basis we
need  
$$
\dim \BM^{11}\ses \dim \BM^{01}\ess ,
$$
or equivalently  
$$
\dim\BM^{11}\ses  {\dim\BM_{\aaa_1}\ssp \over 2}\ess.
\eqno 3.22
$$
Unfortunately, this equality, which has come to be referred to as the $n!/2$ conjecture,
has to this date remained unproved in full generality (even in the ``hook'' case). As a result,
this modified construction of  $\CB_{5111/00}$ only generalizes to a proof
that 
$$
\dim \BM^{11}\geq \dim \BM^{01}\ess .
$$

These observations are essentially all contained in the SF paper [1]. What is new here
is that the introduction of the ``atoms'' ${\bf A}_{ij}^x$ and ${\bf A}_{ij}^y$
leads to a very elegant construction of a basis of $\BM_\mu$ when $\mu$ is
a hook without any need of unproved auxiliary conjectures. We shall illustrate
it here again in the case $\mu=(5,1,1,1,1)$. 
\sas

Since the construction is inductive on the size of $\mu$
we shall again assume that both $\BM_{511}$ and $\BM_{4111}$ have dimension $7!$ and
that we have chosen $\CB_{511}$ and $\CB_{4111}$ as their respective bases.
This given, we may represent our alternate construction of a basis $\tilde \CB_{5111/00}$
for $\BM_{5111/00}$ 
by the diagram
$$
\tilde \CB_{5111/00}= 
\vcenter{\offinterlineskip
\halign{\vrule height12pt depth6pt$\;#\;$\vrule
	&\strut$\;#\;$\vrule
	&\strut$\;#\;$\vrule
	&\strut$\;#\;$\vrule
	&\strut$\;#\;$\vrule\cr
	\multispan1\hrulefill\cr
  \emptyset     \cr
	\multispan1\hrulefill\cr
 \CB_{511}   \cr
	\multispan1\hrulefill\cr
  \CB_{511}  \cr
	\multispan5\hrulefill\cr
 \CB_{511}\cup \cal X &\CB_{4111} &\CB_{4111} &\CB_{4111} &\CB_{4111} \cr
	\multispan5\hrulefill\cr
}}
\eqno 3.23
$$
where $\cal X$ is a suitable collection of monomials. Before we 
exhibit our choice of $\cal X$, it will be instructive to see
that 3.23 gives a basis for $\BM_{5111/00}$ as soon as $\cal X$ satisfies
the following three conditions:
\sas

\itemitem {(i)} $\CX(\del)\DD_{00}$ is an independent set of cardinality $7!$,

\itemitem {(ii)} $D_x \ssp m(\del)\DD_{00}=0\ess\ess \forall \ess \ess m\in \CX$,
\hfill $3.24$

\itemitem {(iii)} For any $0\neq\xi\in \CL[\CX]$ the element $\xi(\del)\DD_{00}$ is never in

\itemitem{} \hfill$
\CL\bigl[\ssp \CB_{4111}(\del)\DD_{01}
		\;\cup\;
		\CB_{4111}(\del)\DD_{02}
		\;\cup\;
		\CB_{4111}(\del)\DD_{03}
		\;\cup\;
		\CB_{4111}(\del)\DD_{04}\ssp \bigr]\ess.
$\hfill\break

\noindent
In fact, suppose that for some $\ess\ess a_0,a_1,a_2\in \CL[\CB_{511}]\ssp ,\ess\ess$ 
$b_1,b_2,b_3,b_4 \in \CL[\CB_{4111}]\ess $
and $\ess \xi\in \CL[\CX]$ we have
$$
\eqalign{
& a_2(\del)\DD_{20}\sps a_1(\del)\DD_{10}\sps a_0(\del)\DD_{00}
\cr
&
\ess\ess\ess\ess\ess\ess\ess\ess\ess\ess\ess\ess\sps \xi(\del)\DD_{00}\sps
 b_1(\del)\DD_{01}\sps b_2(\del)\DD_{02}\sps b_3(\del)\DD_{03}\sps b_4(\del)\DD_{04}\ses 0\ess . 
\cr
}
\eqno 3.25
$$
To show that this forces $a_1,a_2,a_3,\xi,b_1,b_2,b_3,b_4= 0$
we apply $D_x$ to both sides and,
using
the relations in 3.14 and condition (ii) of 3.24,
immediately derive that 
$$
a_2(\del)\DD_{30}\sps a_1(\del)\DD_{20}\sps a_0(\del)\DD_{10}\ses 0\ess .
\eqno 3.26
$$
This given, an application  of $D_x^2$ reduces this to
$$
a_0(\del)\DD_{30}\ses 0
	\ess,
$$
which as we have seen forces $a_0=0$ and 3.26 becomes
$$
a_2(\del)\DD_{30}\sps a_1(\del)\DD_{20}\ses 0\ess .
\eqno 3.27
$$
Applying $D_x\,$, we now get
$$
a_1(\del)\DD_{30}\ses 0\ess ,
$$
which forces $a_1=0\ssp $, reducing 3.17 to
$$
a_2(\del)\DD_{30}\ses 0\ess ,
$$
and this in turn yields
$$
a_2\ses 0\ess .
$$
So 3.15 becomes
$$
\xi(\del)\DD_{00}\sps b_1(\del)\DD_{01}\sps b_2(\del)\DD_{02}\sps b_3(\del)\DD_{03}\sps b_4(\del)\DD_{04}\ses 0\ess .
$$
But then condition (iii) of 3.14 assures that we must separately have
$$
\eqalign{
 \xi(\del)\DD_{00}&\ses 0
\cr
b_1(\del)\DD_{01}\sps b_2(\del)\DD_{02}\sps b_3(\del)\DD_{03}\sps b_4(\del)\DD_{04}&\ses 0
\cr
}
$$
Now, the first equation (using 3.24(i))   yields $\xi=0$, while the second 
yields $b_1,b_2,b_3=0$ 
by successive applications of $D_y^3,D_y^2,D_y $,  
as we have seen before. We are finally left with
$$
 b_4(\del)\DD_{04} \ses 0
	\ess ,
$$
which forces $b_4=0$ and completes the proof of independence of $\tilde \CB_{5111/00}$.
Since by virtue of (i) in 3.14 the cardinality of $\tilde \CB_{5111/00}$ evaluates to $8!\ssp $,
we must conclude that $\tilde \CB_{5111/00}$ must also be a basis.
\sas

It develops that a collection of monomials that satisfies all of the condition in
3.14 is obtained by setting
$$
\CX= \hskip -.2truein \bigcup_{\multi{1\leq i_1<i_2<i_3\leq 7 \cr
1\leq j_1<j_2<j_3<j_4\leq 7\cr  \{i_1,i_2,i_3,j_1,j_2,j_3,j_4\}=\{1,2,3,4,5,6,7\}}}
\hskip -.27truein\bigl\{
\ssp x_{i_1}^{1+\eee_1}
\ssp  x_{i_2}^{1+\eee_2}
\ssp  x_{i_3}^{1+\eee_1}
\ssp  y_{j_1}^{\eta_1}
\ssp  y_{j_2}^{\eta_2}
\ssp  y_{j_3}^{\eta_3}
\ssp  y_{j_4}^{\eta_4}
\ssp :\ssp 0\leq \eee_i\leq i-1\ssp; 
\ssp 0\leq \eta_j\leq j-1\ssp\ssp \bigr\}
\eqno 3.28
$$ 
To see this, note first that since
$$
\DD_{5111/00}\ses \det \pmatrix{
y_1 &y_2 &y_3 &y_4 &y_5 &y_6 &y_7 \cr
y_1^2 &y_2^2 &y_3^2 &y_4^2 &y_5v &y_6^2 &y_7^2 \cr
y_1^3 &y_2^3 &y_3^3 &y_4^3 &y_5^3 &y_6^3 &y_7^3 \cr
y_1^4 &y_2^4 &y_3^4 &y_4^4 &y_5^4 &y_6^4 &y_7^4 \cr
x_1 &x_2 &x_3 &x_4 &x_5 &x_6 &x_7 \cr
x_1^2 &x_2^2 &x_3^2 &x_4^2 &x_5^2 &x_6^2 &x_7^2 \cr
x_1^3 &x_2^3 &x_3^3 &x_4^3 &x_5^3 &x_6^3 &x_7^3 \cr
}
$$
we have
$$
\eqalign{
\del_{x_1}\del_{x_2}\del_{x_3}\DD_{5111/00}  
&\ses
6\times \det
\pmatrix{
1 &1 &1   \cr
x_1  &x_2  &x_3   \cr
x_1^2 &x_2^2 &x_3^2  \cr
} \ssp\times \ssp
\det
\pmatrix{
 y_4 &y_5 &y_6 &y_7 \cr
 y_4^2 &y_5  &y_6^2 &y_7^2 \cr
 y_4^3 &y_5^3 &y_6^3 &y_7^3 \cr
 y_4^4 &y_5^4 &y_6^4 &y_7^4 \cr
} 
\cr
&\ses
6\times y_4 \ssp  y_5 \ssp y_6\ssp y_7  \ssp\times \ssp
\det\pmatrix{
1 &1 &1   \cr
x_1  &x_2  &x_3   \cr
x_1^2 &x_2^2 &x_3^2  \cr
}
 \ssp\times \ssp \det
\pmatrix{
 1 &1 &1 &1  \cr
y_4 &y_5 &y_6 &y_7 \cr
 y_4^2 &y_5  &y_6^2 &y_7^2 \cr
 y_4^3 &y_5^3 &y_6^3 &y_7^3 \cr
} \ess . \cr 
} 
$$ 
\sas

Thus using the notation we introduced at the beginning of the section (see 3.1 and 3.2) we may
write
$$
\eqalign{
\bigl\{\ssp\del_{x_1}^{1+\eee_1}\del_{x_2}^{1+\eee_2}\del_{x_3}^{1+\eee_3}
\del_{y_1}^{\eta_4}\del_{y_5}^{\eta_2}\del_{y_6}^{\eta_3}\del_{y_7}^{\eta_3}
\DD_{5111/00}\ssp :
&\ssp 0\leq  \eee_i\leq i-1\ssp ;\ssp 0\leq  \eta_j\leq j-1\ssp\bigr\}  \cr
&\ess\ess\ess\ess\ess \sdes y_4\ssp y_5\ssp y_6\ssp y_7\ssp
\times \ssp\CB_3(x_1,x_2,x_3)\ssp\times\ssp \CB_3(y_4,y_5,y_6,y_7) 
\cr  
}  
$$
So we see that 
the module $\CL[\CX(\del)\DD_{5111/00}]$ may be expressed as the direct sum
$$
\CL\bigl[\CX(\del)\DD_{5111/00}\bigr]\ses 
\hskip -.3truein 
\bigoplus_{\multi{1\leq i_1<i_2<i_3\leq 7 \cr
1\leq j_1<j_2<j_3<j_4\leq 7\cr  \{i_1,i_2,i_3;j_1,j_2,j_3,j_4\}=\{1,2,3,4,5,6,7\}}}
\hskip -.3truein 
 y_{j_1} \ssp y_{j_2}\ssp y_{j_3}\ssp y_{j_4} \times  
\BM_{1^3}[x_{i_1},x_{i_2},x_{i_3}] \times   \BM_{4}[y_{j_1},y_{j_2},y_{j_3},y_{j_4}]\ssp ,
\eqno 3.29
$$
where the symbols $\BM_{1^3}[x_{i_1},x_{i_2},x_{i_3}]$ and $\BM_{4}[y_{j_1},y_{j_2},y_{j_3},y_{j_4}]$ 
denote $\BM_{1^3}$ and $\BM_{4}$ with $x_s$ replaced by $x_{i_s}$ and $y_r$
replaced by $y_{j_r}$. Now we immediately derive from this that
$$
\dim \CL\bigl[\CX(\del)\DD_{5111/00}\bigr]\ses \Bigl(\ssp {7 \atop 3}\ssp \Bigr)\ssp 3!\ssp 4!\ses 7!
$$
yielding 3.14 (i). Now 3.14 (ii) is immediate since for any choice of $x_{i_1}x_{i_2}x_{i_3}$
we have
$$
D_x\ssp \DD_3(x_{i_1},x_{i_2},x_{i_3})\ses 0\ess .
$$
Finally we note that every one of the determinants $\DD_{i,0}$ is a sum of monomials
only containing three different $\ssp y_i\ssp '\ssp s$ and thus none  of their derivatives
can contain monomials with four different $\ssp y_i\ssp '\ssp s$. Since each element
of $\CL\bigl[\CX(\del)\DD_{5111/00}\bigr]$ has  $y_{j_1} \ssp y_{j_2}\ssp y_{j_3}\ssp y_{j_4}$
as a factor we see that 3.24 (iii) must necessarily hold true precisely as required.
\sas

Now the fact that 3.13, with  $\CX$ given by 3.28, gives a basis for $\BM_{5111/00}$
yields that $\BM_{5111/00}$ has a direct sum decomposition
$$
\eqalign{
\BM_{5111/00} \ses& \BM_{511}(\del)\DD_{5111/20}
\ess \oplus\ess
\BM_{511}(\del)\DD_{5111/10} 
\cr
&
\hskip .5in
 \oplus\ess
\CL[\CX(\del)\DD_{00}]\ess \oplus\ess\BM_{511}(\del)\DD_{00} 
\cr
&
\hskip 1in
\oplus\ess
\BM_{4111}(\del)\DD_{5111/01}
\ess \oplus\ess
\BM_{4111}(\del)\DD_{5111/02}
\cr
&
\hskip 1.5in
\oplus\ess
\BM_{4111}(\del)\DD_{5111/03}
\ess \oplus\ess
\BM_{4111}(\del)\DD_{5111/04}\ess . 
\cr}
\eqno 3.30
$$
Since $\ess\BM_{5111/i,0}=D_x^i\ssp \BM_{5111/00}\ess $ for $i=1,2,3$ and $D_x$
kills $\CL[\CX(\del)\DD_{00}]$ as well as each $\DD_{5111/0,j}$,
we immediately derive, by applying $D_x,D_x^2, D_x^3$ to both sides of 3.30, that
$$
\eqalign{
\BM_{5111/1 0}&\ses
 \BM_{511}(\del)\DD_{5111/10}  
\ess \oplus\ess
\BM_{511}(\del)\DD_{5111/20}
\ess \oplus\ess
\BM_{511} 
\cr
\BM_{5111/2 0}
&\ses
\BM_{511}(\del)\DD_{5111/20}
\ess \oplus\ess
\BM_{511}
\cr
\BM_{5111/3 0}
&\ses
 \BM_{511}\ess ,
\cr}  
\eqno 3.31
$$
where we have used the fact that $\ssp \DD_{5111/30}=\DD_{511}\ssp $.
Similarly by inverting the roles played by the $x$ and $y$ variables we derive 
the direct sum decompositions
$$
\eqalign{
\BM_{5111/0 1}&\ses
\BM_{4111}(\del)\DD_{5111/01}  
\ess \oplus\ess
 \BM_{4111}(\del)\DD_{5111/02}  
\ess \oplus\ess
\BM_{4111}(\del)\DD_{5111/03}
\ess \oplus\ess
\BM_{4111} 
\cr
\BM_{5111/0 2}&\ses
 \BM_{4111}(\del)\DD_{5111/02}  
\ess \oplus\ess
\BM_{4111}(\del)\DD_{5111/03}
\ess \oplus\ess
\BM_{4111} 
\cr
\BM_{5111/0 3}
&\ses
\BM_{4111}(\del)\DD_{5111/03}
\ess \oplus\ess
\BM_{4111} 
\cr
\BM_{5111/0 4}
&\ses
 \BM_{4111}\ess .
\cr}  
\eqno 3.32
$$
Note that 3.30 gives as also that
$$
\eqalign{
\BK_{00}^x\ses &
\CL[\CX(\del)\DD_{00}]\ess   
 \oplus\ess
\BM_{4111}(\del)\DD_{5111/01}
\ess \oplus\ess
\BM_{4111}(\del)\DD_{5111/02}
\cr
&
 \ess \ess \ess \ess \ess \ess  \ess\ess \ess \ess \ess \ess \ess \ess  
\ess\ess \ess \ess \ess \ess \ess  \ess \ess \ess \oplus\ess
\BM_{4111}(\del)\DD_{5111/03}
\ess \oplus\ess
\BM_{4111}(\del)\DD_{5111/04}\ess .
\cr
}
$$ 
In other words
$$
\BK_{00}^x\ses  
\CL[\CX(\del)\DD_{00}]\ess   
 \oplus\ess
\BM_{5111/01}\ess . 
$$
Since $D_x$  kills all of $\ssp \BM_{5111/01}\ssp$. This may be rewritten as
$$
\BK_{00}^x\ses  
\CL[\CX(\del)\DD_{00}]\ess   
 \oplus\ess
\BK_{01}^x \ess ,
$$
yielding that in this case we have
$$
\BA_{00}^x\ses \CL[\CX(\del)\DD_{00}]\ess .
\eqno 3.33
$$
By ``equality,'' we mean that $\CL[\CX(\del)\DD_{00}]$ is a
complement of $\BK^x_{01}$ within $\BK^x_{00}$, thus
forming a system of representatives of the quotient $\BA_{00}^x=\BK_{00}^x/\BK_{01}^x$.
Similarly we can derive from 3.31 and 3.32 that
$$
\BA_{10}^x=\BA_{20}^x=\BA_{30}^x= \BM_{511}
\eqno 3.34
$$ 
and
$$
\BA_{0,i}^x=\BM_{4111}(\del)\DD_{5111/0,i}
\bigsp (\hbox{\ssp for $i=1,2,3,4\ssp $})\ess .
\eqno 3.35
$$

We should point out that analogous results concerning the atoms $\BA_{ij}^y$ 
can be obtained if construct the basis $\tilde \CB_{5111/00}$ according
to the ``transposed'' diagram
$$
\tilde \CB_{5111/00}= 
\vcenter{\offinterlineskip
\halign{\vrule height12pt depth6pt$\;#\;$\vrule
	&\strut$\;#\;$\vrule
	&\strut$\;#\;$\vrule
	&\strut$\;#\;$\vrule
	&\strut$\;#\;$\vrule\cr
	\multispan1\hrulefill\cr
  \CB_{511}   \cr
	\multispan1\hrulefill\cr
 \CB_{511}   \cr
	\multispan1\hrulefill\cr
  \CB_{511}   \cr
	\multispan5\hrulefill\cr
 \CB_{4111}\cup \CY &\CB_{4111} &\CB_{4111} &\CB_{4111} & \;\;\emptyset\;\;  \cr
	\multispan5\hrulefill\cr
}}
\eqno 3.36
$$
with 
$$
\CY= \hskip -.3truein \bigcup_{\multi{1\leq i_1<i_2<i_3\leq 7 \cr
1\leq j_1<j_2<j_3<j_4\leq 7\cr  \{i_1,i_2,i_3;j_1,j_2,j_3,j_4\}=\{1,2,3,4,5,6,7\}}}
\hskip -.3truein\bigl\{
\ssp x_{i_1}^{ \eee_1}
\ssp  x_{i_2}^{ \eee_2}
\ssp  x_{i_3}^{ \eee_3}
\ssp  y_{j_1}^{1+\eta_1}
\ssp  y_{j_2}^{1+\eta_2}
\ssp  y_{j_3}^{1+\eta_3}
\ssp  y_{j_4}^{1+\eta_4}
\ssp :\ssp 0\leq \eee_i\leq i-1\ssp; 
\ssp 0\leq \eta_j\leq j-1\ssp\ssp \bigr\}\ess .
$$ 
It should also be clear that the argument we have illustrated in the case $\mu=5111$ 
can be carried out for all hook partitions. In fact, in this case all our conjectures
can be proved in full including the $C=\TH$  conjecture and the four term recursion.
\sas

For a given  subset $S=\{  i_1< i_2<\cdots <i_k \}$ let $|S|=k$ and set
$$
X(S)=\{x_{i_1},x_{i_2},\ldots ,x_{i_k}\}\ess  \scs\ess\ess
Y(S)=\{y_{i_1},y_{i_2},\ldots ,y_{i_k}\}
\ess .
$$
 Moreover, if $\BM$ is
a space of polynomials in the variables $\ssp x_1,x_2,\ldots , x_k\ssp $, 
let $\BM[X(S)]$  denote the space obtained by replacing $x_s$
by $x_{i_s}$ in all elements of $\BM$. Let $\BM[Y(S)]$
be analogously defined with the $y's$ replacing the $x's$. 
Recall that according to the definitions made in the
introduction,  $\BM_{1^n}$ and $\BM_n$ 
denote the linear spans of derivatives of the Vandermonde determinants
in $x_1,x_2,\cdots ,x_n$ and $y_1,y_2,\ldots ,y_n$ respectively.
With  these conventions, our general result for hooks
may be stated as follows. 
\sas

\heading{\bol Theorem 3.2}

{\ita For $\mu=(n+1-k,1^k)$, set $\aaa=(n+1-k,1^{k-1})$ and $\bbb=(n-k,1^{k})$.
Let
$$
\BX\ses\hskip -.3truein \bigoplus_{\multi{|S|=k\cr|T|=n-k \cr
S+T=\{1,2,\ldots ,n\}\cr}}
\hskip -.1truein 
\Bigl(\ssp \prod_{j\in T} y_j\ssp \Bigr)\times\BM_{1^k}\bigl[X(S)\bigr]\times \BM_{n-k}\bigl[Y(T)\bigr]
\eqno 3.37 
$$
and
$$
\BY\ses \hskip -.3truein \bigoplus_{\multi{|S|=k\cr|T|=n-k \cr
S+T=\{1,2,\ldots ,n\}\cr}}
\hskip -.1truein 
\Bigl(\ssp \prod_{i\in S} x_i\ssp \Bigr)\times\BM_{1^k}\bigl[X(S)\bigr]\times \BM_{n-k}\bigl[Y(T)\bigr]\ess .
\eqno 3.38 
$$
This given, we have the following direct sum decompositions:
$$
\eqalign{
a)\ess\ess \BM_{\mu/00}&\ses \bigoplus_{i=0}^{k-1}\BM_{\aaa}(\del)\DD_{\mu/i,0}
\ess  \oplus \ess \BX \ess  \oplus \ess \bigoplus_{j=1}^{n-k}\BM_{\bbb}(\del)\DD_{\mu/0, j}
\cr
b)\ess\ess \BM_{\mu/00}&\ses \bigoplus_{i=1}^{k}\BM_{\aaa}(\del)\DD_{\mu/i,0}
\ess  \oplus \ess \BY \ess  \oplus \ess \bigoplus_{j=0}^{n-k-1}\BM_{\bbb}(\del)\DD_{\mu/0, j}
\cr}
\eqno 3.39
$$
$$
a)\ess\ess \BM_{\mu/i,0} \ses \bigoplus_{r=i}^{k}\BM_{\aaa}(\del)\DD_{\mu/r,0}
\ess\ess\scs\ess\ess\ess
b)\ess\ess \BM_{\mu/0,j} \ses
  \bigoplus_{s=j}^{n-k}\BM_{\bbb}(\del)\DD_{\mu/0,s}
\eqno 3.40$$
with
$$
\eqalign{
&a)\ess\ess\BA_{00}^x\ses \BX
 \ess\scs\ess\ess\ess
\BA_{i,0}^x\ses\BM_\aaa
\ess\ess\ess\scs\ess\ess\ess\ess\ess\ess\ess
\BA_{0,j}^y\ses \BM_\bbb(\del)\DD_{\mu/0,j}\ess ,
\cr 
&b)\ess\ess\BA_{00}^y\ses \BY
 \ess\scs\ess\ess\ess
\BA_{i,0}^y\ses\BM_\aaa(\del)\DD_{\mu/i,0}
\scs\ess\ess\ess
\BA_{0,j}^y\ses \BM_\bbb\ess . 
\cr 
}
\eqno 3.41
$$
Moreover, the Frobenius characteristics of these modules may be expressed in terms of the 
Macdonald polynomials as follows:
$$
\matrix{
a)\ess\ess \Fch  \BA_{00}^x\ses q^{n-k}\ssp \TH_{1^k}\ssp \TH_{n-k}\cr
\cr
b)\ssp\ess\ess \Fch \BA_{00}^y\ses\ssp  t^{ k}\ssp\ess \TH_{1^k}\ssp \TH_{n-k} \cr
}
\ess\ess\ess\ess\ess\ess\ess\ess\ess c)\ess\ess\ess \Fch \BM_{(n+1-k,1^k)}=\TH_{(n+1-k,1^k)}
\eqno 3.42
$$
}
\heading{\bol Proof}

Formulas 3.39 a) and b) may be obtained by generalizing the argument  that yielded 3.30.
Similarly 3.40 a) and b) can be easily established by the process that gives 
3.31 and 3.32. This given, since $D_x\BX=\{0\}$ and $D_x\DD_{\mu/0,j}=0$,
it follows from
3.39 a) and 3.40 a) that
$$
D_x\ssp \BM_{\mu/00}=  \bigoplus_{i=1}^{k }\BM_{\aaa}(\del)\DD_{\mu/i,0}
= \BM_{\mu/{10}}
\ess\scs\ess\ess
D_x\ssp \BM_{\mu/01}= 0\ess .
$$
Thus 
$$
\BK_{01}^x= \BM_{\mu/{10}}\ess\scs\ess\ess
\BK_{00}^x\ses \BX\ess \oplus\ess \bigoplus_{j=1}^{n-k}\ssp \BM_\bbb(\del) \DD_{\mu/0,j}
\ses \BX\ess \oplus\ess \BM_{\mu/{10}}\ess ,
$$
yielding
$$
\BK_{00}^x\ses \BX\ess \oplus\ess \BK_{01}^x\ess
$$
and 3.41 then follows from the definition I.18. Formula 3.39 b) is 
established in a similar manner. 
The remaining identities in 3.41
follow from 3.39 and the stated properties of $D_x$ and $D_y$. 

Thus it only remains to prove the Macdonality of the Frobenius characteristics
as stated in 3.42. To begin with we note that it is well known (see [2], [10]) that the linear
span of the derivatives of the Vandermonde determinant $\DD_n(x_1,x_2,\ldots ,x_n)$ 
yields a graded version of the left regular representation of $S_n$ with 
Frobenius characteristic given by the symmetric polynomial
$$
(1-t)(1-t^2)\cdots (1-t^n)\ssp h_n\bigl[{\textstyle {X\over 1-t}}\bigr]
\ses \sum_{\la\part n}\ssp S_\la[X]\ssp S_\la[1,t,t^2,\dots ](1-t)(1-t^2)\cdots (1-t^n)\ess .
$$
Now we have shown in [10] that
$$
\TH_{1^n}\ses (1-t)(1-t^2)\cdots (1-t^n)\ssp h_n\bigl[{\textstyle {X\over 1-t}}\bigr]
$$
and
$$
\TH_{n}\ses (1-q)(1-q^2)\cdots (1-q^n)\ssp h_n\bigl[{\textstyle {X\over 1-q}}\bigr]\ess .
$$
Thus formula 3.37 defines $\BX$ as the bigraded module obtained by inducing from
$S_k\times S_{n-k}$ to $S_n$ the tensor product of a representation  with Frobenius
characteristic $\TH_{1^k}$ by one of Frobenius characteristic $q^{n-k}\TH_{n-k}$.
A known result of representation theory  (see [20]) then yields that
$$
\Fch \BX \ses q^{n-k} \ssp \TH_{1^k}\ssp \TH_{n-k}
$$  
and 3.42 a) then follows 3.41 a). Similarly we derive 3.42 b) from 3.38 and
3.41 b).

We should note at this point that the identities we have established 
so far already yield an inductive mechanism for proving the $n!$ conjecture for hooks.
Indeed, making use of I.11 we immediately derive from 3.39 a) and b) that
$$
\eqalign{
a)\ess\del_{p_1}\ssp C_\mu &= C_{\mu/00}= (t+t^2+\cdots +t^{k})\ssp C_{\aaa} 
+ q^{n-k}\ssp \TH_{1^k}\ssp \TH_{n-k}+ 
(1+q+\cdots +q^{n-k-1})\ssp C_{\bbb} 
\cr
\cr
b)\ess \del_{p_1}\ssp C_\mu &= C_{\mu/00}= (1+t +\cdots +t^{k-1})\ssp C_{\aaa} 
+  t^{k}\ssp \TH_{1^k}\ssp \TH_{n-k}+ 
(q+q^2+\cdots +q^{n-k })\ssp C_{\bbb} 
\cr}
\eqno 3.43
$$
Now either of these two equalities yields the implication
$$
\dim \BM_\aaa=\dim\ssp \BM_\bbb=n!\ess\ess\Longrightarrow\ess\ess
\dim \BM_\mu=(n+1)!\ess .
\eqno 3.44
$$
In fact, applying $\del_{p_1}^n$ to both sides of 3.43 a) gives
(using the notation in I.4)
$$
F_{(n+1-k,1^k)}\ses
t\ssp [k]_t \ssp F_{(n+1-k,1^{k-1})} 
\sps q^{n-k}\ssp \Bigl({n\atop k}\Bigr) \ssp [k]_t!\ssp [n-k]_q!\sps 
[n-k]_q\ssp F_{(n -k,1^k)}  
\eqno 3.45
$$
with $\ssp [k]_t=1+\cdots +t^{k-1}\ssp ,\ess$  $\ssp [k]_t!=[{1}]_t[{2}]_t\cdots [k]_t$ and $\ssp [n-k]_q\ssp ,$
$\ssp [n-k]_q!\ssp $ analogously defined. Thus 3.44 follows from 3.45 by setting $t=q=1$. 
\sas

To prove 3.42 c)  we need a few auxiliary identities. To begin with note that subtracting 3.43 b) from 3.43 a)
we obtain that
$$
\TH_{1^k}\ssp \TH_{n-k}\ses { t^k-1 \over t^k-q^{n-k} }\ess C_\aaa
\sps { 1-q^{n-k} \over t^k-q^{n-k} }\ess C_\bbb\ess .
\eqno 3.46
$$
On the other hand, from suitably modified Macdonald Pieri rules (see [6] or [8])
we derive that  
$$
\TH_{1^k}\ssp \TH_{n-k}\ses { t^k-1 \over t^k-q^{n-k} }\ess \TH_{(n+1-k,1^{k-1})}
\sps { 1-q^{n-k} \over t^k-q^{n-k} }\ess \TH_{(n-k,1^{k })}\ess .
\eqno 3.47
$$
Finally,  subtracting 3.47 from 3.46 and recalling  that $\aaa=(n+1-k,1^{k-1})$ and $\bbb=(n-k,1^{k })$
we are led to the recursion
$$
\eqalign{
{ 1-q^{n-k} \over t^k-q^{n-k} }\ssp C_{(n-k,1^{k })}
&\ses
{ 1-q^{n-k} \over t^k-q^{n-k} }\ssp \TH_{(n-k,1^{k })}
\cr
&\ess\ess\ess\ess\ess\ess\ess\ess \sps { t^k-1 \over t^k-q^{n-k} }\ssp \TH_{(n+1-k,1^{k-1})}
\sms{ t^k-1 \over t^k-q^{n-k} }\ess C_{(n+1-k,1^{k-1})}\ess .
\cr
}
\eqno  3.48
$$
This enables us to prove 3.42 c) for each $\ssp n\ssp $ by induction on $\ssp k\ssp $.
Indeed, since $\BM_{(n+1)}$, by definition, is  the linear span of derivatives
of the Vandermonde  determinant in $(y_1,y_2,\cdots ,y_n)$ we necessarily have
$$
\Fch \BM_{(n+1)} \ses \TH_{(n+1)}\ess .  
$$
This gives 3.42 c) for $k=0$. However, if by induction, we assume 3.42 c) for $k-1$,
which is 
$$
C_{(n+1-k,1^{k-1})}\ses \Fch \BM_{(n+1-k,1^{k-1})} \ses \TH_{(n+1-k,1^{k-1})}\ess ,
$$
from 3.48 we immediately obtain that 
$$
  C_{(n-k,1^{k })}
\ses
  \TH_{(n-k,1^{k })} \ess .
$$
Thus 3.43 c) must hold true for all $k$ and our proof is complete.
\sas

\heading{\bol Remark 3.3}

We should point out the remarkable agreement that our 
conjectures  have with the theory of Macdonald polynomials.
To begin with note that substituting 3.46 in 3.43 a) or b) and carrying out the simplifications yields that 
$$
\del_{p_1}\ssp \TH_{(n+1-k,1^k)}\ses {q^{n-k}-t^{k+1}\over q^{n-k}-t^{k}}\ess {1-t^k\over 1-t}\ess \TH_{(n+1-k,1^{k-1})}
\sps {t^k-q^{n+1-k}\over t^k-q^{n-k} }\ess {\ess 1-q^{n-k}\over 1-q\ess \ess}\ess  \TH_{(n-k,1^{k })} 
$$ 
and this is precisely what may be obtained from I.13 and I.14. In the same vein,
we can show that 3.46 itself, which is an instance of higher order Pieri rules,
is in fact a consequence of Conjecture I.16 or the four term recursion
(which are the same  because of Theorem I.1). This can be seen from the following
formula which expresses Frobenius characteristics of atoms directly in terms of   
Macdonald Polynomials.
\sa

\heading{\bol Theorem 3.3}

{\ita Let $l$ and $a$ be the leg and arm of $(i,j)$,  let $\tau$ be the partition in
the shadow of $(i,j)$. As in the proof of Proposition I.8,
let $\ssp x_0^{ij},\ldots ,x_m^{ij}\ess ;\ess u_0^{ij},\ldots
,u_m^{ij}$  be the corner weights of
$\tau$,  
$\rho^{(1)},\rho^{(2)},\ldots ,\rho^{(m)}$ be the predecessors of $\tau$ ordered from
left to right so that $\ssp x_1^{ij},\ldots ,x_m^{ij}$ are the respective weights of the
cells $\tau/\rho^{(1)},\ldots ,\tau/\rho^{(m)}$. Set $\aaa^{(s)}=\mu-\tau+\rho^{(s)}$.
Then on the $C=\TH$ conjecture, we have
$$
{1\over q^a} A_{ij}^x\ses {1\over t^l} A_{ij}^y\ses\sum_{s=1}^m\ssp 
{
\prod_{r=1}^{m-1}(x_s^{ij}- u_r^{ij}) 
\over
\prod_{{r=1\,;\,r\neq s}}^{m }(x_s^{ij}- x_r^{ij})
}
\ess \TH_{\aaa^{(s)}}
	\ess .
\eqno  3.49
$$
 }

\heading{\bol Proof}

Our point of departure is the definition
$$
A_{ij}^x\ses C_{\mu/ij}-t\ssp C_{\mu/i+1,j}-C_{\mu/i,j+1}+t\ssp C_{\mu/i+1,j+1}
\eqno  3.50
$$
with the $C\ssp 's$ computed by means of formula 1.20, that is
$$
C_{\mu/ij}(x;q,t)\ses {1\over M}\ssp \sum_{s=1}^m\ssp {1\over x_s^{ij}}
\ssp {\prod_{r=0}^m \bigl(  x_s^{ij}- u_r^{ij}\bigr)\over \prod_{r=1\,;\,r\neq s}^m \bigl(   x_s^{ij}- x_r^{ij}\bigr) }
\ess \TH_{\aaa^{(s)}}\ess ,
\eqno  3.51
$$
where $M=(1-1/t)(1-1/q)$. For simplicity we shall assume that the shadows of
$(i,j),(i+1,j),(i,j+1)$ and $(i+1,j+1)$ contain the same corners of $\mu$.  
This given, note that for $s\neq 0$ we have the relations
$$
x_s^{ij}\ses t\ssp q\ssp  x_s^{i+1,j+1}
\ess \scs\ess\ess\ess\ess
x_s^{i+1,j}\ses  q\ssp  x_s^{i+1,j+1}
\ess \scs\ess\ess\ess\ess
x_s^{i,j+1}\ses t\ssp q\ssp  x_s^{i+1,j+1}\ess .
$$
Moreover, we recall that
$$
{1\over t}\ssp u_0^{ij}\ses {1\over t}\ssp u_0^{i,j+1}\ses   u_0^{i+1,j}\ses u_0^{i+1,j+1} \ess ,
$$
and
$$
{1\over q}\ssp u_m^{ij}\ses {1\over q}\ssp u_m^{i+1,j }\ses   u_m^{i,j+1}\ses u_0^{i+1,j+1} \ess .
$$
Using these relations in 3.51 written for $(i,j),(i+1,j),(i,j+1)$ and $(i+1,j+1)$,
we obtain from  3.50 that the coefficient of $\TH_{\aaa^{(s)}}$ in $A_{ij}^x$ is
$$
\eqalign{
A_{ij}^x\ess \big|_{\TH_{\aaa^{(s)}}} &= {CF\over M}\ssp
\biggl( 
x_s \ssp t \ssp q\ssp 
\Bigl(  1-{t u_0\over tq x_s} \Bigr)
\Bigl(  1-{q u_m\over tq x_s} \Bigr)
\cr
&
 \ess\ess\ess\ess\ess\ess\ess\ess\ess\ess\ess\ess\ess\ess\ess\ess
\sms
x_s \ssp t \ssp q\ssp 
\Bigl(  1-{ u_0\over  q x_s} \Bigr)
\Bigl(  1-{q u_m\over q x_s} \Bigr)
\cr
&
\ess\ess\ess\ess\ess\ess\ess\ess\ess\ess\ess\ess
\ess\ess\ess\ess\ess\ess\ess\ess\ess\ess\ess\ess\ess\ess\ess\ess\ess\ess\ess\ess\ess
\sms
x_s \ssp t  \ssp 
\Bigl( 1-{t u_0\over t  x_s} \Bigr)
\Bigl( 1-{  u_m\over t x_s} \Bigr)
\cr
&
\ess\ess\ess\ess\ess\ess\ess\ess\ess\ess\ess
\ess\ess\ess\ess\ess\ess\ess\ess\ess\ess\ess\ess
\ess\ess\ess\ess\ess\ess\ess\ess\ess\ess\ess\ess\ess\ess\ess\ess\ess\ess\ess\ess\ess
\sms
x_s \ssp t \ssp 
\Bigl(  1-{  u_0\over   x_s} \Bigr)
\Bigl(  1-{  u_m\over  x_s} \Bigr)
 \biggr)
\cr
}
\eqno 3.52
$$
where for convenience we have set
$$
x_s^{i+1,j+1}= x_s
\ess\scs\ess\ess\ess
u_0^{i+1,j+1}= u_0
\ess\scs\ess\ess\ess
u_m^{i+1,j+1}= u_m
$$
and
$$
CF\ses   
{\prod_{r=1}^{m-1} \bigl(  x_s^{ij}- u_r^{ij}\bigr)\over \prod_{r=1\,;\,r\neq s}^m \bigl(   x_s^{ij}- x_r^{ij}\bigr) }
\ess .
$$
Now a little manipulation simplifies 3.52 to
$$
A_{ij}^x\ess \big|_{\TH_{\aaa^{(s)}}}\ses CF\ess {x_s\ssp t\ssp ({1\over t}-1)(1-q)\over M}
\ess {u_m\over x_s}\ses CF\ssp q\ssp t\ssp u_m
$$
and this is 3.49 since
$$
t\ssp q\ssp u_m\ses t\ssp q \ssp u_m^{i+1,j+1}\ses  q^a \ess .
$$
this completes our proof.
\sa

Note that  for $\mu=(1^k,n+1-k)$,  
and $i=j=0$, formula 3.49 gives
$$
{1\over q^{n-k}}\ess A_{00}^x\ses 
{x_1^{00}-u_1^{00}\over x_1^{00}-x_2^{00} }\ess \TH_{(1^{k-1},n+1-k)}\sps 
{x_2^{00}-u_1^{00}\over x_2^{00}-x_1^{00} }\ess \TH_{(1^{k},n-k)} \ess .
\eqno 3.53
$$ 
Since in this case  
$$
x_1^{00}=t^k 
\ess\scs\ess\ess 
x_2^{00}=q^{n-k} 
\ess\ess \hbox{and } \ess\ess
u_1^{00}=1 \ess ,
$$
substituting this in 3.53 we get that
$$
{1\over q^{n-k}}\ess A_{00}^x\ses 
{t^k-1\over t^k-q^{n-k} }\ess \TH_{(n+1-k,1^{k-1})}\sps 
{q^{n-k}-1\over q^{n-k}-t^k }\ess \TH_{(n-k,1^{k} )} \ess ,
$$
which is in complete agreement with what we obtain by combining 3.42 a) with 
the Macdonald Pieri rule given in 3.47.

\vfill\supereject

\heading{\bol 4. Dimension bounds. }

In this section, we  derive a dimension bound for the
spaces $\BM_{\mu/ij}$. We  begin by reviewing the construction 
that yields the dimension bound of $n!$ for  $\BM_\mu$.
The reader is referred to [10] for proofs and further details.
\sas

Given a finite subset $S$ of  $n$-dimensional Cartesian space, we let $J_S$ 
denote the ideal of polynomials $P(x_1,x_2,\ldots ,x_n)$ which vanish on $S$.
The quotient ring $\BR_S=\BQ[x_1,x_2,\ldots ,x_n]/J_S$
may be viewed as the coordinate ring of the algebraic variety consisting 
of the elements of $S$. This given, it is clear that 
$$
\dim \BR_S\ses |S|\ess .
\eqno 4.1
$$
Although $\BR_S$ is not graded it has a filtration given by the subspaces $\CH_{\leq k}(\BR_S)$
spanned by the monomials $x^p=x_1^{p_1}x_2^{p_2}\cdots x_n^{p_n} $  which are of degree $\ess \leq k$.
A graded version of $\BR_S$ is obtained by setting 
$$
\gr \BR_S\ses \BQ[x_1,x_2,\ldots ,x_n]/\gr J_S 
\eqno 4.2
$$
with 
$$
\gr J_S\ses \bigl(\ssp h(P)\ssp: \ssp P\in J_S\ssp \bigr)
$$
where for a polynomial $\ssp P\ssp $ we let $h(P)$ denote  the  homogeneous component of $\ssp P\ssp $
that is of highest degree. It is also convenient to introduce the space $\BH_S=\bigl(\ssp \gr J_S\ssp \bigr)^\perp\ess ,$
the orthogonal complement of $\gr J_S$ with respect to the scalar product
$$
\langle P\scs Q\rangle \ses P(\del)Q(x)\ssp \big|_{x=0}\ess .
$$ 
We may also define $\BH_S$ as the space of polynomials that are killed by elements of $\gr \ssp J_S$ as differential
operators. In symbols
$$
\BH_S\ses \bigl\{Q(x)\ssp :\ssp P(\del)\ssp Q= 0\ess\ess  \forall \ess \ess P\in \gr J_S\ssp \bigr\}\ess .
\eqno 4.3
$$ 
It is easy to show (see [6]) that  any homogeneous basis
$\CB_S$  for $\BH_S$ is also a basis of $\gr \BR_S$ and $\BR_S$. 
In particular, the dimensions of these three spaces must be the same and thus equal to
$|S|$. In fact, we also have  for all $k\geq 0$ 
$$
\dim \CH_{\leq k}(\BR_S)\ses \sum_{s=0}^k\ssp \dim \CH_{=s}( \BH_S)\ses 
 \sum_{s=0}^k\ssp \dim \CH_{=s}(\gr \BR_S)\ess ,
\eqno 4.4
$$
where $\CH_{=s}(\BH_S)$ and $\CH_{=s}(\gr \BR_S)$ denote the subspaces
of $\BH_S$ and $\gr \BR_S $ consisting of their homogeneous elements
of degree $\ssp s\ssp$.

The natural action of $GL_n$ on polynomials $P(x_1,x_2,\ldots ,x_n)$ is defined by  setting
for an $n\times n$ matrix $A=\|a_{ij}\|_{i,j=1}^n$
$$
T_A\ssp P(x)\ses P(xA)
\eqno 4.5
$$
where $xA$ denotes matrix multiplication of the row vector $x=(x_1,x_2,\ldots ,x_n)$  
by $A$. It is not difficult to show that if $A$ is an orthogonal matrix, then
for all $P,Q\in \BQ[x_1,x_2,\ldots ,x_n]$ we have
$$
\LL T_A P\scs T_A Q\RR\ses \LL   P\scs   Q\RR
\eqno 4.6
$$
If $G$ is a group of $n\times n$ matrices that leave $S$ invariant then
both $J_S$ and $\gr J_S$ remain invariant under $T_A$ for every $A\in G$
 and we can define an action
of $G$ on the two quotient spaces $\BR_S$ and $\gr \BR_S$. It develops that the resulting $G$-modules 
are easily shown to be
equivalent. If in addition  $G$ consists of orthogonal matrices, then from 4.6
it follows that $\BH_S=(\gr J_S)^\perp$ is also $G$-invariant
and equivalent to $\gr \BR_S$ as a graded $G$-module. Moreover we have the following
character identity for all $k\geq 0$:
$$
\ch\CH_{\leq k}\bigl(\BR_S\bigr) \ses  
\sum_{s=0}^k\ssp \ch  \CH_{=s}( \BH_S)\ses 
 \sum_{s=0}^k\ssp \ch  \CH_{=s}(\gr \BR_S)\ess .
\eqno 4.7
$$	

Given a group $G$, the simplest $G$-invariant subsets are its ``{\ita orbits}.''
More precisely, for any point $\rho=(\rho_1,\rho_2,\ldots ,\rho_n)$, we set
$$
[\,\rho\, ]_G\ses \{\ssp \rho A\ssp :\ssp A\in G\ssp \}\ess .
\eqno 4.8
$$
Clearly, $G$ acts on the orbit  $[\,\rho\, ]_G$ as it does on the left cosets
of the subgroup that leaves $\rho$ invariant. It follows from this that
both $\BR_{[\,\rho\, ]_G}$ and $\gr \BR_{[\,\rho\, ]_G}$ afford 
this left coset action; in particular, if $\rho$ is
a regular point
(that is,  $\rho$ has a trivial stabilizer),
then $\BR_{[\,\rho\, ]_G}$ and $\gr \BR_{[\,\rho\, ]_G}$ 
are  versions of the left regular representation of $G$.
Moreover, if $G$ is a group of orthogonal matrices,
then $\BH_{[\,\rho\, ]_G}$  affords a graded version of the left regular representation of $G$ 
and consists of
polynomials that are killed by all $G$-invariant differential operators (see [6]).
In particular, all elements of $\BH_{[\,\rho\, ]_G}$ are harmonic.
\sas

To get our dimension bounds we need to suitably specialize $G$
and the point $\rho$. To this end, given $\mu=(\mu_1\geq\mu_2\geq\ldots \geq\mu_k>0)\part n\ess $ 
let $h=\mu_1$ be the number of parts of the conjugate of $\mu$
and let $(\aaa_1,\aaa_2,\ldots ,\aaa_k;\bbb_1,\bbb_2,\ldots ,\bbb_h)$
be distinct rational numbers. If preferred, the latter may be taken to be two additional sets of  
indeterminates. Recall that an injective  tableau $T$ of shape $\mu\part n$ is a
labeling of the cells of  $\mu$ by the numbers  $\{1,2,\ldots ,n\}$. The collection of
all such tableaux is denoted by ${\cal IT}(\mu)$.  Given a tableau $T\in {\cal IT}(\mu)\, $,
for each  $\, i=1,2,\ldots ,n\, $ we set
$$
a_i(T)=\aaa_r\ess\scs\ess\ess b_i(T)=\bbb_c
\eqno 4.9
$$
if the label $i$ is at the intersection of row $r$ with column $c$. The resulting
point of $2n$-dimensional space will be denoted by $\rho(T)$. In other words we set
$$
\rho(T)=\bigl(a_1(T),a_2(T), \ldots ,a_n(T);b_1(T),b_2(T), \ldots ,b_n(T)\bigr)\ess .
$$ 
For instance, for $\mu=(3,2)\ssp $ and
$$
T\ses 
\tableau[s]{ 
5 & 3  \cr
2 & 1 & 4}
$$
we set
$$
\rho(T)\ses (\aaa_1,\aaa_1,\aaa_2,\aaa_1,\aaa_2\, ;\, \bbb_2,\bbb_1,\bbb_2,\bbb_3,\bbb_1\ssp )\ess .
$$
Note that the collection
$$
\{\ssp \rho (T)\ssp :\ssp T\in {\cal IT}(\mu)\ess \}
\eqno 4.10
$$
consists of $n!$ distinct points. Indeed, since  the $\aa\ssp 's$ and the $\bb\ssp 's$  are
assumed to be  distinct, we can reconstruct the position of any label $\, i\, $ in $\, T\, $ by simply
looking at the $i^{th}$ and the $(n+i)^{th}$ coordinates of $\rho (T)$. Note that
the collection in 4.10 may also be viewed as an $S_n$-orbit under the  diagonal
action.
More precisely, we see that for any $T\in {\cal IT}(\mu)$ 
and $\sig=(\sig_1,\sig_2,\ldots ,\sig_n)\in S_n\,$, we have
$$
\sig   \rho (T) \ses 
\Bigl(a_{\sig_1}(T),a_{\sig_2}(T),\ldots ,a_{\sig_n}(T)\, ;\, 
b_{\sig_1}(T),b_{\sig_2}(T),\ldots ,b_{\sig_n}(T)\Bigr)
\ses \rho (\sig ^{-1}\ssp T )
	\ess ,
$$
where $\sig ^{-1}\ssp T$ is the tableau obtained by replacing the label $i$ in $T$
by the label $\sig ^{-1}_i$. This given, we can consider the collection in 4.10 as the $S_n$-orbit
of a point $\rho_\mu$ corresponding some specially chosen injective tableau of shape $\mu$.
To be specific we may let $T_0$ be the
``superstandard tableau'';
		this is the tableau  
		obtained by labeling the cells of $\mu\part n$ successively
		from $1,\ldots,n$ starting from the bottom row 
		and proceeding on up, from left to right in each row.
Set
$$
\rho_\mu\ses \rho(T_0)\ess .
\eqno 4.11
$$
We can thus apply the theory we have outlined at the beginning of the
section with $G$ specialized 
to the group of matrices yielding the diagonal action of $S_n$ 
and construct the three spaces 
$$
\BR_{ {[\rho_\mu]}}\ess \scs\ess\ess
\gr \BR_{ {[\rho_\mu]}}\ess \hbox{  and }\ess\ess
\BH_{[\rho_\mu]}
$$ 
where $[\rho_\mu]$ denotes the orbit of  $\rho(T_0)$  or, equivalently, 
the subset of $2n$-dimensional space defined by 4.10. We thus obtain three
left regular representations of $S_n$ and in particular we have
$$
\dim \BR_{ {[\rho_\mu]}}= 
\dim \gr \BR_{ {[\rho_\mu]}}=
\dim \BH_{[\rho_\mu]}= n!\ess .
\eqno 4.12
$$
The definition of these spaces suggests that they may depend on
our choice  of the $\aa_i's$ and $\bb_j's$. This is clearly the case for 
the coordinate ring $\BR_{[\pmu]}$. Nevertheless, there is strong
evidence that the space of harmonics $\BH_{[\pmu]}$ as well as the ideal $\gr J_{[\pmu]}$
and the quotient ring $\gr \BR_{[\pmu]}$ only depend on the choice of the partition $\mu$.
The reason for this stems from the following result:
\sas

\heading{\bol Proposition 4.1}

{\ita If $(i,j)$ is an outer corner cell of $\mu$ then for
any $s=1,2,\ldots , n$ the monomial $x_s^{i }y_s^{j }$ belongs to the ideal $\gr J_{[\pmu]}$.
In particular, if a monomial $x^py^q=x_1^{p_1}\cdots x_n^{p_n}y_1^{q_1}\cdots y_n^{q_n}$
does not vanish in $\gr R_{[\pmu ]}$ then all the pairs $(p_s,q_s)$ must give
 cells of  $\mu$. For the same  reason, every polynomial in $\BH_{[\pmu ]}$
must be a linear combination of monomials satisfying the same condition.}
\sas

\heading{\bol Proof}

This result was first proved in [10] (see Proposition 1.2 there). Since the argument
is quite simple and illuminating, we will include a proof here as well. To this end note that
the 
polynomial
$$
P_{(i,j)}(x,y)=\prod_{i'=1}^{i }(x_s-\aa_{i'})
\prod_{j'=1}^{j }(y_s-\bb_{j'})
$$
must necessarily vanish throughout $[\pmu]$. Indeed, for any 
point 
$$
\rho(T)=(a_1,\ldots ,a_n ; b_1,\ldots ,b_n)\in [\rho_\mu] 
$$ 
our definition gives that  $\, a_s =\aaa_{i'}\, $ for some $\ssp i'\leq i\ssp $ if
$\, s\, $ is south of $(i,j)$  in  $T\, $  and 
 $\, b_s =\bbb_{j'}\, $ for some $\ssp j'\leq j\ssp $
if $s$ is west of $(i,j)$. Since every cell of  $\mu$ satisfies at least one of these
conditions we see that at least one of the factors of $P_{(i,j)}$ must  necessarily vanish
for $(x;y)=(a_1,\ldots ,a_n ; b_1,\ldots ,b_n)\, $. This places $P_{(i,j)}$ in
$J_{[\rho_\mu]}$ and  its highest homogeneous component $x_s^iy_s^j$ in $\gr J_{[\rho_\mu]}$.
Thus every monomial which contains $x_s^iy_s^j$ as a factor 
must necessarily vanish in $\gr \BR_{[\rho_\mu] }$
and every polynomial in $\BH_{[\rho_\mu]}$ must be killed by $\del_{x_s}^i\del_{y_s}^j$. 
Since this must hold true for any $s=1,\ldots,n,$ we deduce  
that every element of $\gr \BR_{[\rho_\mu]}$ or $\BH_{[\rho_\mu]}$
must be a linear combination of monomials $x_1^{p_1}\cdots x_n^{p_n}y_1^{q_1}\cdots y_n^{q_n}$
where each pair $(p_s,q_s)$ must be a cell of  $\mu$. 
\sas

This result has the following immediate corollary
\sas

\noindent {\bol Theorem 4.1}

{\ita For any choice of the $\aaa_i$ and $\bbb_j$ we have the containment
$$
\BM_\mu\con \BH_{[\rho_\mu]}\ess .
\eqno 4.13
$$
In particular,
$$
\dim \BM_\mu\leq n!\ess .
\eqno 4.14
$$
Thus on the $n!$ conjecture we have
$$
\BM_\mu\ses\BH_{[\rho_\mu]}\ess\ess\ess \hbox{ and }\ess\ess\ess \gr J_{[\rho_\mu]}\ses I_{\DD_\mu} \ess ,
\eqno 4.15
$$
where $I_{\DD_\mu}$ denotes the ideal of polynomials that kill $\DD_\mu$.
}

\heading{\bol Proof}

These results were first proved in [10] (see Theorems 1.1 and 1.2 there).
We sketch the idea of the argument here. 
Since $\BH_{[\pmu]}$ affords a
version of the left regular representation of $S_n$,
it must contain a polynomial $\DD(x;y)$, unique up to a scalar factor, 
which alternates under the diagonal action. Clearly all the monomials
appearing in $\DD(x;y)$ must be of the form 
$$
x^py^q=x_1^{p_1}x_2^{p_2}\cdots x_n^{p_n}y_1^{q_1}y_2^{q_2}\cdots y_n^{q_n}
$$
with $(p_1,q_1),(p_2,q_2),\ldots ,(p_n,q_n)$ all distinct. On the other hand, Proposition 4.1 
guarantees that each of these pairs must give a cell of  $\mu$.
Combining these two facts yields that  the sequence 
$$
\bigl\{(p_1,q_1),(p_2,q_2),\ldots ,(p_n,q_n)\bigr\}
$$
must be a permutation of the cells of  $\mu$. Thus $\DD(x;y)$ can only be a multiple
of $\DD_\mu(x;y)$ and we must have
$$
\DD_\mu(x;y)\in \BH_{[\rho_\mu]}\ess .
\eqno 4.16
$$  
However, since $\BH_{[\rho_\mu]}$ is derivative closed, we must also have
$$
\BM_\mu\ses \CL_\del[\DD_\mu]\con \BH_{[\rho_\mu]}\ess ,
$$
proving 4.13. This completes our proof since 4.14 and 4.15 are immediate consequences of 4.13.
\sas

Now let $\mu\part n+1$ and $[\rho_\mu]_{ij}$ denote the subset of the orbit $[\rho_\mu]$ consisting of the 
points 
$\rho(T)$ corresponding to tableaux $T$ where $n+1$ lies in the shadow of the cell $(i,j)$.
Clearly the cardinality of this set is 
$$
\big|\ssp [\rho_\mu]_{ij}\ssp \big |\ses \nshadow(i,j)\ssp\times \ssp  n!
\eqno 4.17
$$
where ``$\nshadow(i,j)$'' denotes the number of cells that are in the shadow of $(i,j)$.
Moreover, it is easy to see that under the diagonal action of $S_n$, the set  $[\rho_\mu]_{ij}$
splits into as many as $\nshadow(i,j)$ distinct regular orbits. It follows then that  each of 
the three spaces
$$
\BR_{[\rho_\mu]_{ij}}   \scs\ess\ess\ess
\gr \BR_{[\rho_\mu]_{ij}}   \scs\ess\ess\ess
\hbox{and}
\ess\ess\ess\ess
\BH_{[\rho_\mu]_{ij}} \ess \scs\ess\ess
$$
breaks up into a direct sum of $\nshadow(i,j)$  regular representations of $S_n$. 
These observations yield the following extension of Theorem 4.1.
\sas

\heading{\bol Theorem 4.2}

{\ita For any choice of the $\aaa_i$ and $\bbb_j$ 
and any cell $(i,j)\in \mu$, we have the containment
$$
\CL_\del[\, \del_{x_{n+1}}^i\del_{y_{n+1}}^j\DD_\mu(x;y)\, ]\ess \con\ess  \BH_{[\rho_\mu]_{ij}}\ess .
\eqno 4.18
$$
In particular, 
$$
\dim \CL_\del[\, \del_{x_{n+1}}^i\del_{y_{n+1}}^j\DD_\mu(x;y)\, ]\leq \nshadow(i,j)\ssp \times \ssp n!\ess .
\eqno 4.19
$$
Moreover, equality here forces the equalities
$$
\CL_\del[\, \del_{x_{n+1}}^i\del_{y_{n+1}}^j\DD_\mu(x;y)\, ]= \BH_{[\rho_\mu]_{ij}}
\ess\ess\scs\ess\ess\ess\ess
\gr J_{[\rho_\mu]_{ij}}\ses I_{\del_{x_{n+1}}^i\del_{y_{n+1}}^j\DD_\mu}\ess ,
\eqno 4.20
$$
where $I_{\del_{x_{n+1}}^i\del_{y_{n+1}}^j\DD_\mu}$ denotes the ideal of polynomials
that kill $\del_{x_{n+1}}^i\del_{y_{n+1}}^j\DD_\mu\, $.
But then%
	\break
$\CL_\del[\del_{x_{n+1}}^i\del_{y_{n+1}}^j\DD_\mu(x;y)\, ]$
must necessarily break up into a direct sum of $\nshadow(i,j)$
regular representations of $S_n$.   
}

\heading{\bol Proof}

Note that if $P(x;y)$ is an element 
of the ideal $J_{[\rho_\mu]_{ij}}$ then the
polynomial
$$
Q(x;y)\ses P(x;y)\prod_{i'=1}^i\bigl (x_{n+1}-\aaa_{i'}\bigr)
\prod_{j'=1}^j\bigl (y_{n+1}-\bbb_{j'}\bigr)
$$
must necessarily vanish throughout the orbit $[\rho_\mu]$. In fact, $P(x;y)$ vanishes in 
$[\rho_\mu]_{ij}$ and the product
of the two remaining factors vanishes in the rest of  $[\rho_\mu]$. 
This places $Q(x;y)$ in $J_{[\rho_\mu]}\, $. Denoting as before 
by $h(P)$ and $h(Q)$ the highest homogeneous components of $P$ and $Q\,$, we derive that
$$
h(Q)= x_{n+1}^iy_{n+1}^j\ess h(P)\in \gr J_{[\rho_\mu]}\ess ,
$$
and therefore $h(Q)$ must kill all the elements of $\BH_{[\rho_\mu]}$. In particular,
in view of 4.16 we must also have 
$$
h(P)(\del)\ssp \del_{x_{n+1}}^i\del_{y_{n+1}}^j\DD_\mu\ses 0\ess . 
$$
Since this holds true for any $P\in J_{[\rho_\mu]_{ij}}$ we are brought to the conclusion that
$$
\gr J_{[\rho_\mu]_{ij}}\ess\con \ess  I_{\del_{x_{n+1}}^i\del_{y_{n+1}}^j\DD_\mu}\ess .
\eqno 4.21
$$
Now it is easy to show that
$$
I_{\del_{x_{n+1}}^i\del_{y_{n+1}}^j\DD_\mu}\ses \CL_\del [\del_{x_{n+1}}^i\del_{y_{n+1}}^j\DD_\mu]^\perp
	\ess .
$$
This gives
$$
\bigl(\ssp I_{\del_{x_{n+1}}^i\del_{y_{n+1}}^j\DD_\mu}
\ssp \bigr)^\perp\ses \CL_\del [\del_{x_{n+1}}^i\del_{y_{n+1}}^j\DD_\mu] 
$$
and thus 4.18 follows from 4.21 by taking orthogonal complements.
This given, 4.19 follows from 4.17 and 4.18  since 
$$
\dim \BH_{[\rho_\mu]}\ses \big|\, [\rho_\mu]\ssp \big|\ess .
$$ 
Finally, equality in 4.19 forces equality in 4.18 which in turn can only hold true 
if equality holds in 4.21. This completes our proof since the last assertion 
is a consequence of our preliminary observations.  
\sas

We are now in a position to derive the main result of this section.
\sas

\heading{\bol Theorem 4.3} 

{\ita For any $\mu\part n+1$ and any cell $(i,j)\in \mu$  we have
$$
\dim \BM_{\mu/ij}\ssp \leq \ssp \nshadow(i,j)\ssp \times \ssp n!\ess . 
\eqno 4.22
$$
Moreover, if equality holds here, 
then $\BM_{\mu/ij}$, breaks up into a direct sum of  $\nshadow(i,j)$ 
regular representations 
of $S_n$. 
}

\heading{\bol Proof}

In view of Theorem 4.2 we only need to show that 
$\ssp \BM_{\mu/ij}\ess  \hbox{and }\ess  
\CL_\del[\del_{x_{n+1}}^i\del_{y_{n+1}}^j\DD_\mu(x;y)\, ]
\ssp $ are equivalent as $S_n$-modules under the diagonal action. 
To this end note that from 1.16, we derive that:
$$
\eqalign{
{1\over i!}{1\over j_!}\del_{x_{n+1}}^i\del_{y_{n+1}}^j\DD_\mu(x;y)
&\ses 
\eee_{i ,j }\ssp
\DD_{\mu/i j }(\xon;\yon) \cr
&\ess\ess\ess\ess\ess\ess\ess
\sps  
\sum_{\multi{(i',j')\ssp \in \ssp \mu\cr i'>i\ssp \hbox{~or~} \ssp j'>j}} x_{n+1}^{i'-i} \;  y_{n+1}^{j'-j}\ess
c_{i',j'}\ssp
\DD_{\mu/i'j'}(\xon;\yon)
\cr 
}
$$
\vskip -.3truein
\hfill 4.23
\vskip  .2truein
\noindent
where $\eee_{ij}= \pm 1$ and the $c_{i',j'}$ are suitable constants.
Thus for any $f\in \BQ[\xon ;\yon]$  
we necessarily have
$$
a)\ess\ess f(\del_x;\del_y)\ssp \del_{x_{n+1}}^i\del_{y_{n+1}}^j\DD_\mu(x;y)\ses 0
\ess\ess\ess\ess \longleftrightarrow \ess\ess\ess\ess 
b)\ess\ess f(\del_x;\del_y)\ssp \DD_{\mu/ij}\ses 0\ess .
$$
In fact, we see from  4.23 that b)  immediately follows from a) by setting 
$x_{n+1}=y_{n+1}=0$. Conversely, if b) holds true then applying to it the operator
$$
D_{i'-i,j'-j}\ses \sum_{s=1}^n\ssp \del_{x_s}^{i'-i}\del_{y_s}^{ j'-j}
$$
we obtain that 
$$
f(\del_x;\del_y)\ssp \DD_{\mu/i',j'}\ses 0
$$must hold true for all $(i',j')\in \mu$ that are in the shadow of $(i,j)$
and this forces a) to hold true as well.

Now from the relations in 1.13 it follows that 
$$
\eqalign{
\del_{x_{n+1}}\ssp\bigl(\ssp  \del_{x_{n+1}}^i\del_{y_{n+1}}^j\DD_\mu(x;y)\ssp \bigr)
&\ses  -\sum_{s=1}^n\del_{x_s}\ssp\bigl(\ssp \del_{x_{n+1}}^i\del_{y_{n+1}}^j\DD_\mu(x;y) \ssp  \bigr) \cr
\del_{y_{n+1}}\ssp\bigl(\ssp  \del_{x_{n+1}}^i\del_{y_{n+1}}^j\DD_\mu(x;y)\ssp \bigr)
&\ses  -\sum_{s=1}^n\del_{y_s}\ssp\bigl(\ssp \del_{x_{n+1}}^i\del_{y_{n+1}}^j\DD_\mu(x;y) \ssp  \bigr) \ess .\cr
}
$$
This means that we can construct a basis for 
$$
\CL_\del[\ssp \del_{x_{n+1}}^i\del_{y_{n+1}}^j\DD_\mu(x;y)\ssp ]
$$ 
of the form
$$
\CB_{ij}\ses \bigl\{\ssp b(\del_x;\del_y)\ssp \del_{x_{n+1}}^i\del_{y_{n+1}}^j\DD_\mu(x;y)\ssp :\ssp b\in \CC\ssp \bigr\}
$$ 
with $\CC$ a collection of  polynomials in the variables $\xon;\yon$. 
But then it follows from the observations above that, with the same $\CC$,
the collection  
$$
\CB_{ij}^*\ses \bigl\{\ssp b(\del_x;\del_y) \DD_{\mu/ij}(x;y)\ssp :\ssp b\in \CC\ssp \bigr\}
$$ 
must give a basis for $\BM_{\mu/ij}$. This given, if the elements of $\CC$ are 
chosen to be homogeneous, it follows that the action of $S_n$ on the corresponding homogeneous components
of $\CB_{ij}$ and $\CB_{ij}^*$ must be given by exactly the same matrices, proving that 
$\BM_{ij}$ and $\CL_\del [\ssp \del_{x_{n+1}}^i\del_{y_{n+1}}^j\DD_\mu(x;y)\ssp ]$
must be equivalent also as {\ita graded} $S_n$-modules. This completes our argument.

\vfill \supereject

\def\skewp(#1){{\def\Pscale{.6}\skewptn(#1)}}
\def\twtw{\skewp(2|2)}
\def\ofr{\skewp(1|4)}
\def\snake{\skewp(2|0,1,3)}
\def\otr{\skewp(1|3)}
\def\echair{\skewp(2|0,1,2)}
\def\trtw{\skewp(2|3)}
\def\oooo{\skewp(1|1|1|1)}
\def\ootwtxo{\skewp(1|1|2|0,2,1)}
\def\ootr{\skewp(1|1|3)}
\def\ootw{\skewp(1|1|2)}
\def\otwtw{\skewp(1|2|2)}
\def\oootw{\skewp(1|1|1|2)}
\def\tchair{\skewp(1|2|0,1,1)}
\def\threeone{\skewp(1|1|1|0,1,1)}
\def\twotwo{\skewp(1|1|0,1,1|0,1,1)}
\def\osq{\skewp(1)}
\def\ulchair{\skewp(1|1|3|0,2,1)}
\def\lonchair{\skewp(1|1|2|0,1,1)}
\def\chair{\skewp(1|3|0,2,1)}
\def\otw{\skewp(1|2)}
\def\vtw{\skewp(1|1)}
\def\vtwxtwo{\skewp(1|1|0,1,2|0,2,1)}

\heading{\bol 5. Atoms and further lattice diagram characteristics. }

In [8] Garsia and Haiman call two lattice  diagrams $D_1$ and $D_2$  ``equivalent'' and write
$D_1\approx D_2$ if and only $D_2$ can be obtained from $D_1$ by a sequence of row and column
rearrangements.    
Diagrams that are equivalent to skew diagrams are briefly referred to there  
as ``gistols.'' We should note that it is not visually obvious 
when two diagrams are equivalent. For instance we have
$$
\skewp(1|0,1,1|3) \;\approx\; \skewp(0,1,1|3|1) \;\approx\;
        \skewp(0,2,1|3|1) \;\approx\; \skewp(1|3|0,2,1)
$$
Following standard convention, the ``conjugate'' of a diagram $D$,
denoted by $D'$ is the diagram obtained  by reflecting $D$ across the diagonal line $x=y$.
Similarly, the reflection of a lattice square $s=(i,j)$ across $x=y$  is denoted by $s'=(j,i)$. 
Finally, if $D$ may be decomposed into the union of
two diagrams $D_1$ and $D_2$ in such a manner that no square of $D_2$ is in the same row or 
column of a square of $D_1$, then we shall say that $D$ is ``decomposable'' and we write $D=D_1\times
D_2$.
This given, Garsia-Haiman postulate the existence of a family
of polynomials $\{G_D(x;q,t)\}_D$, and a family of weights $\, w_{ s,D}(q,t)\, $,
 with the following basic properties:
$$
\cases{
(0)\ess\ess G_{D}(x;q,t)=\TH_\mu (x;q,t)& if $\ess D$ is the diagram of $\mu$\cr
\ess &$\ess$\cr
(1)\ess\ess G_{D_1}(x;q,t)=G_{D_2}(x;q,t)& if $\ess D_1\approx D_2$\cr
\ess &$\ess$\cr
(2)\ess\ess G_{D_1}(x;q,t)=G_{D_2}(x;t,q) &if $\ess D_2\approx D_1'$\cr
\ess &$\ess$\cr
(3)\ess\ess G_{D}(x;q,t)=G_{D_1}(x;q,t) G_{D_2}(x;q,t) & if $\ess D\approx D_1\times D_2$\cr\cr 
 (4)\ess\ess\del_{p_1}\ssp G_{D}(x;q,t)\ses \sum_{s\in D}\ess w_{ s,D}(q,t)\ess  G_{D/s}(x;q,t)\ess ,
& with  $\ess D/s\ssp  =\ssp D\ssp $  minus  $s\ess .$\cr}
\eqno 5.1
$$
It should be noted at the onset that these properties overdetermine the family $\{G_D(x;q,t)\}_D$,
so that existence is by no means guaranteed. Nevertheless, all the experimentations so far indicate that
the existence of such a family is consistent with the theory of Macdonald polynomials. 
In particular it was shown in [8] that for any partition $\mu$ we have
$$
\TH_{\mu'}(x;q,t)\ses \TH_{\mu}(x;t,q) 
$$ 
which is in perfect agreement with  condition (2) in 5.1. 

Experimentation suggests that the  weights $\, w_{ s,D}(q,t)\, $ should be monomials in $q,t$,
but there are no conjectured formulas for general lattice diagrams. Nevertheless, we should
point out that if condition (4) holds for the conjugate $D'$ of a diagram $D$,
that is we have
$$
\del_{p_1}\ssp G_{D'}(x;q,t)\ses \sum_{s'\in D'}\ess w_{ s',D'}(q,t)\ess  G_{D'/s'}(x;q,t)\ess ,
\eqno 5.2
$$ 
then, upon interchanging $q$ and $t$, from condition (2) we immediately derive that
we must also have
$$
\del_{p_1}\ssp G_{D }(x;q,t)\ses \sum_{s \in D }\ess w_{ s',D'}(t,q)\ess  G_{D /s }(x;q,t)
	\ess .
\eqno 5.3
$$ 
Thus the conditions in 5.1 force the existence of at least one pair of ``weights'' both yielding the expansion
in 5.1 (4). Now, in the case that $D$ is a skew diagram, representation 
theoretical reasons  suggest that we should use
either one of the following two choices of weights:
\def \spand {\ess\ess\ess \hbox{and }\ess\ess\ess }

$$
a)\ess\ess w[s,D]=t^{l _D(s)}q^{a'_D(s)}
\spand
b)\ess\ess w[s,D]=t^{l'_D(s)}q^{a_D(s)}
\eqno 5.4
$$
where as customary $ l_D(s), l'_D(s)$ denote the number  of cells  
strictly north and south, respectively, of 
$s$ in $D$, and likewise $ a^{\ } _D(s), a'_D(s)$ give the number of cells 
strictly east and west, respectively.
It is easy to see that is  consistent with
the relations given in 5.2 and 5.3. Using these weights, 
we can determine a wide variety of the  polynomials $G_D$,
and each via a number of different independent ways all leading to the same 
final Schur function expansion. Remarkably, all the polynomials 
thus obtained reduce to $h_1^{|D|}$ when we set $t=q=1$. In particular, when $D$ is a skew diagram
or a diagram obtained by removing a cell from a Ferrers diagram, we
invariably obtain an expansion of the form
$$
G_D(x;q,t)\ses \sum_{\la\part n}\ssp S_\la(x)\ssp \TK_{\la ,D}(q,t)
\eqno 5.5
$$
with $\TK_{\la ,D}(q,t)$ polynomials with nonnegative integral coefficients
satisfying
$$ 
\TK_{\la,D}(1,1)\ses f_\la \ses \#\{\ssp \hbox{standard tableaux of shape $\la$}\ssp \}\ess . 
\eqno 5.6
$$
Even more remarkably, all the identities involving Macdonald 
polynomials we have been able to derive  by means of the rules in 5.1 end
up to be computer  verifiable  and/or theoretically provable.
\sas

To get our point across, it will be good to review some of these calculations here.
As a first example, we shall apply  rule (4)  to the diagram 
$D=\{(0,2),(1,0),(1,1),(1,2),(2,0),(3,0)\}$.
In the figure below, the first tableau is obtained by filling the cells of $D$ with
the weights computed according to formula a) of 
5.4 and the second according to formula b).
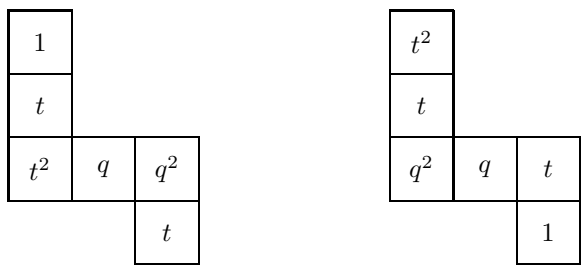
$$
\tableau{ 1 \cr
          t \cr
          t^2 & q & q^2 \cr
          \bl & \bl & t }
\hskip 1in
\tableau{ t^2 \cr
          t \cr
          q^2 & q & t \cr
          \bl & \bl & 1 }
$$
Thus, if we use the diagrams themselves to represent the corresponding
polynomial, rule (4) with the first set of  weights gives
$$
\del_{p_1}\ssp{\ulchair} \ses (1+t)\ess\chair \sps
t^2\ess \vtwxtwo\sps q\ess \lonchair \sps q^2\ess \ootwtxo
\sps t\ess \ootr
$$
while the second set gives
$$
\del_{p_1}\ssp {\ulchair} \ses t(1+t)\ess\chair \sps
q^2\ess \vtwxtwo\sps q\ess \lonchair \sps t\ess \ootwtxo
\sps  \ootr\ess \ess .
$$
We can thus obtain by subtraction that
$$
\chair \ses {t^2-q^2\over t^2-1}\ess
\vtw \times \otw
\sps {q^2-t\over t^2-1}\ess\ootw\times \osq\sps {t-1\over t^2-1}\ess\ootr\ess .
\eqno 5.7
$$
It easily obtained, either by computer or by means of Macdonald  Pieri rules, that
$$
\ootw \times \osq  \ses   {(1-t)(q-t^3)\over (q-t)(q^2-t^3) }\ess\ootr
\sps
{(1-t^2)(q-1)\over (q-t^2)(q-t) }\ess \otwtw 
\ess\sps
{(q-1)(q^2-t^2)\over (q-t^2)(q^2-t^3) }\ess \oootw 
\eqno 5.8
$$
and
$$
\eqalign{\vtw \times \otw \ses  & {(1-t^2)(1-t)\over (q^2-t^2)(q-t) }\ess\trtw 
\sps
{(1-t^2)(q-1)(q-t^2)\over (q-t)^2(q^2-t^3) }\ess \ootr 
\cr
&\ess\ess\ess\ess\ess\ess\ess\ess\sps
{(1-t^2)(q-1)(q^2-t)\over (q-t^2)(q-t)(q^2-t^2) }\ess \otwtw 
\sps
{(q-1)(q^2-t)\over (q-t^2)(q^2-t^3) }\ess \oootw
\cr}
\eqno 5.9
$$
Using 5.8 and 5.9 in 5.7 yields the surprisingly simple final expression
$$
\chair \ses {(1-t)\over (q-t) }\ess  \trtw \sps {(q-1)\over (q-t) }\ess \ootr
\eqno 5.10
$$ 
Applying rule (4) to the diagram $\{ (0,1),(1,0),(1,1),(2,0),(3,0)  \}$
according to the weights
$$
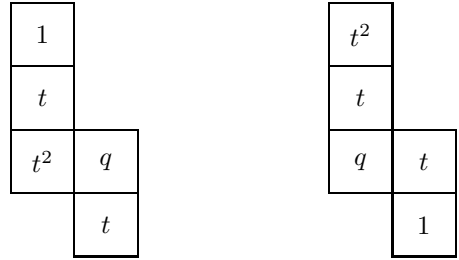
\tableau{1 \cr
	 t \cr
	 t^2 & q \cr
	 \bl & t}
	\hskip 1in
\tableau{t^2 \cr
	 t \cr
	 q & t \cr
	 \bl & 1}
$$ 
gives
$$
{\del_{p_1} \atop \ }\ssp\lonchair\ses (1+t)\ssp
\tchair\sps t^2\ssp \twotwo\sps q\ess \threeone \sps t\ess \ootw 
$$
and
$$
{\del_{p_1} \atop \ }\ssp\lonchair\ses (t^2+t)\ssp
\tchair\sps q\ssp \twotwo\sps t\ess \threeone \sps 1\ess \ootw 
$$
so by subtraction we get
$$
\tchair\ses {q-t^2\over 1-t^2}\ess \twotwo \sps  {t-q\over 1-t^2}\ess \threeone \sps {1-t\over 1-t^2}\ess \ootw
\eqno 5.11
$$
Using Pieri rules again gives
$$
\twotwo\ses   {(1-t)(1-t^2)\over (q-t)(q-t^2) }\ess\twtw 
\sps
{(1-t)(q-1)(1+t)^2\over (q-t)(q-t^3) }\ess \ootw 
\sps
{(q-t)(q-1)\over (q-t^2)(q-t^3)}\ess \oooo
\eqno 5.12
$$
and
$$
\threeone\ses {t^3-1\over t^3-q }\ess \ootw \sps {1-q\over t^3-q }\ess \oooo \ess .
\eqno 5.13
$$
Substituting 5.12 and 5.13 in 5.11 produces
$$
\tchair\ses {1-t\over q-t }\ess  \twtw \sps {q-1\over q-t }\ess \ootw\ess ,
\eqno 5.14
$$
leaving us with the puzzle of explaining why the coefficients we get here are the same we got
in 5.10.
\sas

But we have more surprises coming. We have yet another path that yields an expression for the
polynomial indexed by the diagram
$$ 
D\ses \chair \ess .
\eqno 5.15
$$
This is based on applying rule (4) to the following diagram.
$$
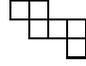
\skewp(2|0,1,3|0,3,1)
$$
Omitting the details, the resulting expansion turns out to be 
$$
\chair\ses {t-q\over 1-q}\ess \osq \times \otr \sps {q^2-t\over 1-q}\ess  
\skewp(2) \times \otw
\sps {t-q^2\over 1-q}\ess \echair \times \osq
\sps {1-t\over 1-q}\ess \snake
\eqno 5.16
$$
Note next that applying rule (2) we can transform 5.14 into
the expansion
$$
\echair\ess \ses {(1-q)\over (t-q) }\ess  \twtw \sps {(t-1)\over (t-q) }\ess \otr\ess .
\eqno 5.17
$$
Omitting again the details, we can show that
$$
\snake\ses {(1-q^2)\over (t-q^2) }\ess  \trtw \sps {(t-1)\over (t-q^2) }\ess \ofr\ess .
\eqno 5.18
$$
Now miraculously, after we substitute 5.17 and 5.18 into 5.16, apply the required Pieri rules
and feed the rather monstrous result into the computer we witness the occurrence of
massive simplifications yielding that 5.16 is yet another way of writing 5.10.
\sas

The reader may find it amusing to play this game by means of Stembridge's
SF Maple package. Seeing is believing that there must be a beautiful 
explanation for all these miraculous identities. 
Now it develops that we can use the present theory to remove the mystery out
of some of them. To see this we start by writing the identities in I.19 
in the form
$$
\matrix{
a)_x \ess\ess C_{\mu/ij}= K_{ij}^x \sps t\ssp C_{\mu/i+1,j}\hfill
&\ess\scs  &
a)_y \ess\ess C_{\mu/ij}= K_{ij}^y \sps q\ssp C_{\mu/i,j+1}
\cr\cr
b)_x\ess\ess K_{ij}^x= A_{ij}^x \sps K_{i,j+1}^x \hfill
&\ess\scs&
b)_y\ess\ess K_{ij}^y= A_{ij}^y \sps K_{i+1,j}^x \hfill
\cr}
\eqno 5.19
$$
Iterating $\ssp a)_x$ and using the fact that $C_{\mu/i,\,  \mu_{j+1}'}^x=\{0\}$
we obtain
$$
C_{\mu/ij}^x=K_{ij}^x+t\, K_{i+1,j}^x+\cdots +t^{l_{ij}}K_{  \mu_{j+1}'-1\, , j}^x 
\eqno  5.20
$$
where $\, l_{ij}=\mu_{j+1}'-i-1\ssp $ is the leg of the cell $(i,j)$.
Similarly from  $\ssp  b )_x$, using  $K_{i,\mu_{i+1}}^x=\{0\}$
we derive that
$$
K_{ij}^x=A_{ij}^x+A_{i,j+1}^x+\cdots +A_{i,\mu_{i+1}-1}^x 
	\ess .
\eqno  5.21
$$
Taking account of 3.47 let us set for each cell $s=(i,j)$
$$
\Xi_{\mu,(i,j)}\ses {1\over q^a} A_{ij}^x\ses {1\over t^l} A_{ij}^y\ses\sum_{s=1}^m\ssp 
{
\prod_{r=1}^{m-1}(x_s^{ij}- u_r^{ij}) 
\over
\prod_{{r=1\ssp   r\neq s}}^{m }(x_s^{ij}- x_r^{ij})
}
\ess \TH_{\aaa^{(s)}}
	\ess .
\eqno  5.22
$$
This given, we may rewrite 5.21 in the form
$$
K_{ij}^x=\sum_{(i,j)\RA \, s'}\ssp q^{a(s')}\ssp \Xi_{\mu,s'}
	\ess ,
$$
where we have used the symbol ``$\ssp (i,j)\RA \, s'\ssp $'' to indicate  that we
are to sum over  cells $s'$ that are directly east of $(i,j)$ including $(i,j)$ itself
and $a(s')$ denotes the arm of $s'$ in $\mu$.
Using such an expression for each of the characteristics $K_{i',j}^x$ occurring in 5.20,
we derive that
$$
C_{\mu/ij}\ses \sum_{(i,j)\leq \ssp (i',j')=s'}\ssp t^{i'-i}\ssp q^{a (s')}\ssp \Xi_{\mu,s'}
\eqno 5.23
$$
where ``$(i,j)\leq \ssp (i',j')$'' is to represent that we are to sum over 
all cells $(i',j')\in\mu$ that are in the shadow of $(i,j)$. Denoting
the partition in the shadow of $(i,j)$ by $\tau_{ij}$, we see that 5.23 
may be rewritten as
$$
C_{\mu/ij}\ses \sum_{(i,j)\leq\ssp   s'}\ssp t^{l'_{\tau_{ij}}(s')  }\ssp q^{a_{\tau_{ij}} (s')}\ssp \Xi_{\mu,s'}
\ess .
\eqno 5.24
$$
On the other hand, we can derive from the recurrence in 5.19 $a)_y$ that we also
have
$$
C_{\mu/ij}= K_{ij}^y+q\, K_{i ,j+1}^y+\cdots +q^{a_{ij}}K_{ i\,  , \, \mu_{i+1} -1}^y 
\eqno  5.25
$$ 
where $a_{ij}=\mu_{i+1}-j-1$ is the arm of $(i,j)$ in $\mu$.
Moreover, from 5.19 $b)_y$ we derive
$$
K_{ij}^y=A_{ij}^y+A_{i+1,j }^y+\cdots +A_{\mu_{j+1}'-1\, ,\, j}^y \ess .
\eqno  5.26
$$
Proceeding as we did above, the identities in 5.22, 5.25 and 5.26 now yield 
that we also have
$$
C_{\mu/ij}\ses \sum_{(i,j)\leq\ssp   s'}\ssp t^{l _{\tau_{ij}}(s')  }\ssp q^{a'_{\tau_{ij}} (s')}\ssp \Xi_{\mu,s'}
\ess .
\eqno 5.27
$$
Since  on the $C=\TH$ conjecture we have (see  I.12) 
$$
C_{\mu/00}\ses \del_{p_1} \ssp \TH_\mu \ess,
$$
we get that the special cases $(i,j)=(0,0)$ of 5.24 and 5.27 yield
$$
 \del_{p_1} \ssp \TH_\mu \ses \sum_{  s\in \mu}\ssp t^{l'_\mu (s )  }\ssp q^{a_\mu  (s )}\ssp \Xi_{\mu,s }
\ess ,
\eqno 5.28
$$
and
$$
 \del_{p_1} \ssp \TH_\mu \ses 
 \sum_{  s\in \mu}\ssp t^{l_\mu  (s )  }\ssp q^{a'_\mu  (s )}\ssp \Xi_{\mu,s }
\ess .
\eqno 5.29
$$
Comparing with 5.1 (4) and 5.3 written for $D=\mu$ and with   $w(s,D)$ and $w(s',D')$
respectively given by the weights in 5.4 a) and b) we come to the inescapable conclusion
that at least for $D$ the diagram of a partition, these mysterious polynomials  $G_{D/s}(x;q,t)$
must be none other than our normalized atom characteristics $\Xi_{\mu,s}\ssp $. 
To be precise, we are thus led to the addition of one further rule to the heuristic apparatus 
exhibited in  5.1, namely that we must also have
$$
(5)\ess\ess\ess\ess  G_{\mu/s} \ses \Xi_{\mu,s}\bigsp\bigsp (\ess \forall\ess\ess s\in \mu\ssp )
	\ess .
\eqno 5.30
$$
It develops that accepting this hypothesis, we can easily explain a wide variety of 
formulas that may be derived from
the rules in 5.1. This is best seen through a few examples. Let us begin with 5.10
which heretofore could only be obtained through the two intricate paths we illustrated above.
Now, we saw at the beginning of the section that we have the equivalence 
$$
\skewp(1|0,1,1|3) \;\approx\; \skewp(1|3|0,2,1) \ess .
$$
Thus from rule (1) and formula 5.22 for $\mu=(3,2,1)$ and $s=(1,0)\,$, we obtain the expansion
$$
\chair \ses {x_1-u_1\over x_1-x_2}\ess \trtw \sps  {x_2-u_1\over x_2-x_1}\ess \otwtw
\eqno 5.31
$$
where $x_1$ and $x_2$ must be the weights of the two corners of the partition
$\;\otw\;$ (which is the shadow of $(1,0)$ in $(3,2,1)$ ), and $u_1$ must be 
the weight of the inner corner. We thus deduce that 5.31 must hold true with
$$
x_1=t\ess \ess \scs\ess\ess\ess
x_2=q\ess \ess \scs\ess\ess\ess
u_1=1
	\ess .
\eqno 5.32
$$
Now we can easily see that making these substitutions in 5.31 immediately yields 
our formula 5.10. For our next example we take 
$$
D\ses \tchair\ess .
$$
In this case we use 5.30 with $\mu=(2,2,1)$ and $s=(0,0)$, obtaining 
that we must have 
$$
\tchair \ses 
{x_1-u_1\over x_1-x_2}\ess \twtw \sps  {x_2-u_1\over x_2-x_1}\ess \ootw\ess .
\eqno 5.33
$$
Here we must take
$$
x_1=t^2\ess \ess \scs\ess\ess\ess
x_2=tq\ess \ess \scs\ess\ess\ess
u_1=t\ess \ess .
\eqno 5.34
$$
Making these substitutions we see that 5.33 reduces to 5.14. The the fact that the
weights in 5.14 are the same as those in 5.10 may be explained from the equivalence
$$
\skewp(1|0,1,1|2) \;\approx\; \skewp(1|2|0,1,1)
	\ess ,
$$
which shows that we could also use 5.30 with $\mu=(2,2,1)$ and $s=(1,0)$, yielding
that we must also have 5.33 with the weights given in 5.32. 
\sas
 
\heading{\bol Remark 5.1}

Incidentally,
the reason that the weights in 5.34 yield the same result  as those in 5.32 is due to a special instance of 
Theorem I.3 stated in the introduction. In fact, we should note that we can easily deduce from formula 5.22 itself that
the normalized atom characteristics $\Xi_{\mu,s}$ must remain constant within any of the
rectangles defined in I.51.
\sa

Formula 5.30 can also yields the expansion in 5.18. Indeed if we use it with
$\mu=(4,2)$ and $s=(0,0)$ we immediately derive that
$$
\snake\ses  
{x_1-u_1\over x_1-x_2}\ess \trtw \sps  {x_2-u_1\over x_2-x_1}\ess \ofr
\eqno 5.35
$$
with
$$
x_1=t\, q\ess \ess \scs\ess\ess\ess
x_2=q^3\ess \ess \scs\ess\ess\ess
u_1=q\ess \ess \scs\ess\ess\ess
$$
or with the weights
$$
x_1=t \ess \ess \scs\ess\ess\ess
x_2=q^2\ess \ess \scs\ess\ess\ess
u_1=1\ess \ess \scs\ess\ess\ess
$$
because of the equivalence
$$
\skewp(2|0,1,3) \;\approx\; \skewp(2|1,1,2) \ess .
$$ 
\sas

We should point out that we haven't proved anything here, since
a number of the above derivations are based on various yet unproven conjectures.
Nevertheless, the variety of identities that may be constructed in this manner  
should be taken as evidence in support of the conjectures. More importantly,
these calculations open up a number of avenues for further investigation.
To begin with, it is difficult to believe that
we could not find some very natural quotients
of subspaces of the modules $\BM_\mu$ whose Frobenius characteristics
may be identified with the conjectured polynomials $G_D(x;q,t)$
(as we have done for the polynomials $A_{ij}^x$ and $A_{ij}^y$). 
Our experience suggests that these subspaces should result from  
restricting to smaller and smaller Young subgroups of $S_n$.
In this vein, just as the characteristics $\Xi_{\mu,s}$ do extend and simplify
the Macdonald (first order) Pieri rules, we would expect that ,using the general
polynomials $G_D(x;q,t)$, we should be able to unravel the combinatorics of 
higher order Pieri rules. From this point of view it appears that we have
uncovered what may be the tip of an iceberg of further research.
Only time will tell the significance of what may  ultimately be found 
in explaining some of the mysteries that stem from the present developments.

\vfill\supereject

\centerline {\bol REFERENCES}
\sa
{\frenchspacing 
\parindent=.25truein

\item{[1]}
F.~Bergeron and A.~Garsia, {\ita {S}cience {F}iction and {M}acdonald
  {P}olynomials}, Algebraic methods and $q$-special functions (L.~Vinet
  R.~Floreanini, ed.), CRM Proceedings \& Lecture Notes, American Mathematical
  Society, to appear.

\item{[2]}
N.~Bergeron and A.~M. Garsia, {\ita On certain spaces of harmonic polynomials},
  Hypergeometric functions on domains of positivity, Jack polynomials, and
  applications (Tampa, FL, 1991), Contemp. Math., vol. 138, Amer. Math. Soc.,
  Providence, RI, 1992, pp.~51--86.

\item{[3]}
C.~Chang, {\ita Geometric interpretation of the Macdonald polynomials and the
  $n!$ conjecture}, Ph.D.\ thesis, University of California, San Diego, 1998.

\item{[4]}
A.~Garsia, {\ita {R}ecent {P}rogress on the {M}acdonald $q,t$-{K}ostka
  conjecture}, Actes du $4^e$ {C}olloque sur les {S}\'eries {F}ormelles et
  Combinatoire Alg\'ebrique, {U.Q.A.M.} (P.~Leroux and C.~Reutenauer, eds.),
  Publications du laboratoire de combinatoire et d'informatique math\'ematique,
  vol.~11, Universit\'e du Qu\'ebec \`a Montre\'eal, 1992, pp.~249--255.

\item{[5]}
A.~Garsia and M.~Haiman, {\ita A random $q,t$-hook walk and a {S}um of {P}ieri
  {C}oefficients}, J. Combin. Theory Ser. A {\bol 82} (1998), no.~1, 74--111.

\item{[6]}
A.~Garsia and M.~Haiman,
  {\ita {O}rbit {H}armonics and {G}raded {R}epresentations}, Laboratoire
  de combinatoire et d'informatique math\'ematique, Universit\'e du Qu\'ebec
  \`a Montr\'eal, in preparation.

\item{[7]}
A.~M. Garsia and M.~Haiman, {\ita A graded representation model for
  {M}acdonald's polynomials}, Proc. Nat. Acad. Sci. U.S.A. {\bol 90} (1993),
  no.~8, 3607--3610.

\item{[8]}
A.~Garsia and M.~Haiman,
  {\ita Factorizations of {P}ieri rules for {M}acdonald polynomials},
  Discrete Math. {\bol 139} (1995), no.~1-3, 219--256, Formal power series and
  algebraic combinatorics (Montreal, PQ, 1992).

\item{[9]}
A.~Garsia and M.~Haiman,
  {\ita A remarkable $q,t$-{C}atalan sequence and $q$-{L}agrange
  inversion}, J. Algebraic Combin. {\bol 5} (1996), no.~3, 191--244.

\item{[10]}
A.~Garsia and M.~Haiman,
  {\ita Some natural bigraded ${S}_n$-modules and $q,t$-{K}ostka
  coefficients}, Electron. J. Combin. {\bol 3} (1996), no.~2, Research Paper 24,
  approx.\ 60 pp.\ (electronic),
  The Foata
  Festschrift,
  {\tt http://www.combinatorics.org/Volume\_3/volume3\_2.html\#R24}\ess.

\item{[11]}
A.~M. Garsia and C.~Procesi, {\ita On certain graded ${S}_n$-modules and the
  $q$-{K}ostka polynomials}, Adv. Math. {\bol 94} (1992), no.~1, 82--138.

\item{[12]}
A.~M. Garsia and J.~Remmel, {\ita {P}lethystic {F}ormulas and positivity for
  $q,t$-{K}ostka {C}oefficients}, {M}athematical {E}ssays in {H}onor of
  {G}ian-{C}arlo {R}ota (B.~E. Sagan and R.~Stanley, eds.), Progress in
  Mathematics, vol. 161, 1998.

\item{[13]}
A.~M. Garsia and G.~Tesler, {\ita Plethystic formulas for {M}acdonald
  $q,t$-{K}ostka coefficients}, Adv. Math. {\bol 123} (1996), no.~2, 144--222.

\item{[14]}
M.~Haiman, {\ita {M}acdonald {P}olynomials and {H}ilbert {S}chemes}, preprint.

\item{[15]}
A.~Kirillov and M.~Noumi, {\ita {R}aising operators for {M}acdonald
  polynomials}, preprint.

\item{[16]}
F.~Knop, {\ita Integrality of two variable {K}ostka functions}, J.\ Reine
  Angew.\ Math. {\bol 482} (1997), 177--189.

\item{[17]}
F.~Knop, {\ita Symmetric and non-symmetric quantum {C}apelli polynomials},
  Comment. Math. Helv. {\bol 72} (1997), no.~1, 84--100.

\item{[18]}
L.~Lapointe and L.~Vinet, {\ita Rodrigues formulas for the {M}acdonald
  polynomials}, Adv. Math. {\bol 130} (1997), no.~2, 261--279.

\item{[19]}
I.~G. Macdonald, {\ita A new class of symmetric functions}, S\'eminaire
  Lotharingien de Combinatoire, Publ. Inst. Rech. Math. Av., vol. 372, Univ.
  Louis Pasteur, Strasbourg, 1988, pp.~131--171.

\item{[20]}
I.~G. Macdonald,
  {\ita Symmetric functions and {H}all polynomials}, second ed., Oxford
  Mathematical Monographs, The Clarendon Press Oxford University Press, New
  York, 1995, With contributions by A. Zelevinsky, Oxford Science Publications.

\item{[21]}
I.~G. Macdonald,
  {\ita Affine {H}ecke algebras and orthogonal polynomials}, Ast\'erisque
  (1996), no.~237, Exp.\ No.\ 797, 4, 189--207, S\'eminaire Bourbaki, Vol.\
  1994/95.

\item{[22]}
E.~Reiner, {\ita A proof of the $n!$ conjecture for generalized hooks}, J.
  Combin. Theory Ser. A {\bol 75} (1996), no.~1, 1--22.

\item{[23]}
S.~Sahi, {\ita Interpolation, integrality, and a generalization of {M}acdonald's
  polynomials}, Internat. Math. Res. Notices (1996), no.~10, 457--471.

\item{[24]}
H.~Weyl, {\ita The {C}lassical {G}roups. {T}heir {I}nvariants and
  {R}epresentations}, Princeton University Press, Princeton, N.J., 1946.

}
\end